# The Topological Weighted Centroid (TWC):
# A topological approach to the time-space structure of epidemic and pseudo-epidemic processes

Massimo Buscema

*Semeion Research Center, Rome, Italy*

*Dept of Mathematical and Statistical Sciences, University of Colorado Denver, USA*

Giulia Massini

*Semeion Research Center, Rome, Italy*

Pier Luigi Sacco

*IULM University Milan, Italy*

*Harvard University, Cambridge MA USA*

*metaLAB (at) Harvard, Cambridge MA USA*

**Abstract**

This paper offers the first systematic presentation of the topological approach to the analysis of epidemic and pseudo-epidemic spatial processes. We introduce the basic concepts and proofs, at test the approach on a diverse collection of case studies of historically documented epidemic and pseudo-epidemic processes. The approach is found to consistently provide reliable estimates of the structural features of epidemic processes, and to provide useful analytical insights and interpretations of fragmentary pseudo-epidemic processes. Although this analysis has to be regarded as preliminary, we find that the approach's basic tenets are strongly corroborated by this first test and warrant future research in this vein.

Keywords

Topological Weighted Centroid; Alpha Point; Alpha, Beta, Gamma, Theta scalar fields; G IN-OUT; Meta-distance.

# The Topological Weighted Centroid (TWC):

# A topological approach to the time-space structure of epidemic and pseudo-epidemic processes

**1. Introduction**

Understanding the spatial and temporal structure of epidemic processes is a theme of major concern in scientific as well as policymaking terms. Properly tackling such a formidable challenge calls for a comprehensive, multidisciplinary approach [1], and for state-of-the-art mathematical tools [2]. Being able to forecast the spatial diffusion patterns of an epidemic and its time dynamics is quintessential for effective action, and mistakes are generally very costly both financially and in terms of loss of living beings. Classical ways to approach these issues have an algebraic base, namely, making use of the already available data to calibrate a set of differential equations to provide the best possible prediction of the epidemic's further unfolding [3]. We can speak therefore of an *algebraic approach* to the mathematical study of epidemics.

The key to this modelling strategy is capturing the structure of the interactions among the main variables involved in the process, and to set free parameters accordingly in order to reproduce the observable dynamic pattern as accurately as possible as the best possible guarantee of future predictive accuracy. For instance, several mathematical models have been proposed in the literature [4] to reconstruct the main features of the inter-pandemic ecology of influenza A in humans. A typical route is modelling the interactions among individuals being (or having been) infected by different viral strains. In this vein, various Susceptible-Infected-Recovered (SIR) models can be abridged through appropriate cross-immunity parameters [5-6].

Despite its intuitive appeal and it wide adoption, this approach has important limitations:

   a. The interaction among the variables needs to be modeled a priori, on the basis of previous data generated by the same process or by affine ones;
   b. Large sets of temporal and spatial data are needed, in order to effectively calibrate the model parameters;

c. The number of the free parameters needed to attain a reasonable data fit could be relatively large;
d. Validation of the model is generally an interpolation based on data from the same epidemic, or an extrapolation based data from similar epidemics.
e. Even when the dynamical modelling provides useful cause-effect relations, calibrating the main parameters and estimating connection strengths between variables are difficult tasks.

A way to amend these limitations is to supplement the algebraic approach with topological considerations, taking the observed spatial properties of the diffusion dynamics into account when tuning the free parameters [7-9].

A somehow complementary approach has been developed in forensic research, and has become known as 'geographical profiling' [10-15]. Albeit tailored to different kind of problems and policy concerns, geographical profiling has been also applied to detect the outbreak points of many epidemics [16], as well as to analyzing other biologically motivated problems such as the spatial patterns of animal foraging [17]. We can label this second stream of literature as a *geometric approach*. Here, researchers try to determine the outbreak area of an epidemic and its likely patterns of spatial diffusion making use of only one snapshot of the past diffusion history of the epidemics. The advantage of this approach over the previous one is apparent: unlike the former, it needs relatively few data and, once more data become available, they can be quickly and easily plugged in to update estimations.

However, this approach as well has substantial limitations:

a. Its validation protocol is statistically weak.
b. The estimated outbreak area depends on the spatial concentration of the event locations. If observed locations are relatively concentrated, the probability of finding the outbreak increases. But when the spatial pattern is more scattered, the method could arrive at ill-founded conclusions.
c. Consequently, the outbreak area estimated by this method generally lies within the convex hull of the distribution of the event locations, and this is a major methodological shortcoming which seriously limits its applicability.

In the last five years, a new approach has appeared, which we could call a *topological approach*[1]. It shares many characteristics of the geometrical approach, and some of the algebraic approach as well. However, the new approach makes it possible to estimate the location of both the outbreak point and the dynamic unfolding of the epidemic by working only on the spatial distribution of the events, without any reference to the chronology and frequency of the observations, that is, on the basis of the purely spatial characteristics of the phenomenon. Moreover, unlike the two previous approaches, the topological one works on the basis of an optimization principle, that is, it adaptively explores the fitness landscape with respect of certain quantities (free energy and entropy), and determines the optimal solution accordingly. The aim of this paper is to present the theoretical foundations of the topological approach in a systematic manner for the first time in the literature, and to illustrate its potential through its application to a few real case studies.

---

[1] The topological approach was developed by M Buscema at Semeion Research Center (Rome, Italy), from 2008 to the present [18-21, 35-37].

Since the topological approach does not need to take into account the actual sequence and frequency of observations to derive the structural features of the underlying process and to predict its future unfolding, it is possible to apply it not only to other kinds of diffusion phenomena, as it was the case for the geographical profiling approach, but also to the analysis of cases where the spatial distribution of the events was not the outcome of a diffusion process in the proper sense of the word. For instance, we can apply the method to characterize the structural features of the distribution of geographical location of points that characterize a certain type of human activity (urban settlements, social and economic activities, and so on) which were not actually shaped by a diffusion process of an epidemic nature, or more generally by the action of a well-recognizable agent, but which nevertheless yield an observed spatial pattern that is observationally indistinguishable from one generated by an epidemic or an agent-driven one. We speak in this case of a pseudo-diffusion pattern or, in case we focus upon epidemic processes as the conceptual benchmark, of a pseudo-epidemic one. As it will be shown in this paper, regarding spatial distributions of human activities as pseudo-epidemic patterns allows us to gain considerable insight into their structural characteristics, as well as to effectively forecast their future unfolding. The topological approach, therefore, has a vast range of application to diverse classes of problems in operation research where the spatial distribution dimension has a key relevance.

## 2. Basic theory

*2.1 The structure of the dataset*

The basic intuition behind the topological approach is that every distribution of point in space has an implicit semantics, provided the following conditions are met:

a.  Each point of the distribution represents a discrete occurrence of the same process;

b.  The distribution of points is statistically representative of the process to be analyzed;

Consequently, we can represent a generic dataset in a two-dimensional space, which is amenable to analysis through the topological approach, as follows:

$$(1) \quad Data = \{P(x_i, y_i)\}_{i=1}^{N};$$
$$where:$$
$$x = Longitude;$$
$$y = Latidute;$$
$$N = Number\ of\ Data\text{-}Points\ (usually\ N > 3);$$
$$i \in \{1, 2, ..., N\}.$$

If the data-points represent a time series, we use the alternative formulation

$$(1a) \quad Data = \{P(x_t, y_t)\}_{t=1}^{T};$$

*where*:

$x = Longitude$;

$y = Latidute$;

$T = $ *Number of Temporal Data - Points;*

$t \in \{1, 2, ..., T\}$.

When to each data-point is attached a vector of values of characteristic attributes, the formulation has to be revised accordingly as

$$(1b) \quad Data = \left\{P\left(x_i, y_i, \{a_{i,k}\}_{k=1}^{M}\right)\right\}_{i=1}^{N};$$

*where*:

$a = $ *specific attribute;*

$M = $ *Number of Attributes for each data point;*

$k \in \{1, 2, ..., M\}$.

In the topological approach, it is assumed that there are several different types of information that are related to the spatial semantics of the dataset:

a. Explicit information, that can be extracted, with different levels of accuracy, by means of a vast range of statistical techniques;
b. Misleading information, that is associated to the noise contained in the dataset, and which can be partially filtered off, with different levels of accuracy, by a vast range of statistical techniques, but whose categorization and analysis can be very different in terms of conceptual depth and accuracy across different techniques;
c. Implicit information, that is hidden in the relations among the data-points, and which is generally difficult to extract and analyze by means of most of the available statistical techniques.

The main goal of the topological approach is the mining of the implicit information embedded into the dataset. For a more thorough mathematical treatment of this methodology see [18-25].

### 2.2 The Alpha Point

The first key concept we need to introduce to shape up our topological approach is the notion of an Alpha Point. The Alpha Point represents a spatial estimate of the hidden (outbreak) point, or area, where the process under study originated. The basic assumption behind this notion, as it is characteristic of the topological approach, is that the geometry of the spatial distribution of the observed events implicitly codes some sort of temporal information, as a consequence of the fact that the observed events have materialized at different points in time. In other approaches, it is implicitly assumed that the regions where a particularly dense occurrence of events is observed is more likely to be the one where the process has initially originated from, i.e. the region that contains the outbreak point or area. The topological approach follows a different line of reasoning: The ideal candidate outbreak point/area is not in the region with the highest concentration of events, but in the one that presents the highest level of organization of their spatial distribution. A highly organized region represents the portion of space from where the information that is needed to code and retrieve the relative position of all the other events is optimized, that is, the location from which the global entropy of the distances from all the other points attains its minimum level. This has an intuitive appeal in terms of both animal and human behavior. For instance, both predators and preys tend to choose their temporary base at locations from where it is easier to maintain a fully sensory control of the surrounding environment, that is, from where they can perceive the environment, and its critical locations, in the most coherent and efficient way in terms of attention resources. Likewise, the optimal location of the military command unit during a battle is generally found at a position where the entire battlefield can be scanned in the most direct and organized way, to track in real time the development of the events and to make rapid, informed decisions.

In formal terms, we can proceed as follows. Regarding all the events in space as points in a scalar field, and, in particular, as local energy minima, we can calculate the attraction strength, $w_i$, of each point as

$$(2) \quad \overline{d_{i,j}} = \frac{1}{N-2} \sum_{\substack{k \neq j \\ k \neq i}}^{N} d_{i,k}$$

$$(3) \quad w_i = \frac{1}{N-1} \sum_{j=1, j \neq i}^{N} e^{-\frac{\overline{d_{i,j}}}{D}};$$

where:

$$d_{i,k} = \sqrt{\sum_{z}^{D} \left(x_{i,z} - x_{k,z}\right)^2} \quad \textit{distance between two points of the dataset;}$$

$$D = \max_{i,j}\{\overline{d_{i,j}}\}.$$

Equation (2) defines 'distance' between two points $i$, $j$ in terms of the average of the distance of one point $i$ from all the others *but* the other one $j$. We call this quantity 'indirect distance'. The more a point is relatively close to most others, therefore, the less it contributes to the distance. A point that is relatively close to most others may therefore be 'close' to another one in terms of indirect distance even if the two are physically distant between each other in the original space. Vice versa, a point that is physically close to another, but relatively distant from most others, will be 'distant' from the other despite being physically close to it in the original space. Therefore, application of the indirect distance 'expands' or 'contracts' the space between two points depending on their relative position in the global distribution. This definition reflects our focus upon the organizational properties of the positions of points in the whole spatial distribution. The indirect distance between two points is therefore a holistic, and not individual property of the points, and consequently indirect distance is not a proper distance, in that it does not satisfy the basic conditions for a metric: The indirect distance between a point and itself need not be null (and generally will not be); the indirect distance between *i* and *j* is generally different from that between *j* and *i*; the triangular inequality needs not be satisfied. Equation (3) quantifies the attraction strength of each given point with respect to all others, where strengths exponentially decay with indirect distance. The stronger the attraction strength at one point, the more such point 'communicates' with the other, that is, the more its location becomes 'meaningful' once put in relation with that of all the other points.

Despite not being a proper distance, indirect distance possesses interesting features from the viewpoint of our topological approach:

   a. As already remarked, indirect distance between two points depends on the distances of each point from all the others (holistic property);
   b. Although indirect distance basically alters the notions of closeness and distance in the Euclidean sense, the sum of all the indirect distances across all points equals that of the classic Euclidean distances between the same points (conservation property):

(i) $\overline{d_{i,j}} = \frac{1}{N-2} \sum_{\substack{k \neq j \\ k \neq i}}^{N} d_{i,k}$; //Indirect Distance

(ii) $d_{i,j} = \sqrt{(x_i - x_j)^2 + (y_i - y_j)^2}$; //Euclidean Distance

(iii) $\Delta d_{i,j} = d_{i,j} - \overline{d_{i,j}}$;

(iv) $Row_i = \sum_{j}^{N} \Delta d_{i,j} \approx 0$;

(v) $Col_j = \sum_{i}^{N} \Delta d_{i,j}$;

(vi) $\sum_{k}^{N} Row_k = \sum_{k}^{N} Col_k \approx 0$.

   c. The recursive calculation of the indirect distance (the indirect meta-distance) presents a global point-wise attractor (convergence property).

Let's present a simple example to fix ideas. Consider the space distribution of 10 points depicted in Figure 1.

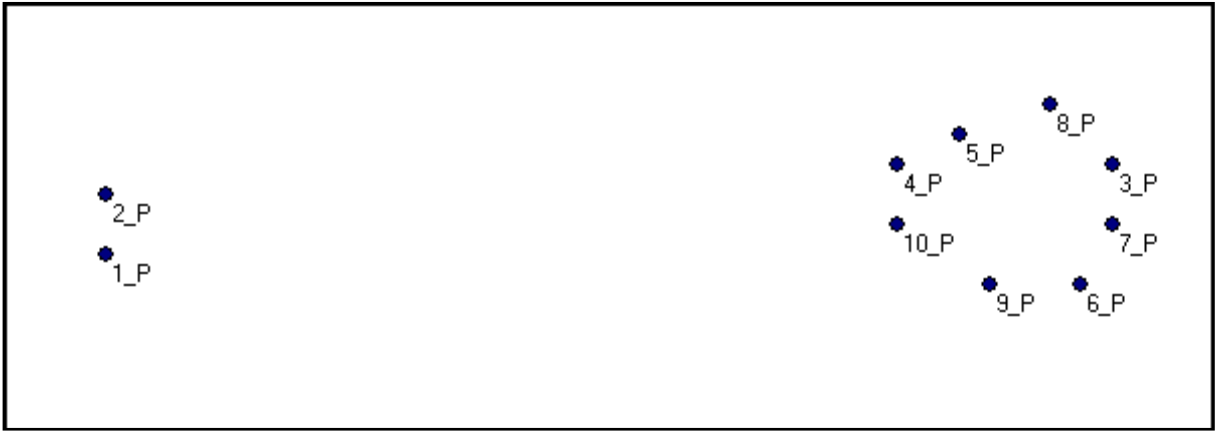

*Figure 1: Ten points with two outliers*

| Deltas Matrix | 1_P | 2_P | 3_P | 4_P | 5_P | 6_P | 7_P | 8_P | 9_P | 10_P | Rows Sum |
|---|---|---|---|---|---|---|---|---|---|---|---|
| 1_P | 0 | -422.463509 | 108.265198 | -10.432159 | 25.563843 | 89.166443 | 106.203796 | 78.683716 | 38.057556 | -13.044922 | 0.0000 |
| 2_P | -421.553017 | 0 | 107.114288 | -12.13443 | 22.844818 | 92.202423 | 107.114288 | 75.233734 | 41.312195 | -12.13443 | -0.0001 |
| 3_P | 375.162865 | 373.101494 | 0 | -70.338531 | -102.741631 | -119.376522 | -155.565819 | -141.444885 | -93.233055 | -65.563942 | 0.0000 |
| 4_P | 278.879486 | 276.266731 | -47.924561 | 0 | -129.127937 | -44.32621 | -43.149956 | -75.450165 | -82.01548 | -133.151833 | 0.0001 |
| 5_P | 315.846245 | 312.216705 | -79.356895 | -128.157204 | 0 | -57.127823 | -66.880775 | -112.369556 | -79.356911 | -104.81374 | 0.0000 |
| 6_P | 357.204445 | 359.329918 | -118.236217 | -65.599876 | -79.372238 | 0 | -150.401646 | -84.832962 | -137.380005 | -80.711464 | 0.0000 |
| 7_P | 374.648933 | 374.648933 | -154.018379 | -64.016502 | -88.71804 | -149.994495 | 0 | -111.879669 | -111.879662 | -68.791092 | 0.0000 |
| 8_P | 351.139809 | 346.779335 | -135.886475 | -92.30574 | -130.195866 | -80.414856 | -107.868713 | 0 | -76.293396 | -74.954147 | 0.0000 |
| 9_P | 321.638275 | 323.982414 | -76.550026 | -87.746422 | -86.058586 | -121.837265 | -96.744072 | -65.168762 | 0 | -111.515442 | 0.0001 |
| 10_P | 276.355431 | 276.355431 | -43.061272 | -133.063133 | -105.695789 | -59.349113 | -47.835876 | -58.009903 | -105.695831 | 0 | -0.0001 |
| Cols Sum | 2229.32247 | 2220.21745 | -439.654339 | -663.793997 | -673.501426 | -451.057418 | -455.128773 | -495.238452 | -606.484589 | -664.681012 | -0.0001 |

*Table 1: Each cell reports the difference between the Euclidean and Indirect Distances. Negative signs stand for space 'expansion', positive signs for space 'contraction'.*

Table 1 shows the deformation of space induced by the contractions and expansions that take place once the Euclidean distances are substituted by the indirect distances. Table 2 reports the different attraction strengths of the 10 points, as computed according to equation (3), whereas Figure 2 shows the 3d representation of the scalar field interpolated by the attraction strength of the ten points.

| Points | Points Strenght |
|---|---|
| 5_P | 0.967892 |
| 10_P | 0.967723 |
| 4_P | 0.967706 |
| 9_P | 0.96662 |
| 8_P | 0.964513 |
| 7_P | 0.963755 |
| 6_P | 0.963677 |
| 3_P | 0.963462 |
| 2_P | 0.914437 |
| 1_P | 0.914274 |

*Table 2: Strength of attraction of the ten points*

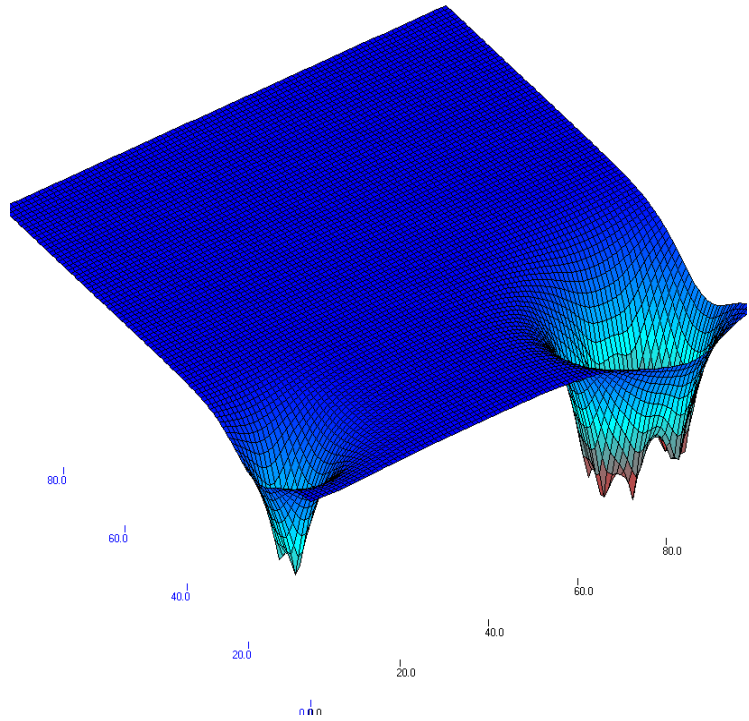

*Figure 2: The 3d representation of the attraction strength of the 10 points.*

Figure 3 shows the dynamic of the recursive application of indirect distance algorithm, where the values of indirect distances gradually collapse, converging to a single point.

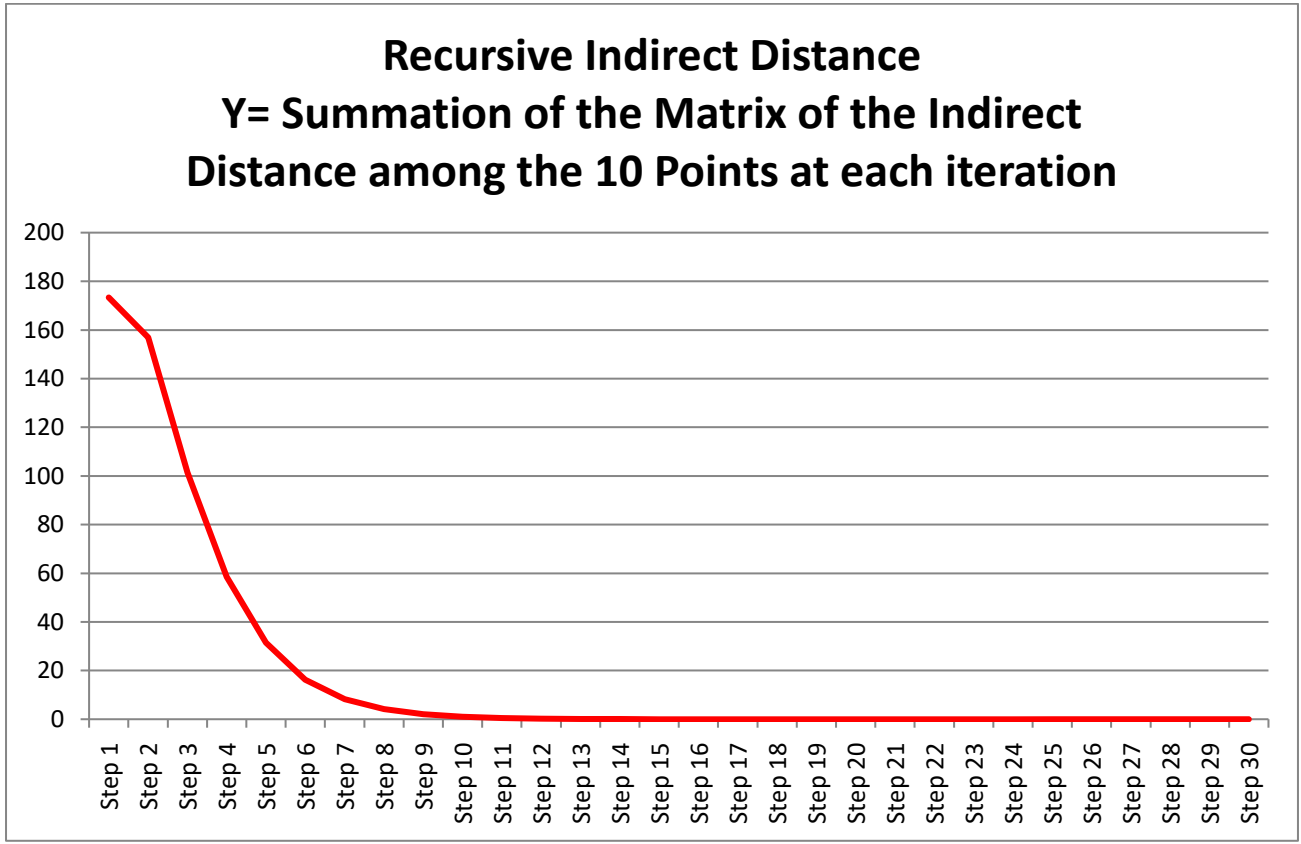

*Figure 3: The dynamic of recursive application of indirect distance algorithm on the 10 points of Figure 1.*

We can think of recursive indirect distance as a measure of the energy needed for a given distribution of points in space to collapse their spatial organization into a single point. The deeper the level of recursion needed for the collapse to take place, the higher the level of energy required to the purpose in view of the spatial structure of the distribution.

Let us now introduce a control parameter $\alpha_n$, where $\alpha_n \in [0,\infty]$, which modulates the attraction strength of points. We may then rewrite equation (3) as:

$$(3a) \quad w_i(\alpha_n) = \frac{1}{N-1} \sum_{j=1,\ j \neq i}^{N} e^{-\frac{\overline{d_{i,j}}}{D} \alpha_n}.$$

Think of the α parameter as the inverse of the temperature of a thermodynamic system, that is, define $T \equiv 1/\alpha$. When α is small, the average temperature of the system is high, and vice versa when α is large, the temperature of the system is close to zero. This analogy will be convenient for the intuitive interpretation of the quantities that we are going to introduce below.

Once we have reformulated our space topology as defined by the contraction/expansion of original space operated by indirect distances, we are now in the position to introduce, in analogy to what is done in the geometric approach, our own definition of a 'centroid' that represents our estimate of the outbreak point of the spatial distribution. The major difference is that our space averaging is not conducted on the original space, but rather in the modified space that takes into account the holistic structure of the distribution and in particular its varying local levels of organization, which represent for us a key information to infer the temporal unfolding of the process. Moreover, the averaging weights reflect the attraction strengths of the various points. We then define the Topological Weighted Centroid (TWC) of all the points of the dataset for any specific value of the free parameter α as:

$$(4a) \quad TWC_x(\alpha_n) = \frac{1}{\sum_{i=1}^{N} w_i(\alpha_n)} \sum_{i=1}^{N} w_i(\alpha_n) \cdot x_i;$$

$$(4b) \quad TWC_y(\alpha_n) = \frac{1}{\sum_{i=1}^{N} w_i(\alpha_n)} \sum_{i=1}^{N} w_i(\alpha_n) \cdot y_i;$$

When $\alpha_{n=0} = 0$, then the TWC coincides with the center of mass of the space distribution, that is, it boils down to the average of the coordinates of the points. Consequently, when all the weights are equal to 1, equations (4a and 4b) simply determine the center of mass. Monotonically increasing the control parameter by small quantities ($\alpha_{n+1} = \alpha_n + \varepsilon.$), the TWC moves away from the center of mass, and for each iteration, we may calculate the entropy, $H$, of the attraction strength of all the points, by means of equations (5-6):

$$(5) \quad p_i(\alpha_n) = \frac{w_i(\alpha_n)}{\sum_{k}^{N} w_k(\alpha_n)};$$

$$(6) \quad H(\alpha_n) = -\sum_{i}^{N} p_i(\alpha_n) \cdot \log_2(p_i(\alpha_n)).$$

Equation (5) yields the relative attraction strength of each point in relation to the total attraction strength of the distribution, whereas (6) computes the entropy $H$ for any specific value of $\alpha_n$. The value of the entropy at the center of mass is $H(\alpha_{n=0}) = \log(N)$, which corresponds to the maximum level of entropy. Moving away from the center of mass as $\alpha_n$ varies, the entropy decreases monotonically until it reaches a minimal value for $H(\alpha_{n=\infty})$, where the TWC coincides with the point where the total indirect distance $\overline{d_{i,j}}$ from all other points of the distribution is minimal (for details, see appendix B in [21]).

### 2.3 Free energy

In this section, we suitably adapt from physics our concept of 'free energy'. In physics, free energy (*F*) is a measure of the amount of energy that is available to perform work. Whenever a change in energy takes place, the total amount of energy remains constant, as with the increase of entropy, disordered energy increases and free energy decreases (for more details on the physics interpretation see Appendix C in [21]). We borrow this conceptualization as a natural basis for our mathematical modelling of the topological approach, without making use of these concepts in a literal sense, but as convenient references that facilitate by analogy the interpretation of our formalism and quantities.

If we denote by Z the partition function of a thermodynamic system at 'temperature' ('thermal energy') $1/\alpha$ ($T \equiv 1/\alpha$), using this partition function we can compute free energy ***F*** as in equation (7):

$$(7) \quad F(\alpha_n) = \frac{-\ln(Z)}{\alpha_n};$$

$$-\infty \leq F(\alpha_n) \leq 0.$$

*where*:

$$Z = \sum_{i}^{N} w_i(\alpha_n).$$

We can thus single out a 'high temperature' phase (for values of $\alpha_n$ such that *TWC($\alpha_n$)* is close to the center of mass), and a 'low temperature' phase (for values of $\alpha_n$ such that *TWC($\alpha_n$)* is close to the region with the maximum density of points), which are in turn separated by a region where the free energy is maximum at least in one point. This point, as defined by equation (8), can be proven to be unique.

$$(8) \quad F(\alpha^*) = \underset{\alpha}{Max}\{F(\alpha_n)\};$$

The maximum value of free energy, that is reached when the free parameter $\alpha_n$ equals $\alpha^*$, represents a singular point (see Appendix D in [21]: For larger values of $\alpha_n$, entropy and 'free energy' jointly decrease, as 'temperature' monotonically decreases. We can therefore characterize $\alpha^*$ as an optimal value for α, in that in $\alpha^*$, 'free energy' (that is, the potential attraction strength among the points of the distribution), is maximized. Figures 4 and 5 report the variation of 'free energy' and entropy for the example in Figure 1 as α increases.

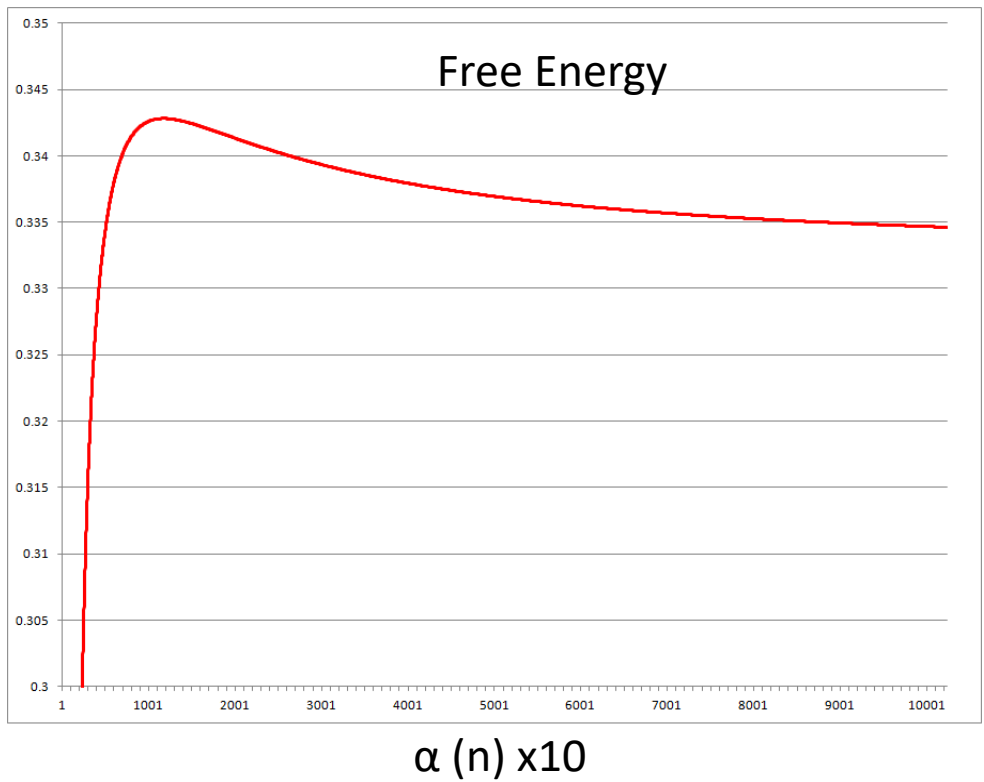

*Figure 4: Free Energy as function of α (example).*

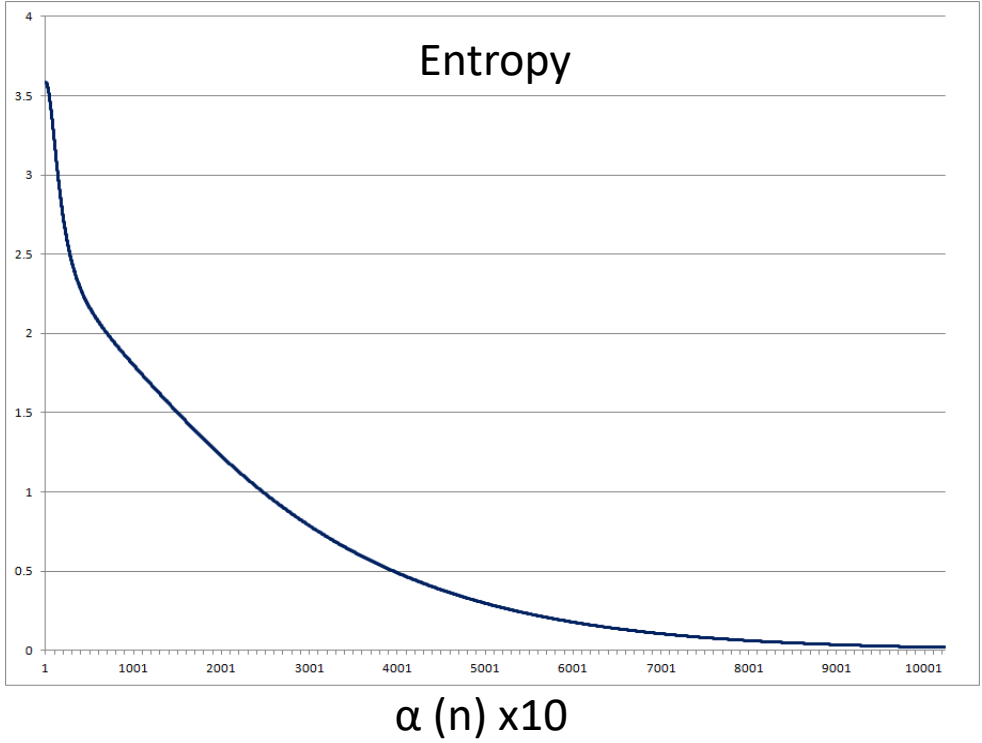

*Figure 5: Entropy as function of α (example).*

Once the optimal value of α, that is, $α^*$, has been determined, we can compute our best estimate of the outbreak point of our distribution in terms of the topological weighted centroid (TWC) at $α^*$. We call this point the Alpha Point, and we then reformulate equations (4a-4b) as:

$$(9a) \quad TWC_x(α^*) = \frac{1}{\sum_{i=1}^{N} w_i(α_n)} \sum_{i=1}^{N} w_i(α^*) \cdot x_i;$$

$$(9b) \quad TWC_y(α^*) = \frac{1}{\sum_{i=1}^{N} w_i(α_n)} \sum_{i=1}^{N} w_i(α^*) \cdot y_i;$$

Let us elaborate a bit about these quantities and the point in space they represent, which can be thought of as the 'hidden' point of the distribution for reasons that will become clear below. The value of α at which 'free energy' attains its maximum value corresponds to the α where the points of the distribution manifest the maximum level of attraction strength among them. The corresponding TWC, the Alpha Point, is the point in space from which the whole distribution of the points reaches its highest level of organization. Being the total attraction strength at maximum value, this is the TWC which, in the whole family of centroids defined by the variation of α, ensures the maximum level of communication among all points. Taking such the Alpha Point as a reference, the location of all other points is the most meaningful possible with respect to the given spatial distribution of the points. The Alpha Point is therefore the ideal vantage point to appreciate the spatial complexity of the distribution in all of its informational richness. If there is a 'pattern' in the distribution, this is the best location to look at it and to discern it. As α moves from its null value, with the corresponding TWC coinciding with the center of mass, toward the optimal value α*, the TWC moves away from the center of mass and describes a trajectory across the space. For a given sequence of values of α, $α_n$, we therefore have a vector of TWC($α_n$) points, which represents the discretization of the continuous trajectory described by TWC($α$) as $α$ varies along the real interval. As the sequence converges toward the Alpha Point TWC(α*), the relative distance between the TWC($α_n$) positions at two subsequent iterations decreases.

In our topological approach, we make use of this vector of subsequent positions of TWC($α_n$) to transform the space where the distribution of points sits into a scalar field, in the way that will be described below.

To sum up, we have so far introduced a few key basic concepts on which our topological approach is built:

a. The concept of indirect distance (Equation 2);
b. The concept of attraction strength that each point exerts on the others (Equation 3);
c. The concept of entropy (Equation 6);
d. The concept of 'free energy' (Equations 7-8);
e. The concept of Alpha Point (Equations 9a-9b);
f. And the concept of Vector of Alpha(n) Points (Equations 4a-4b).

*2.4 The Beta Parameter*

Once introduced the basic battery of concepts we can now work on some refinements of our topological approach. A first one is the Beta parameter, which builds on the same logic already introduced for the Alpha one. Specifically, we reformulate equation (3a) as follows:

$$(10) \quad w_i(\beta_n) = \frac{1}{N-1} \sum_{j=1}^{N} e^{-\frac{\overline{d_{i,j}}}{D} \beta_n}.$$

The main difference between (3a) and (10) is that in the latter, each point interacts also with itself, that is, we also consider the indirect distance between one point and itself which, as already discussed, is generally not null. As in the previous case, we can study how the attraction strength (10) varies along a sequence of small increases of β, that is, $\beta_{n+1} = \beta_n + \varepsilon$. Accordingly, for each value of β we can compute the position of the corresponding TWC in the space. In this case, for reasons that are now intuitive, we speak of Self-TWC (STWC), which is formally defined as:

$$(11a) \quad STWC_x(\beta_n) = \frac{1}{\sum_{i=1}^{N} w_i(\beta_n)} \sum_{i=1}^{N} w_i(\beta_n) \cdot x_i;$$

$$(11b) \quad STWC_y(\beta_n) = \frac{1}{\sum_{i=1}^{N} w_i(\beta_n)} \sum_{i=1}^{N} w_i(\beta_n) \cdot y_i;$$

The STWC has specific properties that differentiate it from the TWC introduced above. In particular, when β is small, as for the case of TWC the STWC moves away from the center of mass. However, when β reaches a specific value, STWC, unlike the TWC, heads back toward the center of mass. The reason is that, when the free parameter β becomes large enough, the interaction of each point with itself takes over and the attraction strengths tend to collapse into self-attraction, that is, focus more and more, for each point, upon its own position rather than on its relative position within the global distribution of points across the space. This specific value of beta beyond which the trajectory switches back, $\beta^*$, is the optimal β that represents the logical analogue of $\alpha^*$. $\beta^*$ is the boundary value beyond which each point starts to decrease its strength of interaction with the others. Consequently, β* represents the value of β at which the bell curve that describes the radius of interaction of each point with the others reaches its maximum width:

$$(12) \quad \beta^* = \underset{\beta_n}{Max}\left\{D\left(STWC_{\beta_n}, MassCentre\right)\right\};$$

$D$ = Euclidean Distance.

An example of a parametric family of bell curves as β varies for the example in Figure 1 above is shown in Figure 4 below.

The interpretation of the STWC is basically different than that of the TWC. Here, the point that is selected is the one that optimally balances the self-referential perspective of the single points and the global perspective of the whole distribution of points: the local information contained in any specific position and the non-local information contained in the whole distribution. Each point is itself a potential candidate to serve as the centroid of the distribution, that is, each point allows to look at the distribution from a perspective that carries a specific value and meaning of its own. At the STWC, the balance between local and global information reaches an optimal tradeoff, and is therefore telling in terms of how, given a certain distribution of points, each already existing point (that is, each event that has already occurred) is going to play a role in the future unfolding of the phenomenon once the local/global tradeoff is further influenced by the emergence of new points (the occurrence of further events).

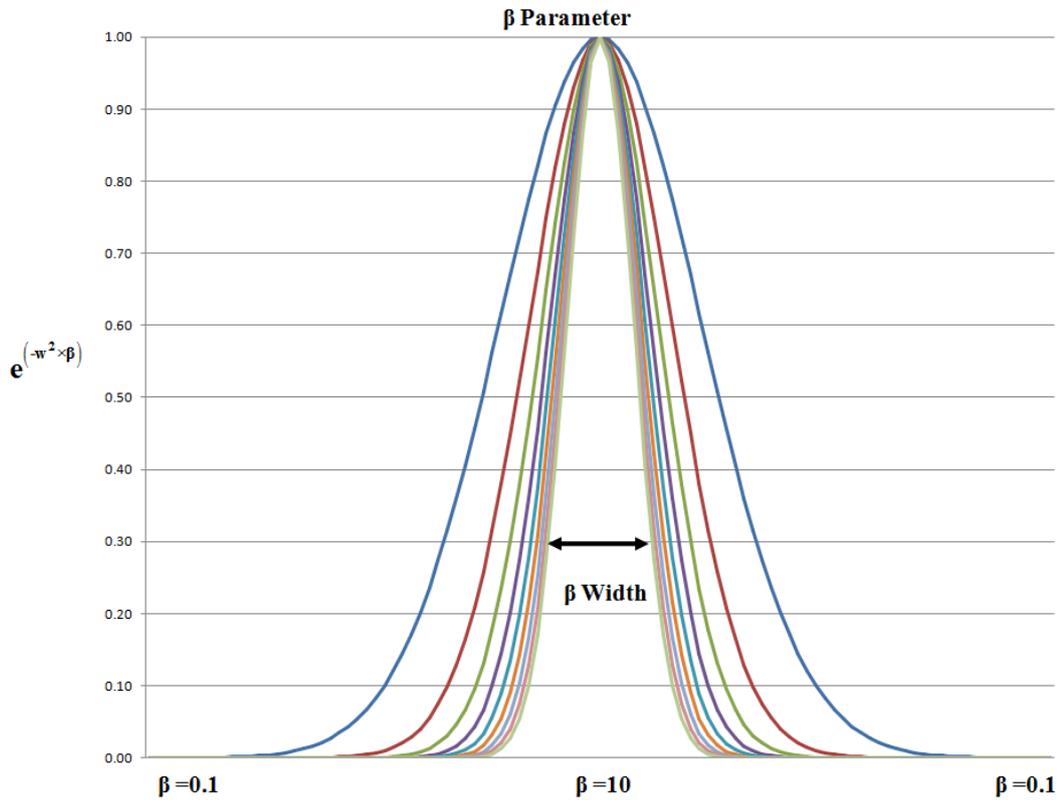

*Figure 4: Width of bell curves as β varies. **w** is the strength of each node ( $-10 \leq w \leq +10$ ).*

## 3. Scalar fields

Now that we have defined the notions of the TWC and of the STWC, we are ready to explore how the discretization of their trajectories drawn out for given sequences of values of the corresponding free parameters can be used to generate scalar fields with very interesting properties and implications.

### 3.1 The Alpha scalar field: An application to a pseudo-epidemic process

Making use of the Vector of Alpha Points (TWC(α(n))) and of the (optimal value of the) Beta parameter (β*), our topological approach allows us to transform the assigned two-dimensional data set of *N* points/events into two different scalar fields:

a. The Alpha Map
b. The Beta Map.

The Alpha Map quantifies the closeness of each point of the space to the parametric family of TWCs (and thus, to the Alpha Point), that is, the estimated outbreak from which the process has spatially originated. The more a given point is close to the trajectory defined by the Vector of Alpha Points, the more that point has a strong activation value – that is the more it is structurally connected to the deployment of the underlying phenomenon causing the observed distribution of events.

To fix ideas, we consider an example of a real pseudo-epidemic process for which we have specific historical evidence, namely, the spatial organization of Etruscan towns during the Roman Iron Age. Other applications of TWC theory in the archaeological field can be found in [39] and [40]. Figure 5 reports the position of the twelve main Etruscan towns during the Roman Iron Age. In this case, the Alpha Point may be regarded as the most meaningful location to observe the spatial pattern drawn out by the observed settlements, that is, a place of special significance for the socio-spatial organization defined by the distribution of the twelve cities. Figure 6 shows the Alpha Point for this spatial distribution, as well as the Vector of Alpha Points (TWC(α(n))). Figure 7 shows the corresponding Alpha Map, whereas Figures 8a-8b show the projection of the Alpha Map on the real geographical map. The Alpha Point is located between the site of Monte Becco and the coast of Bolsena Lake, in an area that, according to local popular tradition, is the location of the Fanum Voltumnae, the temple which hosted the yearly meetings of the federation of the twelve Etruscan cities [22], and which therefore constitutes the natural hub of the system of local territorial relations. According to state of the art research, the actual location of the Fanum is still controversial [23], but it is significant that the Alpha Point identifies with impressive accuracy the most credited location according to orally transmitted local knowledge.

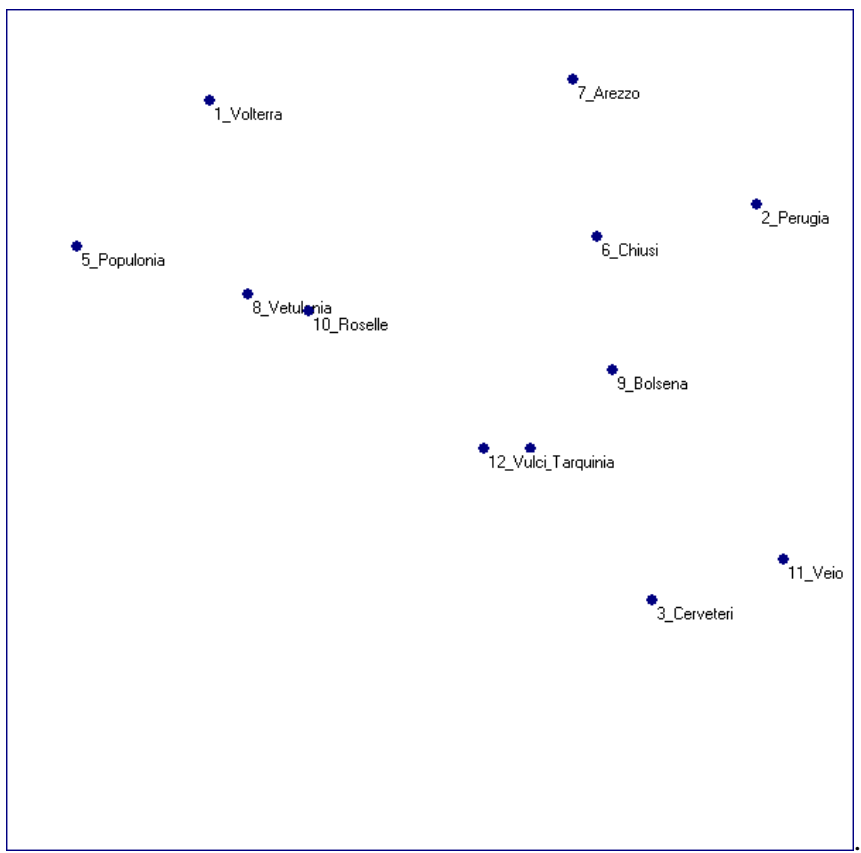

*Figure 5: Location of the twelve Etruscan cities during the Roman Iron Age.*

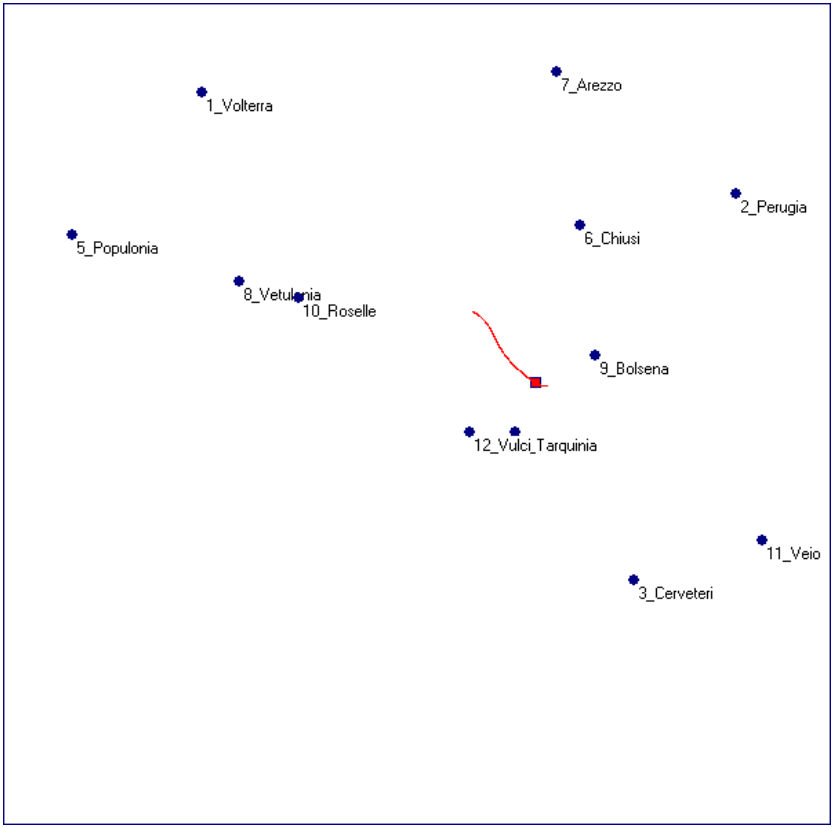

*Figure 6: Alpha Point (red square) and Vector of Alpha Points (red line).*

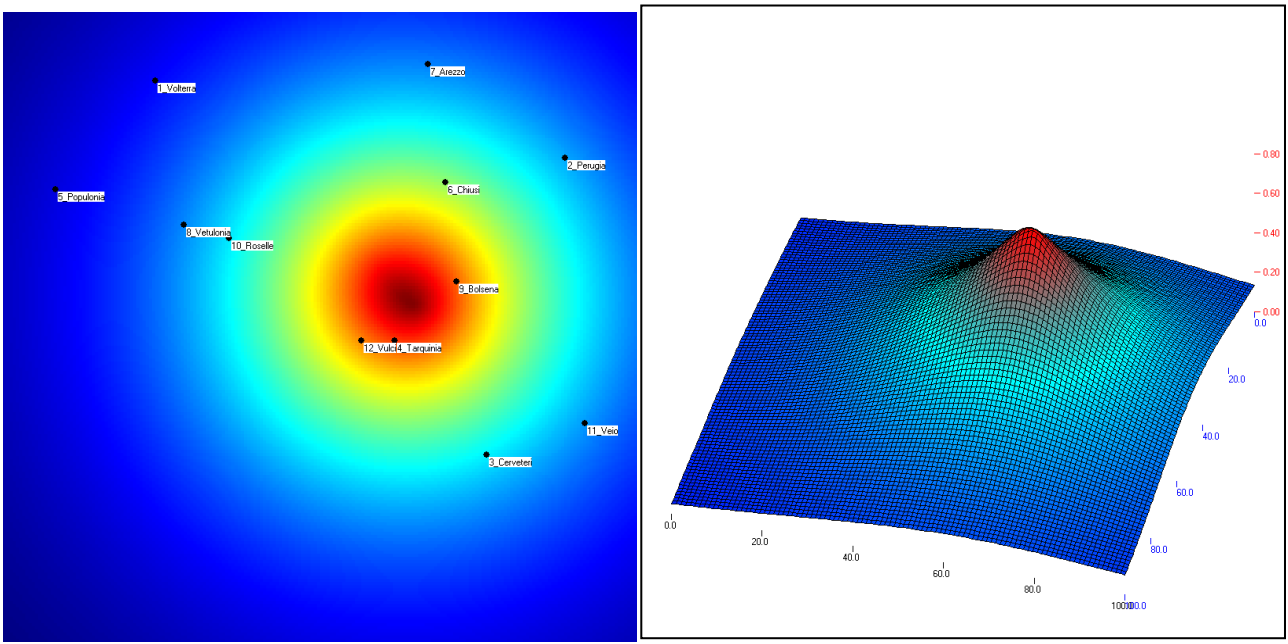

*Figure 7: Alpha Map (2D and 3D) (red color corresponds to higher activation).*

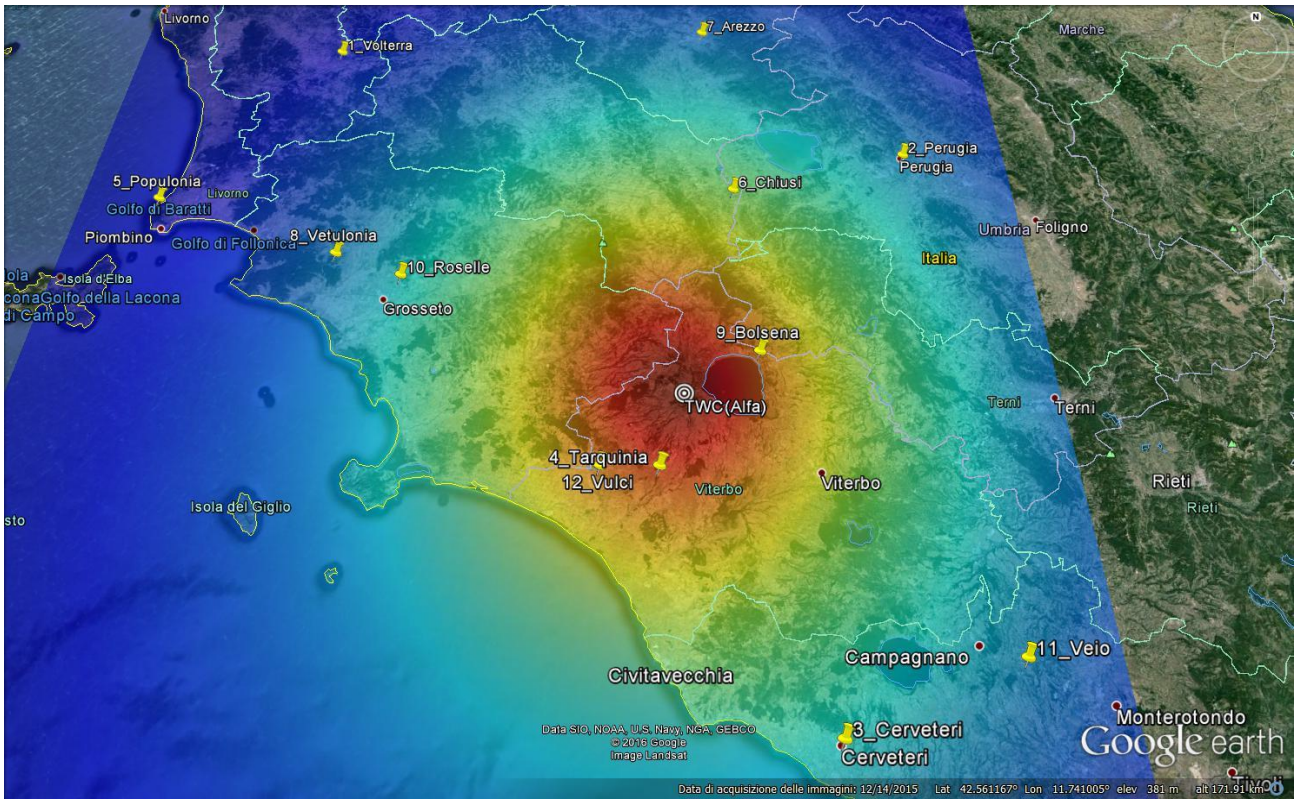

*Figure 8a: Color map overlays of the Alpha Map using Google earth mapping.*

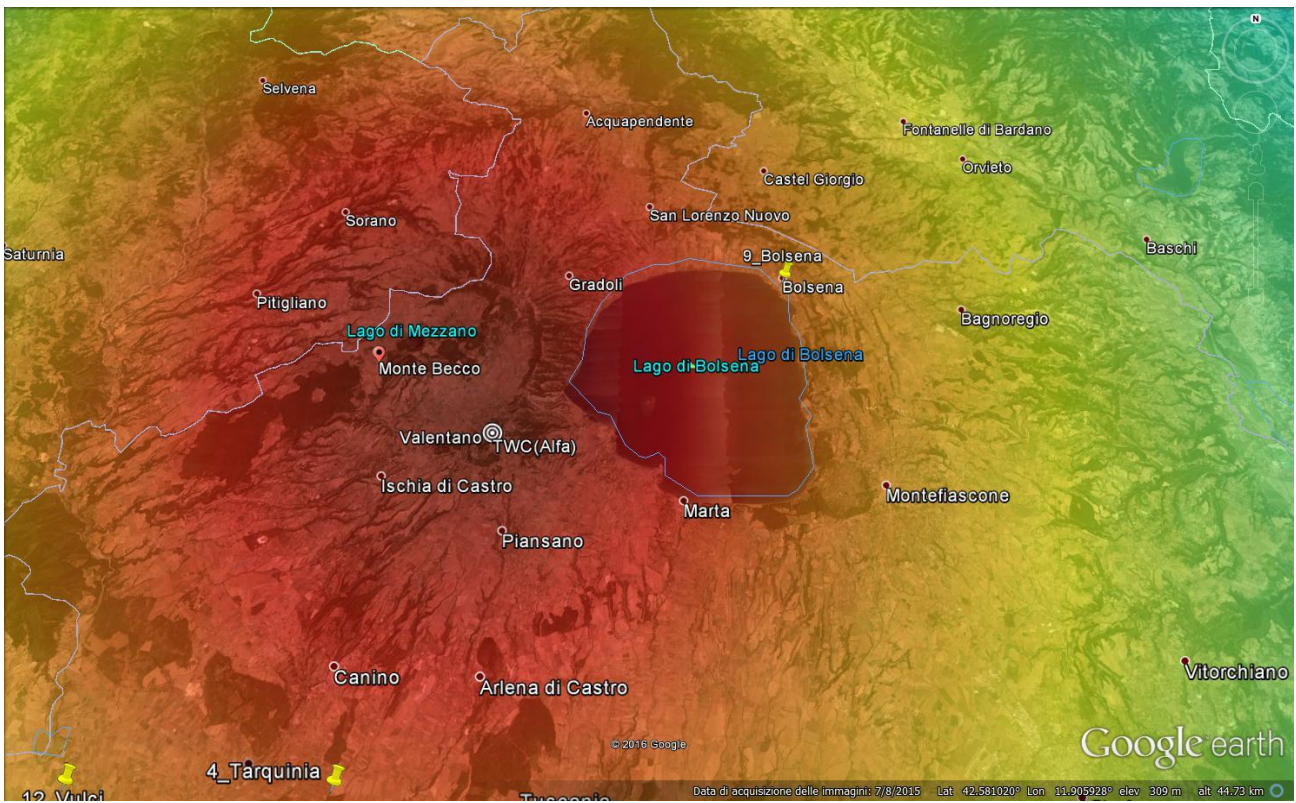

*Figure 8b: Zoom of the Alpha Map overlay based upon Google earth mapping. The Alpha Point is considerably close to Monte Becco and the coast of Bolsena Lake.*

The Alpha Map measures the degree of 'activation' of different regions of the space as to their proximity (in terms of indirect distance) to the path leading to the Alpha Point. The more active a region, the more likely it is to play a role in the further unfolding of the process. The Alpha Map therefore represents a projection of the likely patterns of further spatial diffusion of the process on the basis of its observed, current spatial 'grammar'. It is important to stress that such an inference may only be correct if the observed distribution of points represents a statistically significant sample of the underlying phenomenon, in the absence of spurious observations.

The Alpha Map therefore estimates the short-term unfolding of the process by means of the observed distribution of past events. The estimation process is based upon the Vector of Alpha Points, as follows:

$\{$Euclidean distance, $\{m_{k,j}\}$, between each point of 2D space and each Twc Alfa Point$\}$

(13) $m_{k,j} = \sqrt{\left(x_k - twc_x(\alpha_j)\right)^2 + \left(y_k - twc_y(\alpha_j)\right)^2}$;

$j =$ index for each of the TWC Alpha Points;
$k =$ index for each point of the plane.

$\{$Activaction strength of each point of the 2D plane$\}$

(14) $a_k = \dfrac{1}{V} \sum_{j=1}^{V} e^{-\dfrac{m_{k,j}}{D}\beta^*}$;

$V =$ the total number of TWC points.
$D = \max_{i,j}\{d_{i,j}\}$; Maximum distance among the points of the assigned dataset;
$\beta^* =$ optimal value of $\beta$ according to equation (12);
$a_k =$ Value of activation of $k$-$th$ point of the plane; $a_k \in [0,1]$.

The activation value $a_k$ is an element of the scalar field that describes the responsiveness of that specific location to the diffusion process under way. If we consider a specific portion of the space (such as the box window of Figure 7) we can interpret such activation values as probabilities of occurrence of process-related events. The more circumscribed the high activation area (the red one), the more precise the probabilistic estimation. Consequently, we can compute a p-value for each area of the box window, as shown in the Alpha Map for the Etruscan cities example (see Table 2 and Figure 9). When a new event occurs, the lower the p-value of the corresponding area, the more safely we can conclude that the observed event is connected to the process under analysis.

| Activation Value (a) | 0.0000=<a<0.1000 | 0.1000=<a<0.2000 | 0.2000=<a<0.3000 | 0.3000=<a<0.4000 | 0.4000=<a<0.5000 | 0.5000=<a<0.6000 | 0.6000=<a<0.7000 | 0.7000=<a<0.8000 | 0.8000=<a<0.9000 | 0.9000=<a<1.0000 |
|---|---|---|---|---|---|---|---|---|---|---|
| Probabily | 0.18095 | 0.277175 | 0.2025 | 0.12845 | 0.0813 | 0.05285 | 0.034675 | 0.022125 | 0.013375 | 0.0066 |

*Table 2: Probability that a new event may belong to any of the areas derived by the space splitting according to the activation value of the points (connect this table to Figure 9).*

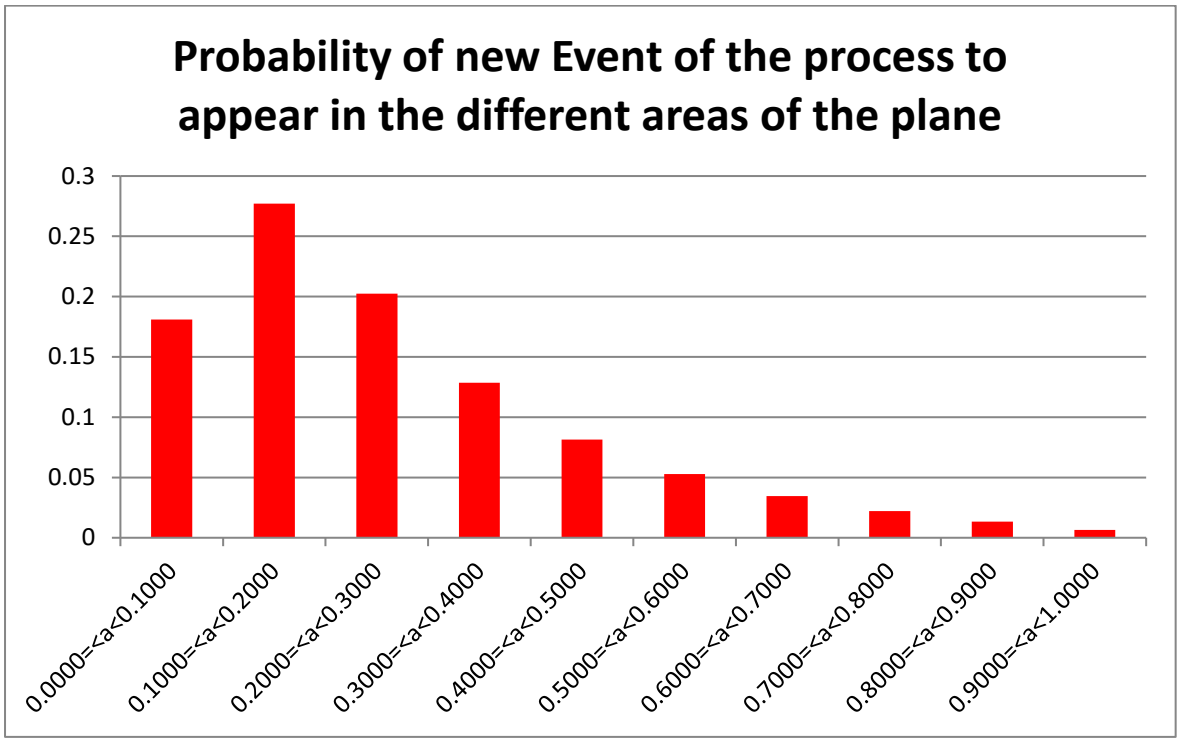

*Figure 9: Histogram of Table 2*

### 3.2 The Beta scalar field

A second notable scalar field that can be built within the context of the topological approach is derived from the Beta Map. In the Alpha Map, the points with high activation value represent the highly responsive ones at the moment the spatial distribution of points (events) was collected. The Beta Map, instead, defines the probability to observe new points (events) that belong to the same distribution (that is, to the same data generating process). The Beta Map can therefore be regarded as a snapshot of a probability density function which is implicit into the structure of the data set. It is generated on the basis of the weighted distance between each generic point of the space and all of the points belonging to the observed distribution, as detailed in equations (15-18). Specifically, (16) calculates the attraction strength of each of the observed points of the distribution with respect to all the other points of the distribution, and (17) computes the level of activation of each point of the space in terms of its weighted distance from all of the points of the observed distribution. Therefore, the strength of attraction of each of the observed points is the weighing factor used to compute (indirect) distances.

{Euclidean distance, $\{n_{k,j}\}$, between each point of 2D space and each point of the dataset}

(15) $n_{k,j} = \sqrt{(x_k - x_j)^2 + (y_k - y_j)^2}$;

$i, j$ = index for each point of the assigned dataset;
$k$ = index for each point of the plane.

{Strength of attraction of each point of the assigned dataset}

(16) $F_i = \dfrac{1}{N-1} \cdot \sum_{\substack{j=1 \\ j \neq i}}^{N-1} e^{-\dfrac{\overline{d_{i,j}}}{D}}$;

{Activaction strength of each point of the 2D plane}

(17) $b_k = \dfrac{1}{N} \sum_{j=1}^{N} F_j \cdot e^{-\dfrac{n_{k,j}}{D} \beta^*}$;

$N$ = all the points of assigned dataset;
$D = \max_{i,j}\{d_{i,j}\}$; Maximum distance among the points of the assigned dataset;
$\beta^*$ = optimal value of $\beta$ according to equation (12);
$b_k$ = Value of activation of $k$-th point of the plane; $b_k \in [0,1]$.

Once again, it is convenient to fix ideas to refer to our example regarding Etruscan towns. The Beta Map presents a distinct pattern of activation, as shown in Figures 9a-c. The local maxima measure the capacity of attraction of each city in the area, according to the main variables that drive the quasi-epidemic process under study, where the 'events' are observed urban settlements. In this specific case, such variables are the socio-economic forces that shape the relationship between the different towns, and therefore the most active areas are the one where socio-economic activation on the basis of the observed spatial distribution is estimated to be stronger, with a consequently higher probability of future, close-by new settlements. If the distribution of points concerned observed occurrences of an epidemic process, the Beta Map would instead represent the probability of observing future cases of the same epidemic in any given point of the space – that is, of observing not previously known events which had nevertheless already occurred.

The interpretation of the Beta Map, therefore, and more generally of the maps derived from the topological approach, depends on the specific nature of the process under study, be it quasi-epidemic or epidemic, and of which kind. Our example of the spatial distribution of Etruscan towns is of interest in guiding our understanding of some of the basic features of the Beta Map. In the case of ancient settlements like the Etruscan ones, and unlike what happens for events related to an epidemic, determining the exact timing of the 'events' is difficult for at least two separate reasons. In the first place, historical data and information are lacking and fragmentary; even basic data such as the town's population need to be estimated on the basis of indirect evidence. Moreover, being the 'events' settlements, there is not a single moment in time that can be chosen to mark the event 'occurrence': A settlement develops in time, and it is difficult to tell exactly when it ceases to be a small cluster of human dwellings to become a town, and even more to be recognized, as in our case, as one among the twelve major Etruscan towns that made the so-called Dodecapolis among the many Etruscan settlements scattered over the surrounding region and beyond. Consequently, how can we interpret the Beta Map in this case specifically, what is the 'next future' the Map is providing information about? An interesting element here is the fact that the very identification of the towns that were part of the Dodecapolis, which represents the urban network where economic trade and social exchange was most intense in the whole Etruscan world, is itself partly conjectural. For instance, there is controversy about whether the list should include Orvieto (Volsinii Veteres) or Bolsena (Volsinii Novi) [24]. Our selection includes Bolsena and not Orvieto, and being Bolsena a settlement that acquired importance in the later phase of Etruscan civilization, this should imply that the corresponding pattern should refer to a late rather than an early phase of the quasi-epidemic process.

In our Beta Map, the Etruscan towns with the highest activation are, primarily, Vulci and Tarquinia, and to a slightly lesser degree Bolsena; there is a second circle of less activated but still responsive cities including Roselle, Vetulonia and Chiusi, a somehow intermediate activation level in Cerveteri (Caere), and an outer circle of low-activation cities including the five remaining ones. It must be remarked that the spatial distribution of cities provides no specific information as to their population, or to its economic and social importance. Nevertheless, the mere structure of the 'space grammar' of the distribution should allow to derive relevant aspects of such underlying quasi-epidemic phenomena as if they were the outcome of a proper epidemic process.

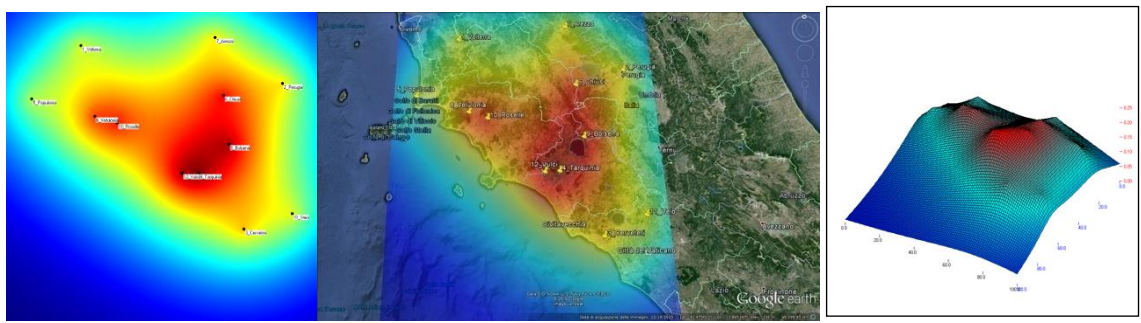

*Figure 9a,b,c: Beta Map(a), its geographic projection using Google earth(b) and 3D graphic (c).*

From what we know of the population of the twelve Etruscan towns in the Archaic period at the peak of their success [25], Veio was by far the most populous, followed by Cerveteri and Populonia, then by Tarquinia, and further by Vetulonia, Volterra, and Vulci. This hierarchy does not reflect at all in the Beta Map shown in Figure 9. However, starting from the Fourth Century BC, the hierarchy influence changes. Veio, which reaches the apogee of its influence between the Seventh and the Fifth century BC is totally destroyed by the Romans in 396 BC and never rebuilt. Populonia loses momentum after the Sixth century BC, and also for Voleterra the peak of influence is reached around the Fifth century BC. Cerveteri is still influent in the Fourth century BC but once again its apogee is in the previous two centuries. On the contrary, in the Fourth century Vulci and Tarquinia are gaining momentum in relative terms. Tarquinia, in particular, reaches its apogee in the Fourth century [24]. The map seems therefore to refer to the transition between the Fifth and Fourth century scenarios. The inclusion of Bolsena instead of Orvieto probably contributes to this time focus, defining a geography of influence that tracks the late rather than early Dodecapolis.

In quasi-epidemical models, therefore, the inherent ambiguity of the time coordinates of the observed 'events' calls for some specific intelligence like in the example above, further underlining the necessity of thinking of the topological approach in terms of a contextually situated analysis. In the specific example of the Etruscan towns, we do not consider of course our analysis as a scientific contribution to Etruscan studies, as a detailed interpretation of the maps would call for the profound insight only accessible to disciplinary specialists. We only consider this exercise an interesting benchmark to help us understand the temporally and contextually situated nature of the inferences of the topological approach for quasi-epidemical processes. In this vein, it is meaningful to remark that the mere analysis of the topological structure of the spatial distribution of the Etruscan towns reflects with such accuracy certain trends that involve the working of complex economic, social and even military factors.

### *3.3 The Gamma scalar field*

We now introduce another quantity which is helpful in the analysis of the distribution of points in space in epidemic and quasi-epidemic processes. When considering distances on the two-dimensional plane, generally we refer to Euclidean distances as measured by a straight line connecting the two points. This notion of distance is clearly insensitive to the distance of any of the two points with respect to their distances to other points in the plane. Therefore, if the position of another point C in the space changes, the Euclidean distance between A and B remains unaffected.

However, in terms of our topological approach, as the position of each point in the spatial distribution carries a strong meaning in terms of the global organization, we have to take into account how variations in the positions of certain points may influence the whole spatial structure from the perspective of any single point belonging to it, in its relations with the other points. To this purpose, we re-define the weights for the computation of the TWC in terms of a new free parameter $\gamma$ which acts as a modifier of the Euclidean distance between points, so as to suitably 'tune' the optimal deformation needed to capture the actual spatial organization of the points according to a logic that is similar to the one already followed for the construction of the Alpha and Beta Maps.

By computing a sequence of values of the $\gamma$ parameter, ranging from 1 (no modification of the Euclidean distance) to infinity, we see how, as $\gamma_i(t)$, that is, the deformation parameter assigned to the *i*-th point in the distribution, increases, the influence of the *i*-th point in affecting the global organization of distances increases. In other words, a high $\gamma_i(t)$ signals that the position of that point is highly critical in determining the expected evolution of the spatial distribution. Think for instance, in quasi-epidemic terms, of how the geography of trade and social exchange determining a certain distribution of urban settlements would be modified if one of those would suddenly acquire the status of a State Capital. All at a sudden, the relative distances of all the other points with respect to the new Capital would matter much more than in the past in determining whether another given settlement is now considered as 'far' or 'central' in the now re-defined spatial organization. If we take the center of mass of the distribution as the natural reference, we can therefore consider how the distance between the center of mass and any given point is influenced by the relative distance of all the other points in the distribution.

We therefore define this new notion of distance in terms of a vector of points, computed for a given sequence of values of $\gamma(t)$, which need not lie on a straight line but can also draw a curvilinear shape. Formally:

(18)
$$TWC_x(\gamma_i(t)) = \frac{1}{\sum_{j=1}^{N} p_{i,j}(\gamma_i(t))} \sum_{j=1}^{N} p_{i,j}(\gamma_i(t)) \cdot x_j;$$

$$TWC_y(\gamma_i(t)) = \frac{1}{\sum_{j=1}^{N} p_{i,j}(\gamma_i(t))} \sum_{j=1}^{N} p_{i,j}(\gamma_i(t)) \cdot y_j;$$

Where:

$$p_{i,j}(\gamma_i(t)) = e^{-\frac{d_{i,j}}{D} \gamma_i(t)};$$

$$\gamma(t+1) = \gamma(t) + \varepsilon \quad \varepsilon = \text{small positive quantity;}$$

$$\lim_{\gamma(t) \to \infty} TWC_x(\gamma_i(t)) = x_i$$
$$\lim_{\gamma(t) \to \infty} TWC_y(\gamma_i(t)) = y_i$$

$$D = \max_{i,j}\{d_{i,j}\}; \quad \{\text{Maximum distance among the assigned data set points}\}.$$

Applying (18) to each point in the distribution will generate a vector of coordinates of $N_i$ components for each of these points. It is important to point out that each point can be represented by a vector of different cardinality, $M_i$, because any vector of coordinates, $TWC_x(\gamma_i(t)), TWC_y(\gamma_i(t))$, may converge in a different number of steps (see Table 3).

$$\begin{Bmatrix} TWC_x(\gamma_1(1)) \\ TWC_y(\gamma_1(1)) \end{Bmatrix} \begin{Bmatrix} TWC_x(\gamma_1(2)) \\ TWC_y(\gamma_1(2)) \end{Bmatrix} \begin{Bmatrix} TWC_x(\gamma_1(3)) \\ TWC_y(\gamma_1(3)) \end{Bmatrix} \begin{Bmatrix} TWC_x(\gamma_1(...)) \\ TWC_y(\gamma_1(...)) \end{Bmatrix} \begin{Bmatrix} TWC_x(\gamma_1(M_1)) \\ TWC_y(\gamma_1(M_1)) \end{Bmatrix}$$

$$\begin{Bmatrix} TWC_x(\gamma_2(1)) \\ TWC_y(\gamma_2(1)) \end{Bmatrix} \begin{Bmatrix} TWC_x(\gamma_2(2)) \\ TWC_y(\gamma_2(2)) \end{Bmatrix} \begin{Bmatrix} TWC_x(\gamma_2(3)) \\ TWC_y(\gamma_2(3)) \end{Bmatrix} \begin{Bmatrix} TWC_x(\gamma_2(...)) \\ TWC_y(\gamma_2(...)) \end{Bmatrix} \begin{Bmatrix} TWC_x(\gamma_2(M_2)) \\ TWC_y(\gamma_2(M_2)) \end{Bmatrix}$$

$$\begin{Bmatrix} TWC_x(\gamma_3(1)) \\ TWC_y(\gamma_3(1)) \end{Bmatrix} \begin{Bmatrix} TWC_x(\gamma_3(2)) \\ TWC_y(\gamma_3(2)) \end{Bmatrix} \begin{Bmatrix} TWC_x(\gamma_3(3)) \\ TWC_y(\gamma_3(3)) \end{Bmatrix} \begin{Bmatrix} TWC_x(\gamma_3(...)) \\ TWC_y(\gamma_3(...)) \end{Bmatrix} \begin{Bmatrix} TWC_x(\gamma_3(M_3)) \\ TWC_y(\gamma_3(M_3)) \end{Bmatrix}$$

$$\begin{Bmatrix} TWC_x(\gamma_{...}(1)) \\ TWC_y(\gamma_{...}(1)) \end{Bmatrix} \begin{Bmatrix} TWC_x(\gamma_{...}(2)) \\ TWC_y(\gamma_{...}(2)) \end{Bmatrix} \begin{Bmatrix} TWC_x(\gamma_{...}(3)) \\ TWC_y(\gamma_{...}(3)) \end{Bmatrix} \begin{Bmatrix} TWC_x(\gamma_{...}(...)) \\ TWC_y(\gamma_{...}(...)) \end{Bmatrix} \begin{Bmatrix} TWC_x(\gamma_{...}(...)) \\ TWC_y(\gamma_{...}(...)) \end{Bmatrix}$$

$$\begin{Bmatrix} TWC_x(\gamma_N(1)) \\ TWC_y(\gamma_N(1)) \end{Bmatrix} \begin{Bmatrix} TWC_x(\gamma_N(2)) \\ TWC_y(\gamma_N(2)) \end{Bmatrix} \begin{Bmatrix} TWC_x(\gamma_N(3)) \\ TWC_y(\gamma_N(3)) \end{Bmatrix} \begin{Bmatrix} TWC_x(\gamma_N(...)) \\ TWC_y(\gamma_N(...)) \end{Bmatrix} \begin{Bmatrix} TWC_x(\gamma_N(M_N)) \\ TWC_y(\gamma_N(M_N)) \end{Bmatrix}$$

*Table 3: The coordinates of TWC($\gamma_i(t)$) vector for each point of the assigned dataset generated by (18).*

By applying (18), we thus generate a matrix whose columns may have different cardinality, as each of them represents the specific trajectory connecting the center of mass with each specific point of the distribution, as shown in the example in Figure 10. Notice how the different trajectories present their own, specific curvilinear patterns depending on the interaction (a kind of 'gravitational effect') between the given point and the spatial distribution of all the others. Figure 11 reports some more detailed examples of some such TWC($\gamma_i(t)$) trajectories as $\gamma_i(t)$ varies.

The structure of the curvilinear patterns allows us to cluster the points into sub-groups depending on the composition of the attraction effects. A cluster is formed when all the points whose trajectories present a curvature that signals a gravitational effect in a certain direction have a curvature that heads toward other points of the cluster. In Figure 10b, we thus have two clusters, a 'Northern' and a 'Southern' one. Point 12, whose trajectory is a straight line, does not belong to any of the two clusters and therefore functions as a boundary between the two.

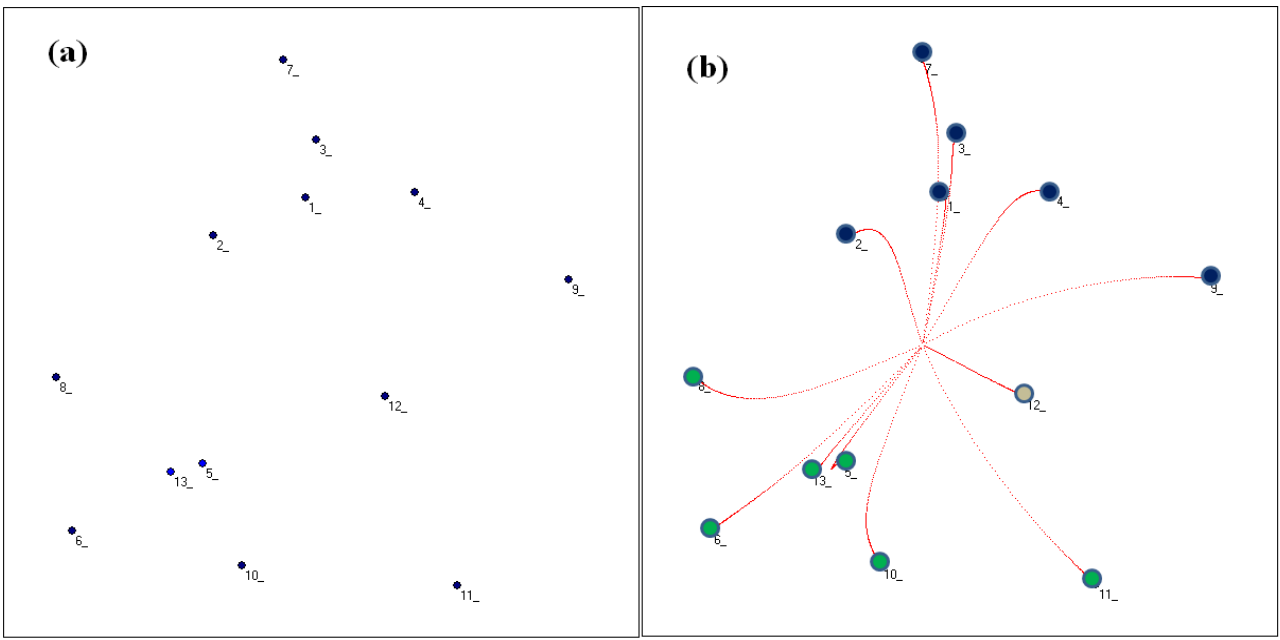

*Figure 10: (a) Initial Positions of 13 points and (b) new points connecting the Centre of Mass to each one of the original points using (18). Green and Blue colors distinguish the implied clusters of attraction. Grey color denotes points lying on the border of the two clusters.*

(a)              (b)

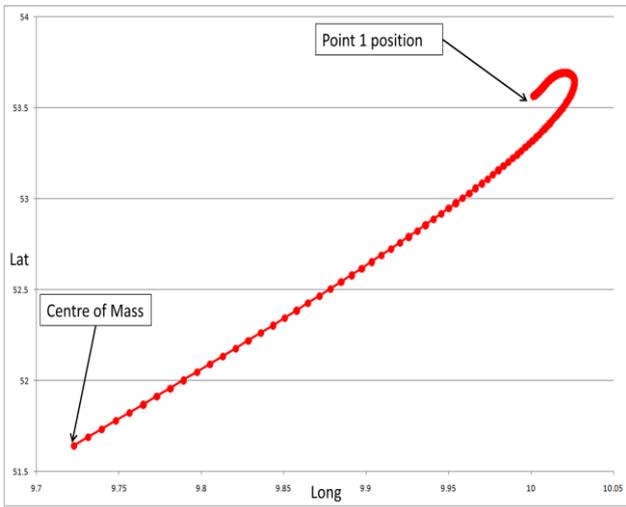
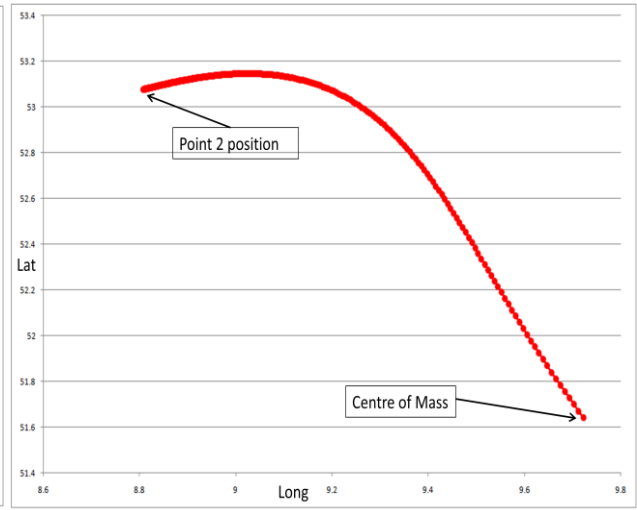

(c)                                                                 (d)

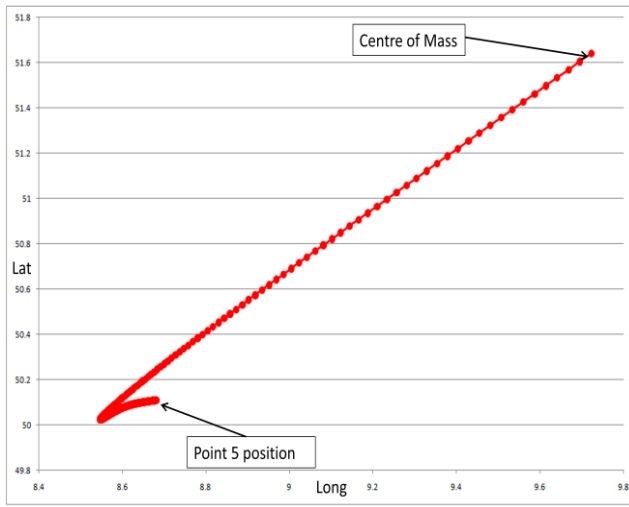
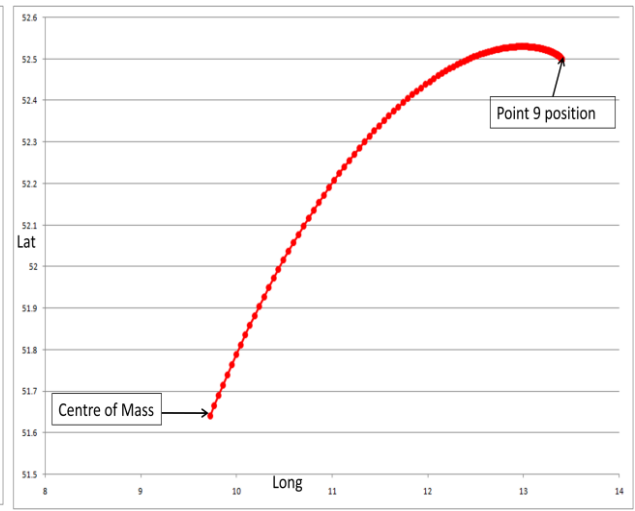

*Figure 11(a,b,c,d): Details of some trajectories of $TWC(\gamma_i(t))$ as $\gamma_i(t)$ varies.*

Let us discuss the example in more detail. Figure 10a presents a spatial distribution made of a sparse collection of isolated points. As these points are 'events' related to a same, underlying phenomenon, we can conjecture that each observed occurrence establishes some form of interaction with the other ones, which is not explicitly represented in our database, but which nevertheless drives the diffusion process. In the case of an epidemic process, this may be the movement of some individual carrying the epidemic agent from one location to another, or in the case of a quasi-epidemic process such as a territorial network of towns, the interaction may consist in trade and social exchanges, migrations of families and workers, and so on. At the early stages of the process, when the phenomenon is still emerging, one may assume that the interactions are mostly stochastic, i.e. the fact that the infectious agent travels in a certain direction or that certain trade relations begin to be established among certain small villages rather than others may be due to random factors. If this is the case, we can conjecture that, statistically, most of such random early interactions are spatially mediated by the center of mass, which represents the average of the coordinates of the all the points belonging to the spatial distribution. For this reason, the matrix of all the TWC(γ$_i$(t)) vectors describing the trajectories of convergence to the various points of the distribution moving away from the center of mass provide a good approximation of the statistics of the interaction among the points at early phase. Moreover, the matrix of trajectories gives useful information about the different attraction strengths of the various points: The more a trajectory is curved, the more it is affected by the attraction of other points located in the corresponding direction (as a sort of 'gravitational' effect). Therefore, the trajectories whose shape is close to a straight line identify the mainstream of the diffusion process, as they are relatively less affected by the attraction strength of the other points while at the same time possibly affecting the trajectories of such points (see e.g. points 3 and 13 in Figure 10b).

The Gamma scalar field will then describe the extent to which each point of the space is activated by its closeness to any of the points belonging to any of the TWC(γ$_i$(t)) trajectories. Formally:

*Legend* :

$M = $ Number of all the points of the plane (window box Rows x Columns);

$N = $ Number of the points of the data set;

$Q_i = $ Number of the new points belonging to the vector $TWC(\gamma_i(t))$, *where* $i \in N$;

$t \in \{1, 2, \ldots, Q_i\}$; indeces of the new points of the vector $TWC(\gamma_i(t))$;

$k \in M$; indeces of the points of the window box on the plane;

$D = \max_{s,z} \{d_{s,z}\}$; Maximum distance among the points of the dataset;

$m_{k,\gamma_i(t)} = $ distance between the $k$-*th* point of the plane and a $TWC(\gamma_i(t))$ point;

$c_k = $ activation of the k-th point of the window box of the plane;

$\beta^* = $ optimal parameter found for $\beta$ (see equation 12).

(19) $\quad m_{k,\gamma_i(t)} = \sqrt{\left(x_k - TWC_x((\gamma_i(t)))\right)^2 + \left(y_k - TWC_y((\gamma_i(t)))\right)^2}$;

(20) $\quad c_k = \dfrac{1}{\sum\limits_i^N Q_i} \sum\limits_{i=1}^{N} \sum\limits_{t}^{Q_i} e^{-\dfrac{m_{k,\gamma_i(t)}}{D}\beta^*} \quad c_k \in [0,1].$

The Gamma Map for the twelve Etruscan towns is shown in Figure 12.

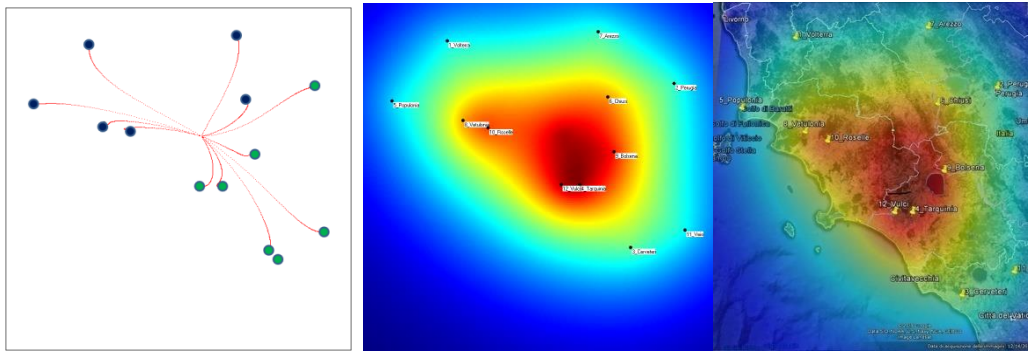

*Figure 12: (a) the Gamma trajectories with two clusters of attraction (Blue, Green), (b) The Gamma scalar field and (c) the same map projected using Google Earth.*

The Gamma Map here presents significant differences with respect of the Beta Map computer for the same spatial distribution. Vulci and Tarquinia are further consolidated as the most activated areas, whereas Vetulonia and Roselle lose in terms of relative intensity of activation, and this also occurs to some extent to Bolsena as well. All towns in the 'outer ring' now present low levels of activation, including not only Veio, but also Cerveteri in the Southern quadrant. In terms of clusters, the Gamma Map separates between the Tuscan towns (Roselle, Vetulonia, Populonia, Volterra, Arezzo, Chiusi) and the others. Perugia, which is the furthermost element of the Southern cluster, is weakly attracted whereas the mutual attraction between Tarquinia, Vulci and Bolsena is particularly strong. Cerveteri and Veio belong to the same cluster but appear less responsive, coherently with the previous analysis from the Beta map.

The Gamma field may be taken as an estimate of the intensity of exchange among the points of the spatial distribution at the epoch that is pertinent to the observed distribution (which, as discussed above, in the case of our example is roughly the Fourth century BC). Our historical evidence concerning the volumes of trade among Etruscan cities in the Fourth century BC is very fragmentary. However, an important correlative of exchange networks is provided by the network of the Etruscan roads connecting the various towns. The documented network of Etruscan roads reported in [26], once the proper superposition with our scalar field is made, seems to be particularly concentrated in the area of highest activation of the Gamma Map, as shown in Figure 13.

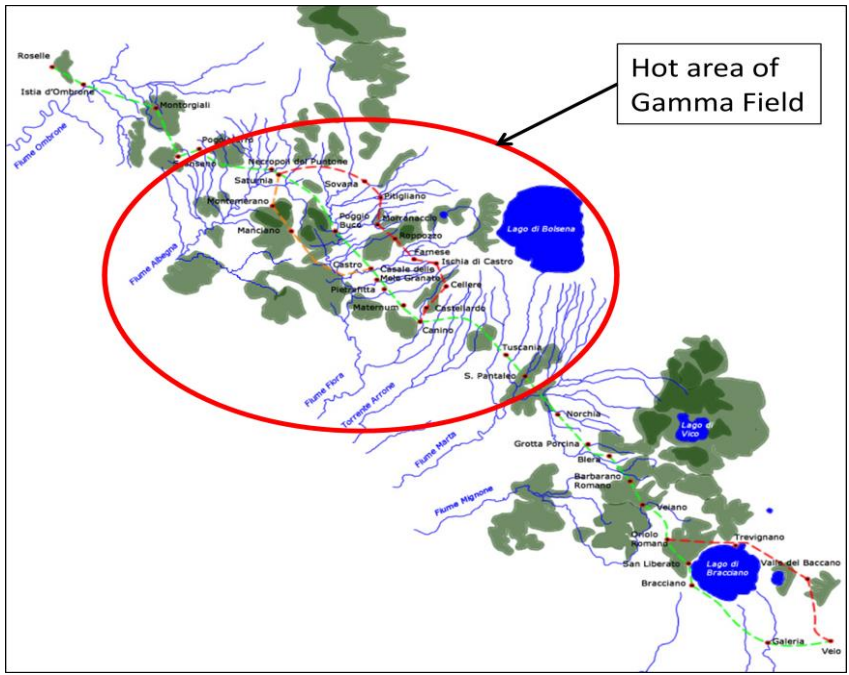

*Figure 13: The networks of the Etruscan roads (copyright Associazione Culturale Etruscan Corner. http://etruscancorner.com/it/associazione-culturale-etruscan-corner/).*

The Gamma scalar field may therefore be thought of as a further qualification of the Beta field, which transforms attraction strengths into intensities of network interaction and therefore highlights longer-term activation patterns within the spatial distribution. As the Beta field is a 'prediction' of the evolution of the Alpha field, the Gamma field may in turn be regarded as a 'prediction' of the evolution of the Beta field. The time scale on which such evolution takes place, as explained in the previous discussion, has to be suitably determined in the specific context and on the basis of the specific nature of the phenomenon under study. In the case of the Etruscan towns, the Gamma field implies that the emerging Tarquinia-Vulci central axis as emerging from the Beta map representation was going to be further consolidated, whereas the already declining major Etruscan towns like Cerveteri and Volterra (not to speak of the already destroyed Veio) were not going to roll down their declining trend even further.

### *3.4 The Theta scalar field*

Whereas the Gamma paths establish nonlinear distances between each pair of points in the space distribution through the mediation of the center of mass, it is also possible to define, as a relatively straightforward extension, nonlinear distances directly between each pair of points. We speak in this case of Theta distance.

We define Theta distance, $TWC\_\theta_{(i,j)_x}(t)$, in terms of the trajectories (that is, the vector of the respective TWC coordinates) that connect directly each pair of points of the data set, according to their nonlinear distances from the center of mass (namely, TWC($\gamma_i(t)$) and TWC($\gamma_j(t)$)). In particular, we can formulate such distances as follows:

$$\text{(21)} \quad \begin{aligned} TWC\_\theta_{(i,j)_x}(t) &= x_j + \left(TWC_{i_x}(\gamma_i(t)) - TWC_{j_x}(\gamma_j(t))\right); \\ TWC\_\theta_{(i,j)_y}(t) &= y_j + \left(TWC_{i_y}(\gamma_i(t)) - TWC_{j_y}(\gamma_j(t))\right). \end{aligned}$$

Applying formula (21), we obtain as in the case of Gamma a new matrix, that we call Theta, that contains all the nonlinear distances between all the points. Each component of the matrix, in turn, is a vector of positions that represent the nonlinear distance between the two points, built following the same logic as for Gamma. The module $|Z_{i,j}|$ of each vector that is part of the Theta matrix is in principle different for each pair of points:

$$\text{(22)} \quad \begin{aligned} |Z_{i,j}| &= Min\{Q_i, Q_j\}; \\ Q_i, Q_j &= \text{ see equation (19)}. \end{aligned}$$

At this stage, we can finally transform the three-dimensional matrix Theta into a direct distance between each pair of points, as follows:

$$\text{(23)} \quad \theta_{i,j} = \sum_{t=1}^{Z_{i,j}-1} D\left(TWC\_\theta_{(i,j)}(t), TWC\_\theta_{(i,j)}(t+1)\right);$$

*where*:

$D() = $ *Distance between two points of the same path.*

Once derived the nonlinear direct distances between each couple of points, we have that the Theta matrix transforms the points of the spatial distribution into a complete regular graph with nonlinear edges, as exemplified in Figure 14 in the case of our Etruscan towns. In this way, for each possible interaction between any couple of towns, we can determine on the basis of the curvature of the corresponding nonlinear edge what is the spatial direction that exerts the most significant interference on that interaction.

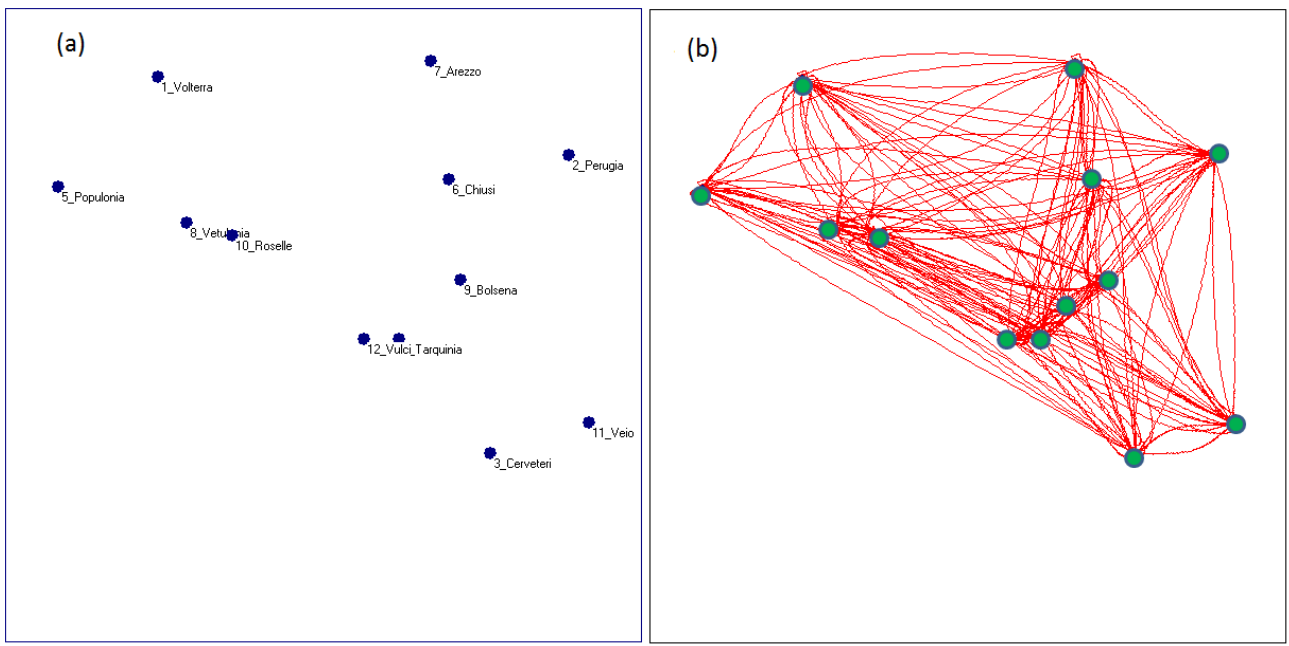

*Figure 14: The 12 Etruscan cities: source data set (a) and complete and regular graph (b) generated from the Theta matrix.*

As the Theta Distance $\theta_{i,j}$ quantifies all the nonlinear connections between the points of the spatial distribution, we can now build from this information the Non-Linear Minimum Spanning Tree (NL-MST) that represents the minimum-energy tree configuration that connects all the points of the distribution, pruning away all the edges that are not indispensable to maintain the connectedness among points by means of the most relevant interactions [27]. In particular, it is obtained as follows:

$$Direct\theta_z = \theta_{i,j} \quad if\ (i \wedge j) \in Mst(\theta_{i,j}) \quad // \text{ Direct distance between } i \text{ and } j;\text{ in Mst path}$$

(24) $\quad DirectPath\_\theta_{z_{x,y}}(t) = TWC\_\theta_{(i,j)_{x,y}}(t) \quad if\ (i \wedge j) \in Mst\left(TWC\_\theta_{(i,j)_{x,y}}(t)\right)$

$\quad // \text{ coordinates of the new points between } i \text{ and } j;$

$\quad where\ z \in \{1,2,...,N-1\}$

(25) $\quad NL\_Mst = \bigcup Direct\theta_z.$

Equation (25) determines NL_MST as the optimal path connecting all the points of the spatial distribution, taking the *Theta* scalar field as the energy function to be minimized. We exemplify NL_MST for our Etruscan towns example in Figure 15. We notice how this graph adds further structural insights into the spatial organization of the main Etruscan settlements. Tarquinia, Vulci and Bolsena are all directly connected to the Alpha Point (whose location, as seen above, is closely connected to the historically more consolidated conjectural location of the Etruscan federation clearinghouse institution, the Fanum Voltumnae). Each of the three towns sits in the MST organization as an autonomous entity, which means that their most significant interaction is with the central clearinghouse rather than with other towns in the federation. This property is coherent with the information arriving from the Alpha and Beta maps as to the centrality of these towns in the global spatial organization at the time of reference. There are, moreover, three distinct branches moving from the Alpha point that represent three different 'peripheral' subsystems: the Tuscan one, mediated by the Maremma towns of Roselle and Vetulonia, with Volterra and Populonia as the faraway 'leaves' of the tree; the Valdichiana one, with the small town of Chiusi as the nexus and other marginal towns like Arezzo and Perugia as leaves; and the Latium one, with Cerveteri as the nexus and Veio as the faraway leaf. The Valdichiana branch is the only one that contains towns belonging to different clusters according to the partition determined by the Gamma Map. This branch is the most marginal one, as it contains towns that even at their peak of influence only played a secondary role in the federation. The Latium branch was the dominant one in the earlier period, thanks to the strategic geographical position of Cerveteri and Veio that allowed to profit from the most active and rewarding channels of trade, but this source of advantage turned into the main threat when Rome rose to power, as demonstrated by the destruction of Veio and by the downplaying of Cerveteri. On the other hand, the towns along the Tuscan branch suffered from the fact that the Maremma region, where the nexus cities of Roselle and Vetulonia sat, was severely threatened by the unhealthy swamps that plagued the region until its final reclamation many centuries after, and this clearly also compromised the influence of Volterra and Populonia within the federation, which despite their strategical position of nodes toward other neighbor populations were poorly connected to the other towns. In conclusion, the triad Vulci-Tarquinia-Bolsena not only promited from the geographical proximity of the three towns to the Alpha point, but also on the fact that they were the only towns in the federation clear of all of the major threats affecting the other towns on the external branches: socio-economic marginality (the Valdichiana branch), environmental barriers (the Tuscan branch), excessive closeness to the increasingly hostile foreign superpower (the Latium branch). What is particularly remarkable is that this wealth of multi-layered information may simply be deduced through several rounds of structural analysis of the spatial distribution of the points of the quasi-epidemic process.

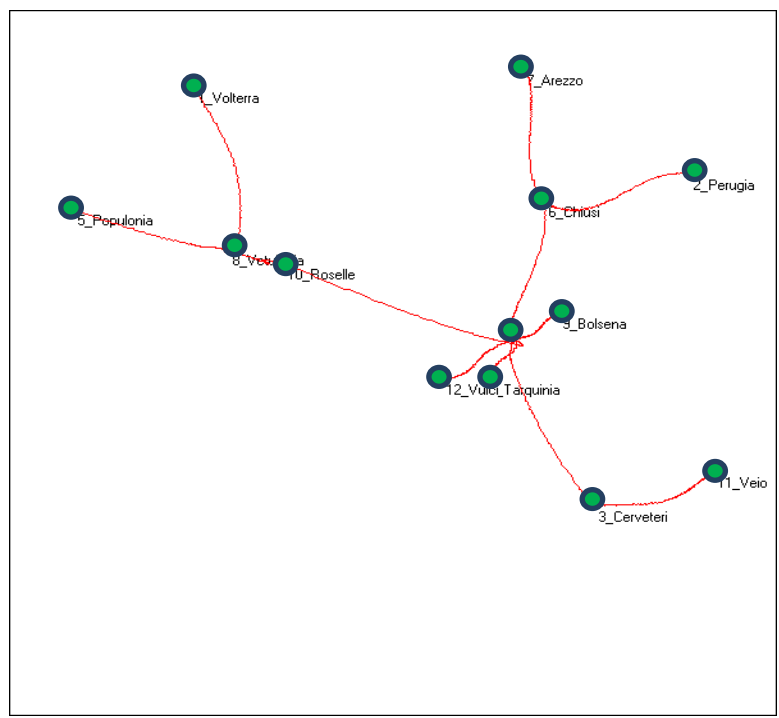

*Figure 15: Non Linear Minimum Spanning Tree of the 12 Etruscan towns.*

The sequence of the Alpha, Beta and Gamma Maps have gradually revealed deeper and deeper structural layers of the spatial organization of points. As a consequence, their predictive power has gradually moved toward the longer term, making Beta the 'future' (i.e. the probability density function) of Alpha, and Gamma the 'future' (likewise) of Beta. The same can be said of Theta, which can be seen in turn as the 'future' of Gamma, with the usual qualification that the appropriate time scales have to be determined and calibrated on a case by case basis. However, at this stage this remains a fascinating theoretical conjecture that has to be proven and tested further. At any rate, the fact that the Theta Map allows us to build the MST that links together the points of the spatial distribution in terms of a minimal network of influence corresponds to the deepest structural layer explored so far.

The *NL_MST* serves as the basis to build the Theta scalar field, as follows:

*Legend* :

$M =$ Number of all the points of the plane (window box Rows x Columns);

$N =$ Number of the points of the data set;

$P_z =$ Number of the new points belonging to the vector $DirectPath\_\theta_{z_{x,y}}(t)$, where $z \in \{1, 2, ..., N-1\}$;

$t \in \{1, 2, ..., P_z\}$; indices of the new points of the vector $DirectPath\_\theta_{z_{x,y}}(t)$;

$k \in M$; indices of the points of the window box on the plane;

$D = \max_{s,z}\{d_{s,z}\}$; Maximum distance among the points of the dataset;

$l_{k,\theta_z(t)} =$ distance between the *k*-*th* point of the plane and a $DirectPath\_\theta_{z_{x,y}}(t)$ point;

$e_k =$ activation of the k-th point of the window box of the plane;

$\beta^* =$ optimal parameter found for $\beta$ (see equation 12).

$$(26) \quad l_{k,\theta_z(t)} = \sqrt{\left(x_k - DirectPath\_\theta_{z_x}(t)\right)^2 + \left(y_k - DirectPath\_\theta_{z_y}(t)\right)^2} ;$$

$$(27) \quad e_k = \frac{1}{\sum_{z}^{N-1} P_z} \sum_{z}^{N-1} \sum_{t}^{P_z} e^{-\frac{l_{k,\theta_z(t)}}{D} \beta^*} \quad e_k \in [0,1].$$

Figure 16 presents the Theta Map of the twelve Etruscan towns. The highest activation area of the Theta map appears smaller than the corresponding area of the Gamma map, and much more concentrated upon an area that contains a particularly dense roads networks, that is the main transportation hub between the Northern and the Southern cluster [28] (keep in mind that information about the road network is not part of the database on which the various Maps of the topological approach are computed). The highest activation in the Theta scalar field focuses upon Tarquinia, which emerges as the town among those of the Dodecapolis which reached its apogee in the Fourth century BC and not in previous periods as the other major Etruscan towns.

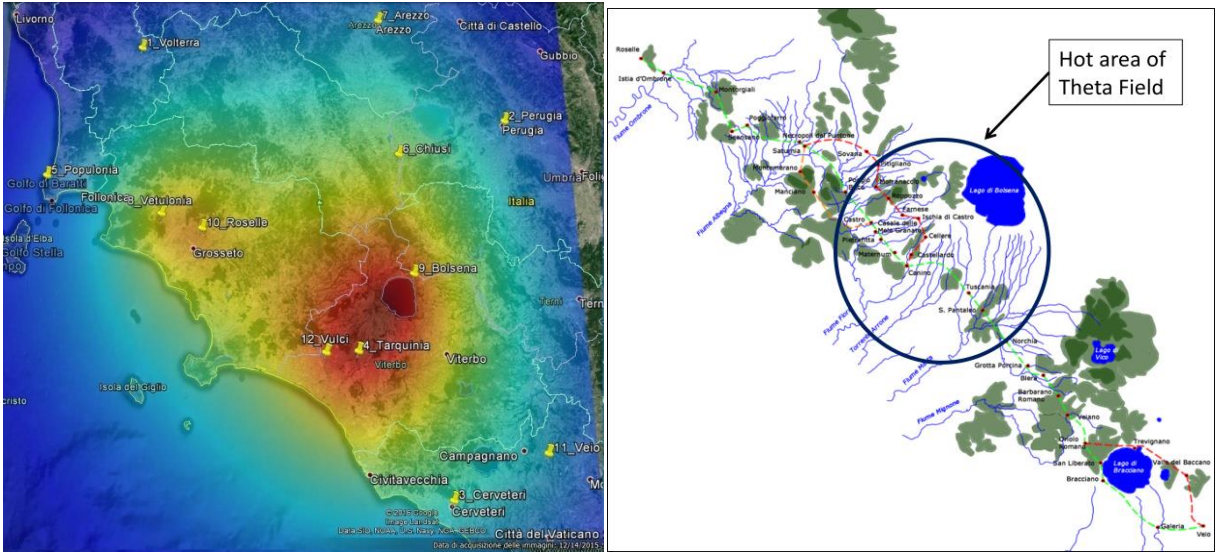

*Figure 16: Theta map of the 12 Etruscan towns and the networks of the Etruscan roads*

### 3.5 Theta Paths as a Probabilistic Machine

The distance matrix $\theta_{i,j}$ (see Equation 23) may be easily transformed into a probability matrix, $p_{i,j}$, which estimates the probability any given point in the spatial distribution communicates with another one as the effect of some random external shock (for instance, a local shortage of some key resource for quasi-epidemic processes, or the occurrence of a mutation of the infective agent in an epidemic process). The closer two points in the spatial distribution, the more likely that the random shock will force an interaction between them:

$$(28) \quad \tilde{\theta}_{i,j} = Max\_\theta - \theta_{i,j};$$
$$where: Max\_\theta = Max\{\theta_{i,j}\};$$

$$(29) \quad p_{i,j} = \frac{\tilde{\theta}_{i,j}}{\sum_{k}^{N} \tilde{\theta}_{i,k}}$$
$$where: \tilde{\theta}_{i,i} = 1.0;$$

The Matrix, $p_{i,j}$, thus measures the probability of transition of a given state from one point to another, when the system is randomly perturbed and when every point exchanges information with the others, where as anticipated closeness implies higher probability of information exchange. This assumption can be taken as the basis for building an explicit diffusion model for our space distribution. Probability values for our Etruscan towns example are reported in Table 4.

| P(i,j) Transition Matrx | | Virtual Time (n+1) | | | | | | | | | | | | |
|---|---|---|---|---|---|---|---|---|---|---|---|---|---|---|
| | | 1_Volterra | 2_Perugia | 3_Cerveteri | 4_Tarquinia | 5_Populonia | 6_Chiusi | 7_Arezzo | 8_Vetulonia | 9_Bolsena | 10_Roselle | 11_Veio | 12_Vulci | T.W.C.($\alpha$) | Row Sum |
| Virtual Time (n) | 1_Volterra | 0.175168 | 0.044575 | 0.019377 | 0.060623 | 0.121947 | 0.071908 | 0.079821 | 0.122558 | 0.059143 | 0.11312 | 0.006215 | 0.060035 | 0.065509 | 1.00000 |
| | 2_Perugia | 0.03892 | 0.152944 | 0.065606 | 0.082269 | 0.01626 | 0.106556 | 0.099499 | 0.047542 | 0.095697 | 0.055342 | 0.071335 | 0.075584 | 0.092446 | 1.00000 |
| | 3_Cerveteri | 0.017071 | 0.066199 | 0.154327 | 0.10303 | 0.016012 | 0.076896 | 0.047089 | 0.04902 | 0.098561 | 0.056264 | 0.112236 | 0.097398 | 0.105895 | 1.00000 |
| | 4_Tarquinia | 0.042167 | 0.065538 | 0.081341 | 0.12184 | 0.039702 | 0.083128 | 0.060843 | 0.065384 | 0.096063 | 0.072499 | 0.071027 | 0.096529 | 0.103939 | 1.00000 |
| | 5_Populonia | 0.13664 | 0.020866 | 0.020364 | 0.063956 | 0.196273 | 0.055449 | 0.053556 | 0.139655 | 0.055065 | 0.128801 | 0 | 0.064119 | 0.065256 | 1.00000 |
| | 6_Chiusi | 0.051905 | 0.088091 | 0.063001 | 0.086267 | 0.035721 | 0.12644 | 0.08782 | 0.061897 | 0.093101 | 0.068683 | 0.061465 | 0.081893 | 0.093715 | 1.00000 |
| | 7_Arezzo | 0.068086 | 0.097203 | 0.04559 | 0.074613 | 0.04077 | 0.103778 | 0.149414 | 0.065813 | 0.083087 | 0.073163 | 0.044468 | 0.070985 | 0.08303 | 1.00000 |
| | 8_Vetulonia | 0.098143 | 0.043603 | 0.044556 | 0.075275 | 0.099809 | 0.068669 | 0.061786 | 0.140273 | 0.070242 | 0.116221 | 0.030211 | 0.074669 | 0.076543 | 1.00000 |
| | 9_Bolsena | 0.041138 | 0.076236 | 0.077814 | 0.096064 | 0.034183 | 0.089716 | 0.067754 | 0.061013 | 0.121842 | 0.06796 | 0.072512 | 0.092421 | 0.101347 | 1.00000 |
| | 10_Roselle | 0.086039 | 0.048209 | 0.048573 | 0.079278 | 0.087432 | 0.072373 | 0.065239 | 0.110388 | 0.074314 | 0.133232 | 0.036429 | 0.077824 | 0.080669 | 1.00000 |
| | 11_Veio | 0.006088 | 0.08004 | 0.124805 | 0.10004 | 0 | 0.083423 | 0.051073 | 0.03696 | 0.10213 | 0.046922 | 0.171608 | 0.092139 | 0.104772 | 1.00000 |
| | 12_Vulci | 0.043259 | 0.062377 | 0.079659 | 0.099999 | 0.041234 | 0.08175 | 0.059965 | 0.067188 | 0.095742 | 0.073727 | 0.067769 | 0.126219 | 0.101112 | 1.00000 |
| | T.W.C. | 0.043783 | 0.070765 | 0.080333 | 0.099874 | 0.038925 | 0.086774 | 0.065058 | 0.063884 | 0.097382 | 0.070885 | 0.071477 | 0.093786 | 0.117074 | 1.00000 |

*Table 4: Transition matrix, $P_{i,j}$, applied to the twelve Etruscan cities.*

This Markovian process allows us to attain yet another level of depth in our analysis, in that we not only have now an understanding of the minimal connectivity structure of the urban network, but also an estimate of the intensity of the exchange flows along the structure. Clearly, flows are not symmetric, and therefore the transition from A to B need not be as likely as the transition from B to A, etcetera. The transition matrix can thus be transformed into a directed weighed graph. We call it G IN-OUT, and we regard it as the most essential representation of the actual flows generated by the underlying process (of a commercial, social, infective nature etcetera) that generates the observed spatial distribution of points. Figure 17(a,b) shows an example of computation of G IN-OUT and its projection on a geo-referenced map for our example of the Etruscan towns.

Figure 17a partitions the twelve towns into three sub-graphs, further building up on the intuition deriving from the NL_MST. One groups together the Tuscan branch of the NL_MST as an isolated graph. A second groups the Valdichiana branch of the NL_MST as a branch of the almost fully connected subgraph that includes the Tarquinia-Vulci-Bolsena triad and the Latium branch of the NL_MST. The G IN-OUT shows therefore how the emerging leaders (Tarquinia, Vulci) are strongly connected to the past leaders (Cerveteri, Veio) and to the Alpha Point, which is structurally part of this subgraph. The other subgraphs are either satellites (as in the case of the small towns of Valdichiana) or basically separated clusters (as the Tuscan towns divided by the dangerous environmental barrier of the Maremma swamp). This explains in an even more profound way the reason of the takeover of Tarquinia and Vulci once Veio and Cerveteri have been gradually overcome by the threat of the rise to power of Rome. They are basically the towns that, being closer to the major town in decline, may interact with them more intensively, thus reaping the benefits of their safer position and exploiting this opportunity to gain influence and resources. The other clusters are either too small and remote (Valdichiana) or too segregated for environmental reasons (Tuscany) to compete for this opportunity. The Alpha Point, therefore, stably belongs to the main sub-graph, although lying in the direction of the other sub-graphs, also due to its function as the federal clearinghouse.

Whether this reconstruction is accurate or not from the point of view of our state of the art knowledge of Etruscan history is an open point that needs the opinion of specialists in the field. However, we think that the illustration of the full spectrum of layers of structural analysis that has been possible through this example, and the tentative inferences that can be derived from them, some of which possibly posing questions that have not been formulated or tackled by Etruscan studies so far, is a promising result in terms of the evaluation of the potential of the topological approach for the analysis of quasi-epidemic processes.

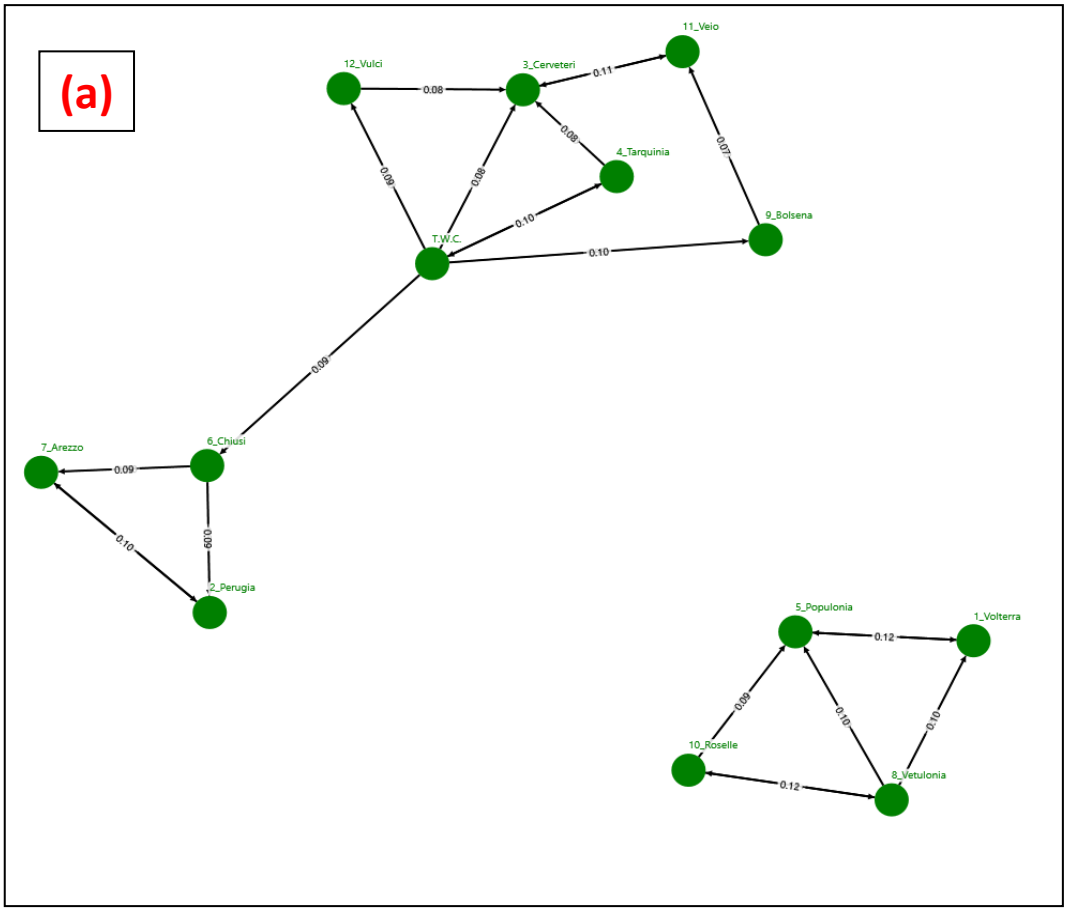

*Figure 17(a): G IN OUT among the 12 Etruscan cities*

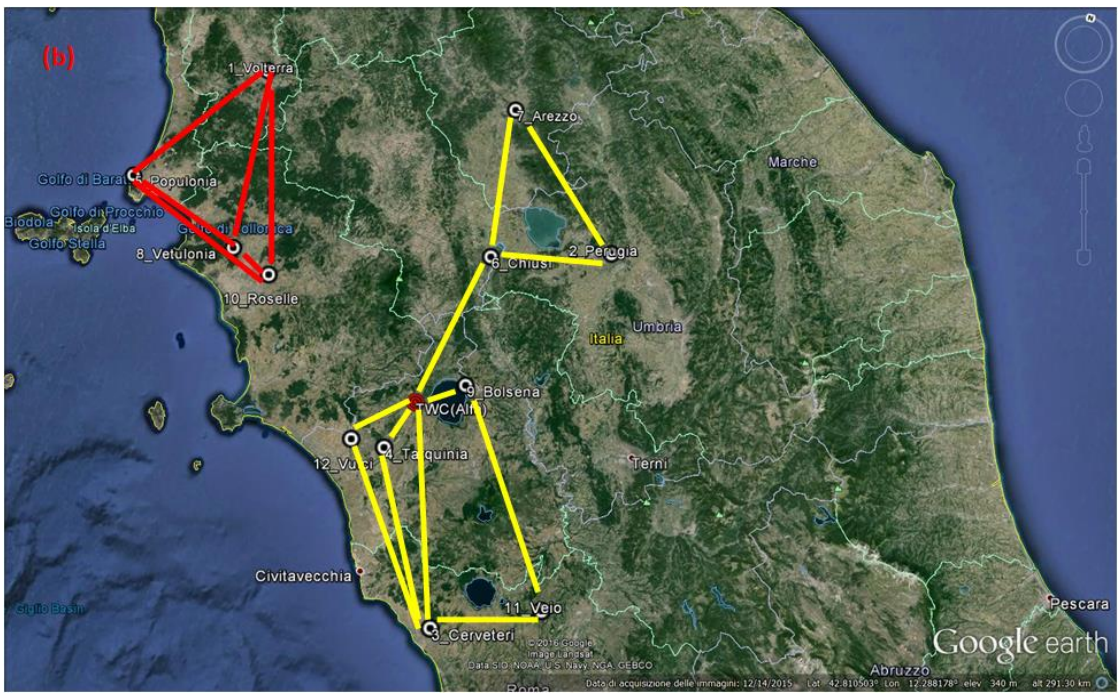

*Figure 17(b): G IN OUT among the 12 Etruscan cities projected on the Google map*

The matrix, $p_{i,j}$ moreover provides us with all the needed information to build up a Discrete Time Markov Chain (DTMC) [29]. We can therefore study an explicit dynamical model that provides still another layer of interpretation of our spatial distribution, and in particular allows us to compute the attractor of the dynamics.

In the example of the 12 Etruscan towns, we have run 10000 random chain simulations. Table 5 summarizes the results. Three different attractors emerge from the simulation:

a. The transition from TWC to Tarquinia and back (53.85%);
b. The transition from Roselle to Vetulonia and back (30.77%);
c. The transition from Veio to Cerveteri and back (15.38%).

| Chains | Trf 1 | Trf 2 | Trf 3 | Trf 4 | Trf 5 |
|---|---|---|---|---|---|
| Chain 1 | 1_Volterra->8_Vetulonia | 8_Vetulonia->10_Roselle | 10_Roselle->8_Vetulonia | | |
| Chain 2 | 2_Perugia->6_Chiusi | 6_Chiusi->9_Bolsena | 9_Bolsena->T.W.C. | T.W.C.->4_Tarquinia | 4_Tarquinia->T.W.C. |
| Chain 3 | 3_Cerveteri->11_Veio | 11_Veio->3_Cerveteri | | | |
| Chain 4 | 4_Tarquinia->T.W.C. | T.W.C.->4_Tarquinia | | | |
| Chain 5 | 5_Populonia->1_Volterra | 1_Volterra->8_Vetulonia | 8_Vetulonia->10_Roselle | 10_Roselle->8_Vetulonia | |
| Chain 6 | 6_Chiusi->9_Bolsena | 9_Bolsena->T.W.C. | T.W.C.->4_Tarquinia | 4_Tarquinia->T.W.C. | |
| Chain 7 | 7_Arezzo->6_Chiusi | 6_Chiusi->9_Bolsena | 9_Bolsena->T.W.C. | T.W.C.->4_Tarquinia | 4_Tarquinia->T.W.C. |
| Chain 8 | 8_Vetulonia->10_Roselle | 10_Roselle->8_Vetulonia | | | |
| Chain 9 | 9_Bolsena->T.W.C. | T.W.C.->4_Tarquinia | 4_Tarquinia->T.W.C. | | |
| Chain 10 | 10_Roselle->8_Vetulonia | 8_Vetulonia->10_Roselle | | | |
| Chain 11 | 11_Veio->3_Cerveteri | 3_Cerveteri->11_Veio | | | |
| Chain 12 | 12_Vulci->4_Tarquinia | 4_Tarquinia->T.W.C. | T.W.C.->4_Tarquinia | | |
| Chain 13 | T.W.C.->4_Tarquinia | 4_Tarquinia->T.W.C. | | | |

*Table 5: The proto typical chains after 10000 run of DTMC. In red color the TWC(α) point.*

The attractor with the largest basin is by far, coherently with the previous analysis, the one connecting the Alpha Point with Tarquinia. On the one hand, this further confirms the appropriateness of the determination of the relevant time period chosen to contextualize the map, as well as the effect of the strategic spatial positioning of Tarquinia in reaping the benefits of the decadence of the nearby major towns in the South at the beginning of the Fourth century BC. A second relevant attractor focuses upon the Roselle-Vetulonia link, that is, the major alternative urban hub of the system, which at its apogee has an even greater population than Tarquinia, but is marginalized by the environmental barrier of the Maremma. And finally, a trace of the major Veio-Cerveteri axis still remains as a minor but still relatively significant attractor of the dynamics, in that the influence of the two major towns is considerably weakened but they nevertheless maintain a relevance in the spatial organization of the Dodecapolis. On the other hand, the Valdichiana cluster is too remote and small to become an attractor of the dynamics. The dynamic model therefore confirms that the spatial organization that is inscribed in the observed distribution of points is the one describing the transition between the old Southern Veio-Cerveteri axis to the new Tarquinia-Vulci one, with Tarquinia as the main settlement, and with the transitional interference of the Tuscan nexus Roselle-Vetulonia which mostly does not succeed in providing an alternative hub for the urban federation despite the potential of its localization as the nexus toward exchanges with other, less threatening neighbors than the Romans who were undermining the power of the Southern larger towns. In particular, the organization still reflects the fading influence of the latter to some extent. The sequential insights made possible by the analysis of the Alpha, Beta, Gamma and Theta Maps gradually allowed us to unveil the emergent preponderance of Tarquinia from the transitional spatial organization depicted in the Alpha Map and therefore to project the eventual development of the hierarchy of influence across the urban federation.

As a general point, we can therefore define a temporal implication among the four maps (α,β,γ, and θ), that can summarized in Figure 17.

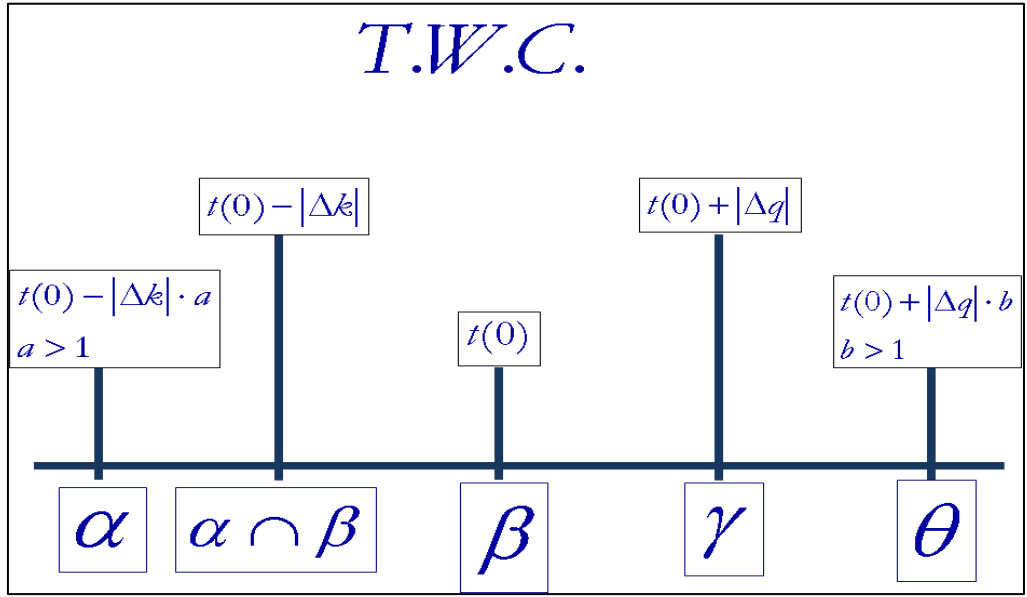

*Figure 17: Temporal implication between the four TWC maps.*

Figure 17 describes as a general rule the time relationships between the four maps:

    a. t(0) indicates the time when data points have been collected (to be properly contextualized for pseudo-epidemic processes);
    b. α∩β indicates the intersection between α-Map and β-Map;
    c. *a, b, Delta k* and *Delta q* are free parameters to be defined experimentally.

As already emphasized above, the specific relationship between the maps has to be determined on a case by case basis from the available information and from a careful analysis of the 'border conditions' defined by the data, as shown in the Etruscan towns example.

## 4. Analysis of epidemic processes

Now that we have tested our approach on the difficult benchmark of a quasi-epidemic process, and in particular on one, such as the case of the Etruscan towns, for which available information was fragmentary and the time frame of reference especially ambiguous, we are in the position to put the topological approach at work on some cases of epidemic processes. Here, unlike the previous example, data are generally abundant and the time frame of the process precisely determined. This allows us to test the approach in another demanding sense, that is, in terms of the precision of the inferences that can be made against a clear, measurable benchmark.

*4.1 The German HUS epidemics, May 2011*

An unusually high number of cases of Haemolytic Uremic Syndrome (HUS) had been observed in Germany since early May 2011. HUS is a serious and sometimes deadly complication that can occur in bacterial intestinal infections with Shiga toxin (syn. verotoxin) producing Escherichia coli (STEC/VTEC). The complete clinical picture of HUS is characterized by acute renal failure, haemolytic anemia and reduction of circulating platelets number (thrombocytopenia). Typically, it is preceded by diarrhea that is often bloody. According to statistics generated by the Robert Koch Institute each year, on average one thousand symptomatic STEC-infections and approximately sixty cases of HUS are reported in Germany, affecting mostly young children under five years of age [30]. STEC are of zoonotic origin and can be transmitted directly or indirectly from animals to humans. Ruminants, especially cattle, sheep, and goats, are believed to be the reservoir. Transmission occurs via the faecal-oral route through contact with animals (or their faeces), by consumption of contaminated food or water, or by direct contact from person to person (smear infection). The incubation period of STEC is between two and ten days, with a latency period from the beginning of gastrointestinal symptoms to enteropathic HUS of approximately one week.

It has to be remarked that, during May 2011, data regarding the geographic locations of registered cases were officially unavailable on public website domains. The authors received confidential information regarding the GPS data of the first 13 locations with at least one proved HUS case on May 30, from a person who was attending a congress on environmental toxicity epidemiology in Europe at that time, and was in contact with German epidemiologists following the outbreak (see Table 6).

| ID | STATE | City Used | WHY USED | LAT | LONG | Q |
|---|---|---|---|---|---|---|
| 1 | Hamburg | Hamburg | Exact match | 53°33'55"N | 10°00'05"E | 59 |
| 2 | Bremen | Bremen | Exact match | 53°4'33"N | 8°48'27"E | 11 |
| 3 | Schleswig-Holstein | Kiel | Capital | 54°19'31"N | 10°8'26"E | 21 |
| 4 | Mecklenburg=Vorpommern | Schwerin | Capital | 53°38'0"N | 11°25'0"E | 10 |
| 5 | Hesee - largest city is Frankfurt | Frankfurt | Largest city | 50°6'37"N | 8°40'56"E | 31 |
| 6 | Sarrland | Saarbrücken | Capital | 49°14'0"N | 7°0'0"E | 5 |
| 7 | Lower Saxony | Hanover | Capital | 52°22'N | 9°43'E | 28 |
| 8 | North Rhine - Westphalia | Duesseldorf | Capital | 51°14'N | 6°47'E | 31 |
| 9 | Berlin | Berlin | Exact match | 52°30'2"N | 13°23'56"E | 3 |
| 10 | Baden0Wuttemberg | Stuttgart | Capital | 48°46'43"N | 9°10'46"E | 8 |
| 11 | Bavaria | Munic (München) | Capital | 48°31'52"N | 11°57'50" | 5 |
| 12 | Thuringa | Erfurt | Capital | 50°59'0"N | 11°2'0"E | 1 |
| 13 | Rhineland-Palatinate | Mainz | Capital | 50°0'0"N | 8°16'16"E | 1 |

*Table 6: The 13 locations with at least one proved HUS case on May 30, 2011.*

Preliminary results of a case–control study conducted by Hamburg health authorities demonstrated a significant association between disease and the consumption of raw tomatoes, cucumbers and leafy salads and the attention was focused upon the area of Lower Saxony (between Hamburg and Hannover) where the highest density of cases was registered.

On June 18, 2011 the source of the epidemic outbreak was found. The deadly 0104:H4 strain of E. coli that claimed the lives of nearly 40 Germans was found in the northeast of Frankfurt, in the Erlenbach stream, in the evening of June 17,2011 (for a more in-depth analysis of German HUS using TWC see [21]).

Upon closer scrutiny, two main outbreak points were eventually found: The major one close to Hannover, and a second one in Frankfurt [31, p.25]. Figure 18 reports the final map of the epidemic, as edited by Robert Koch Institute, and Figure 19 presents the digital map of the 13 locations.

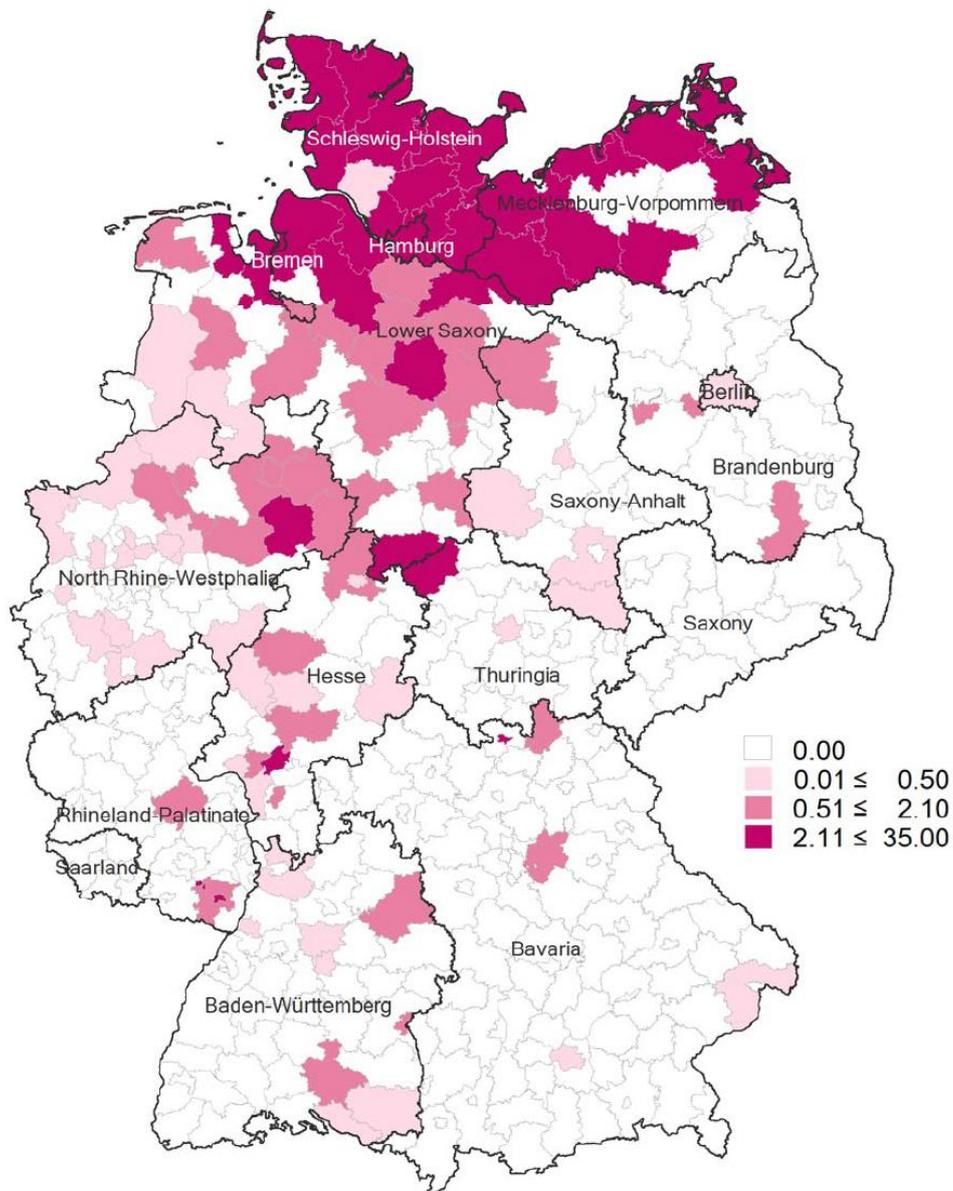

Figure 3: Incidence (cases per 100,000 persons) of HUS during the outbreak, illustrated by county, in which infection likely occurred (residence county, or in cases with travel history, the county of presence at the time of infection).

*Figure 18: Final map of the HUS epidemic, as edited by the Robert Koch Institute.*

Figure 20 shows the Alpha Map and the two TWC(α) points found (the estimated first outbreak). The TWC(α)-determined points locate the epidemic outbreaks in two areas in the vicinities of Hannover and Frankfurt. TWC(α) estimated the point close to Hannover as the more relevant (90%) for the outbreak, and the point close to Frankfurt as the second possible option (10%). This double outbreak was generated using a Leave One Out protocol: We calculated the Alpha point several times by taking away each one of the 13 points at a time.

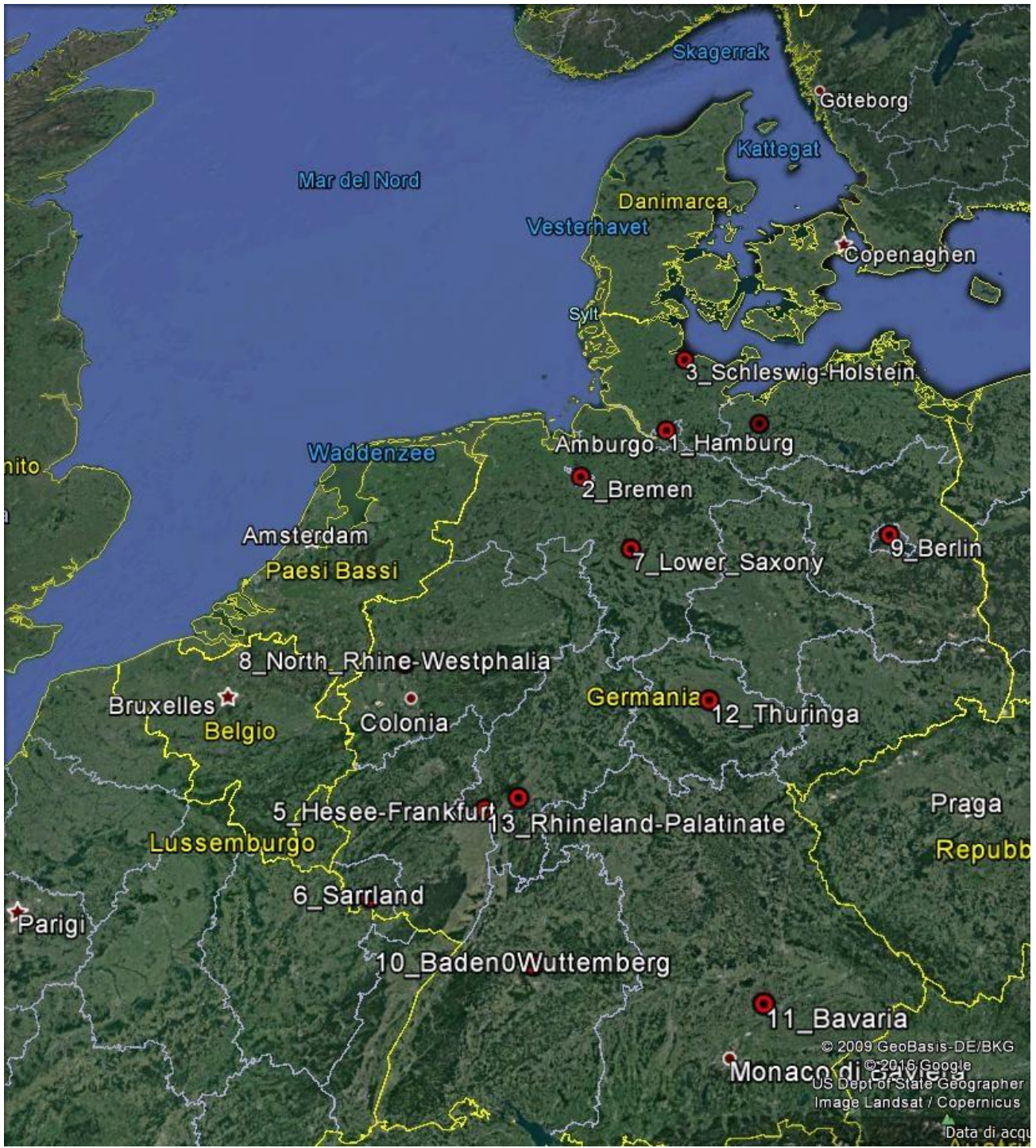

*Figure 19: The digital map of the 13 locations where German HUS was proved on May 30, 2011.*

Figure 21 shows the Beta Map (we remind that TWC algorithms neither consider the frequency of cases observed at a given location nor the chronology of their observation).

The TWC(β), instead, shows two main areas (again, one around Frankfurt and one around Hamburg) as the riskier ones in terms of epidemic diffusion at the time corresponding to the collection of all observed cases at the 13 locations on which the elaboration was based.

Figures 22 and 23 show the Gamma and Theta Maps, respectively. TWC(γ) estimates the diffusion of the epidemic after the data collection time, and TWC(θ) estimates the possible trend of the epidemic after the situation depicted in the Gamma Map. The Gamma Map highlights a rootedness of the epidemic in North and Central Germany. Whereas the Theta map shows a tendency of the epidemic to come back to its two original points, close to Frankfurt and Hannover.

The global dynamics of the TWC(α,β,γ, and θ) may also explain the shifting conjectures about the outbreak locations formulated by the epidemiologists as the epidemic unfolded. They initially observed the epidemic at the time frame corresponding to the Gamma Map (end of May),and as its further evolution led to the situation depicted in the Theta Map, they focused upon the Frankfurt and Hannover areas, eventually discovering, in mid-June, the actual location of the main outbreak.

Further insights are provided by the Non-Linear MST built using the Theta Map parameters (see Figure 24). As the Non-Linear MST overlaps quite precisely with the main highways connecting the 13 observed locations, we have a corroboration of the likely spatial pattern and actual vectors of transmission of the disease.

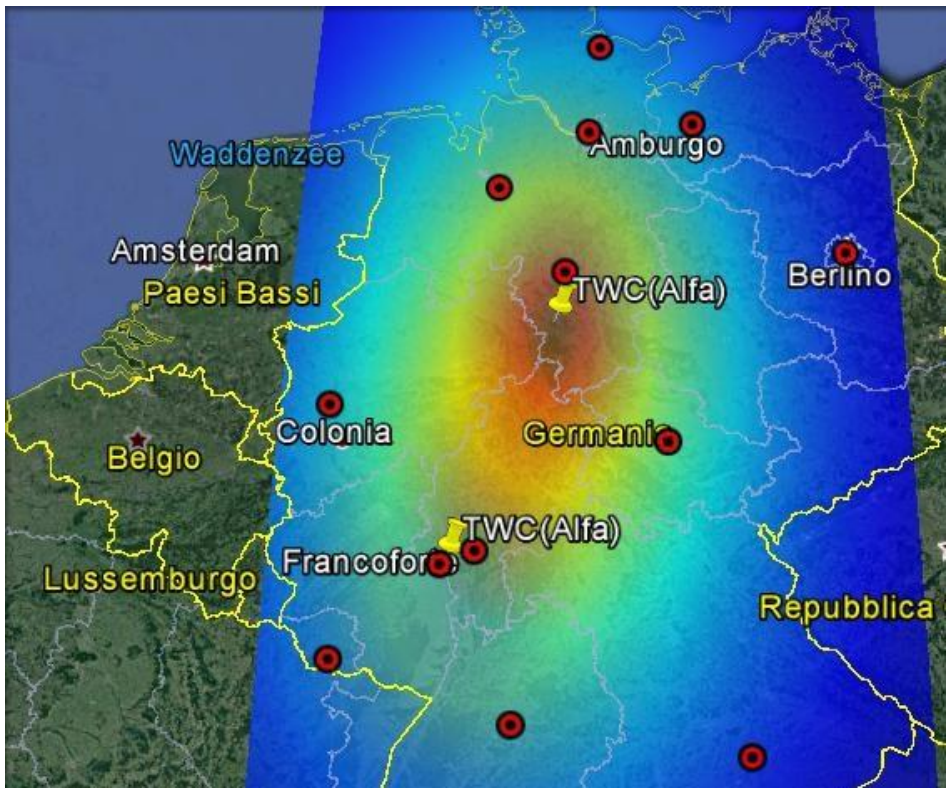

*Figure 20: The Alpha Map, and the two Alpha points.*

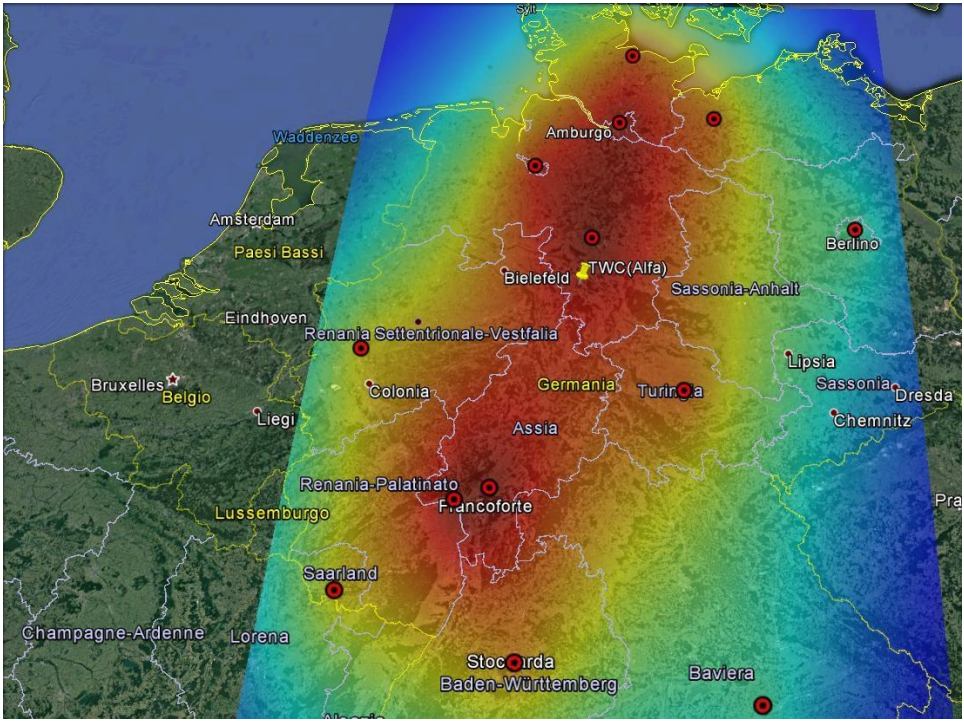

*Figure 21: The Beta Map and the diffusion pattern of the epidemic.*

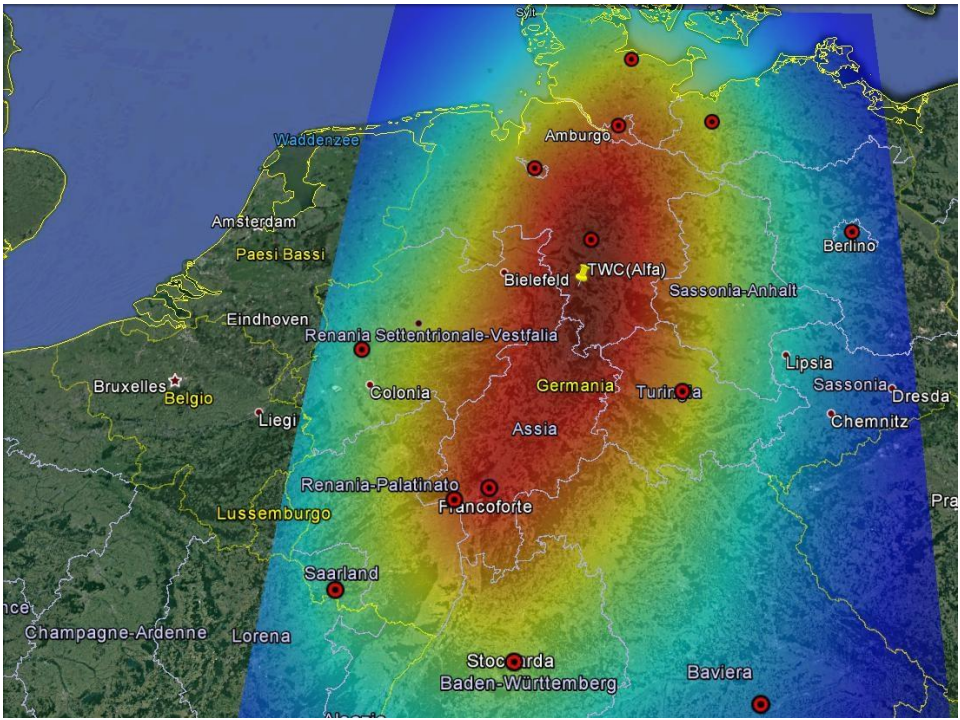

*Figure 22: The Gamma Map and the further diffusion of the epidemic.*

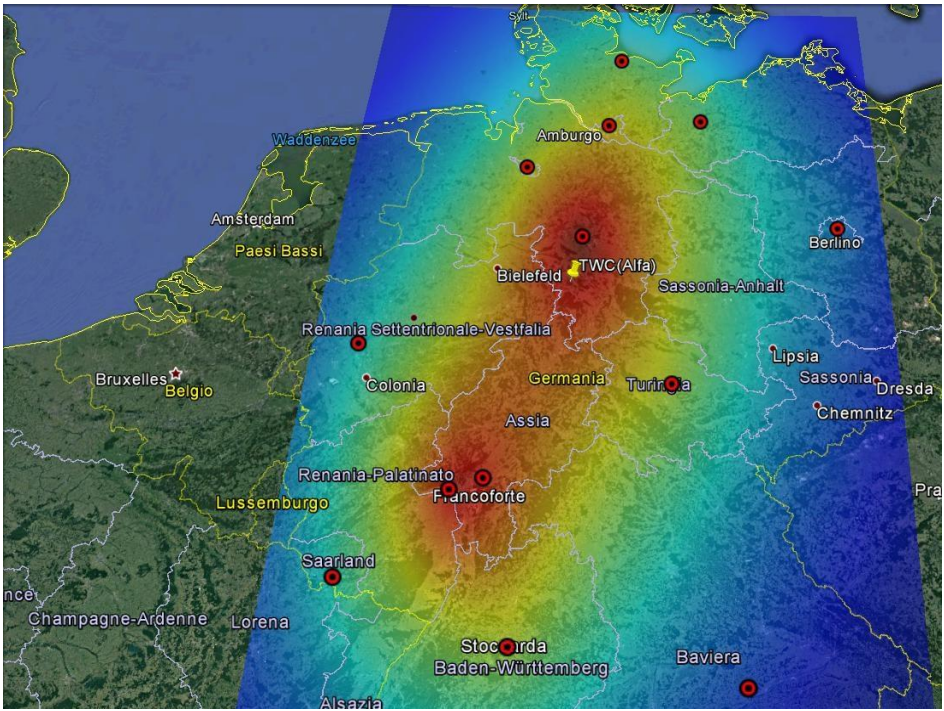

*Figure 23: The Theta Map, and the epidemic 'fold back' to its original locations.*

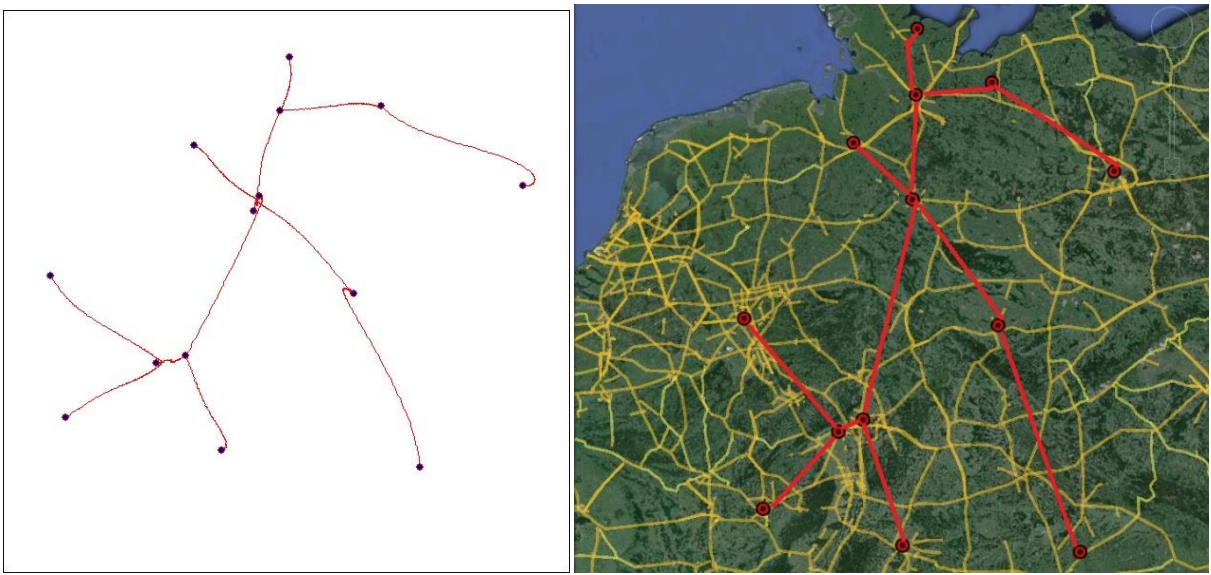

*Figure 24: The Non-Linear MST (on the left), and the real map of the highways connecting the 13 locations (on the right).*

If we now apply the G IN-OUT filter to the transition matrix provided by the Theta parameters (see equations 28 and 29 above), the results are the weighted and direct graphs in Figures 25-26, where the cause-effect structure of the German HUS is shown. We have two main clusters, the first in Lower Saxony and Hamburg and the second around Frankfurt (see the direction of the arrows in Figure 25). The TWC point functions as a sort of an 'epidemic engine', pumping the disease from the first Northern German cluster toward Southern Germany. The orientation of the arrows shows that the Frankfurt outbreak is independent of the main outbreak, and merges with it as the former comes through Lower Saxony.

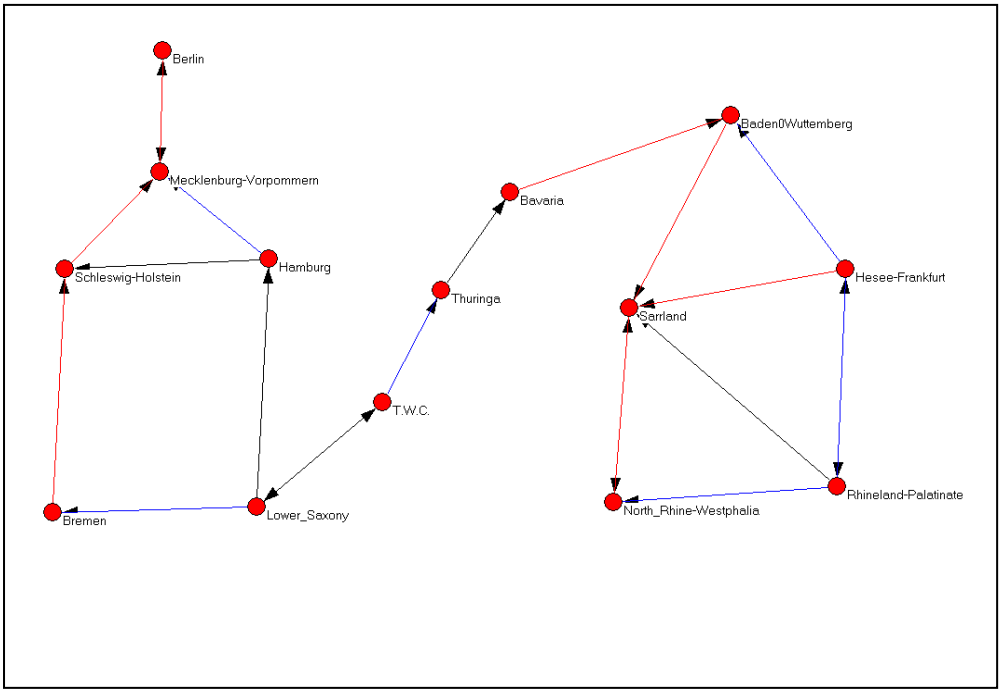

*Figure 25: The G IN OUT graph that explains the cause-effect structure of the German HUS epidemic.*

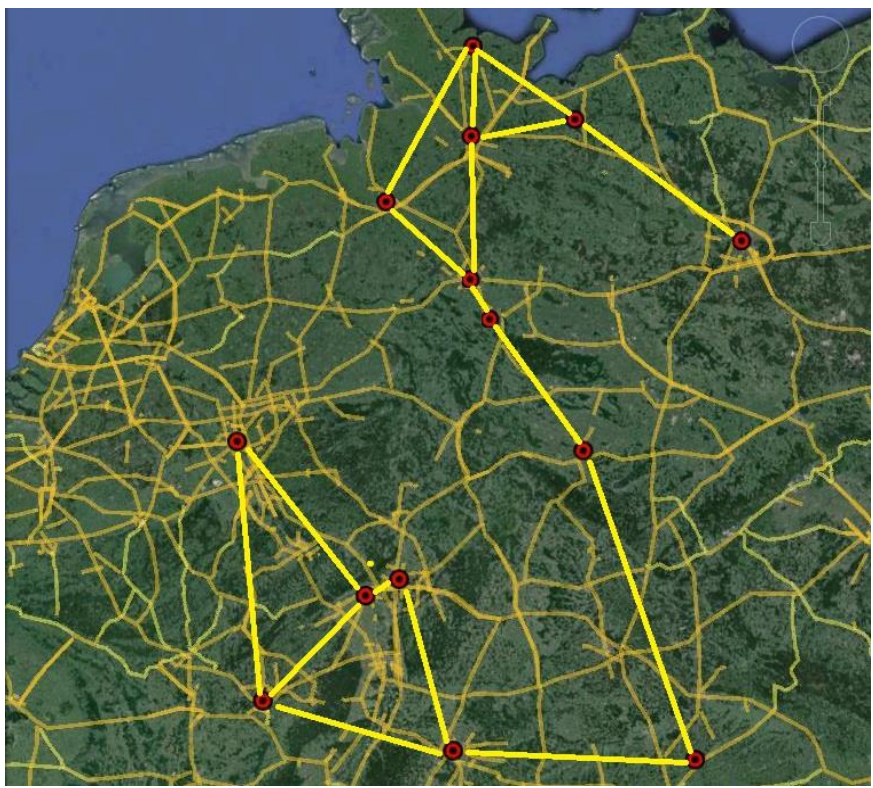

*Figure 26: The G IN-OUT graph connections as applied to the real map.*

The prototypical chains of the DTMC show us the attractors of the epidemics, namely the places where the German Hus tends to be resilient (see Table 7).

| Chains | Trf 1 | Trf 2 | Trf 3 | Trf 4 | Trf 5 |
|---|---|---|---|---|---|
| Chain 1 | Hamburg->Schleswig-Holstein | Schleswig-Holstein->Lower_Saxony | Lower_Saxony->Schleswig-Holstein | | |
| Chain 2 | Bremen->Hamburg | Hamburg->Schleswig-Holstein | Schleswig-Holstein->Lower_Saxony | Lower_Saxony->Schleswig-Holstein | |
| Chain 3 | Schleswig-Holstein->Lower_Saxony | Lower_Saxony->Schleswig-Holstein | | | |
| Chain 4 | Mecklenburg-Vorpommern->Hamburg | Hamburg->Schleswig-Holstein | Schleswig-Holstein->Lower_Saxony | Lower_Saxony->Schleswig-Holstein | |
| Chain 5 | Hesee-Frankfurt->T.W.C. | T.W.C.->Hesee-Frankfurt | | | |
| Chain 6 | Sarrland->11Rhineland-Palatinate | Rhineland-Palatinate->Hesee-Frankfurt | Hesee-Frankfurt->T.W.C. | T.W.C.->Hesee-Frankfurt | |
| Chain 7 | Lower_Saxony->Schleswig-Holstein | Schleswig-Holstein->Lower_Saxony | | | |
| Chain 8 | North_Rhine-Westphalia->Rhineland-Palatinate | Rhineland-Palatinate->Hesee-Frankfurt | Hesee-Frankfurt->T.W.C. | T.W.C.->Hesee-Frankfurt | |
| Chain 9 | Berlin->Mecklenburg-Vorpommern | Mecklenburg-Vorpommern->Hamburg | Hamburg->Schleswig-Holstein | Schleswig-Holstein->Lower_Saxony | Lower_Saxony->Schleswig-Holstein |
| Chain 10 | Baden0Wuttemberg->T.W.C. | T.W.C.->Hesee-Frankfurt | Hesee-Frankfurt->T.W.C. | | |
| Chain 11 | Bavaria->Thuringa | Thuringa->Hesee-Frankfurt | Hesee-Frankfurt->T.W.C. | T.W.C.->Hesee-Frankfurt | |
| Chain 12 | Thuringa->Hesee-Frankfurt | Hesee-Frankfurt->T.W.C. | T.W.C.->Hesee-Frankfurt | | |
| Chain 13 | Rhineland-Palatinate->Hesee-Frankfurt | Hesee-Frankfurt->T.W.C. | T.W.C.->Hesee-Frankfurt | | |
| Chain 14 | T.W.C.->Hesee-Frankfurt | Hesee-Frankfurt->T.W.C. | | | |

*Table 7: The prototypical transitions of the German HUS epidemic according to 1000 runs of DTMC (in red the attractors for each chain).*

Table 7 highlights the two key transitions that are the attractors of the epidemic dynamics:

a. The transition between Frankfurt and Rhineland-Palatinate, and vice versa (29%);
b. The transition between the TWC point and Lower Saxony, and vice versa (61%).

This application, with its entirely known benchmarks, clearly exemplifies the amount of key information that it is possible to extract from a small dataset using TWC algorithms.

### 5. The TWC algorithms: a synthesis

The topological approach can thus generate a number of useful quantities for the analysis of epidemic and quasi-epidemic processes from the mere information embedded in the spatial distribution of the observed events related to the phenomenon under study, namely:

a. The Alpha point, TWC($\alpha^*$), that gives an estimation of the possible location of the outbreak on the basis of the spatial pattern of the observed events (see Equations 9a-9b);
b. An Alpha vector, TWC $\alpha$ (n), that describes the transition from the center of mass to the Alpha point (see Equations 4a-4b);
c. A Beta parameter, $\beta^*$, that represents the optimal value of the bell width of the interaction of each point with the others (see Equations 11a-11b, 12);
d. The Alpha map, $\alpha$-Map, that describes the high activation area where the process that has generated the observed events is likely to have begun (see Equations 13-14);
e. The Beta map, $\beta$-Map, that describes the probability density function of near future events on the basis of the spatial distribution of the observed ones (see Equations 15, 16, 17);
f. The Gamma map, $\gamma$-Map, that, making use of the Gamma Paths, approximates the evolution of the probability density function as determined by the Beta-Map (see Equations 18,19, 20);
g. The Theta map, $\theta$-Map, that, making use of the Theta NL-MST, approximates the further evolution of the probability density function as estimated by the Gamma Map (Equations 21-27);

h. The Theta Transition Model, that is, the *matrix $P_{i,j}$* (Equations 28-29), that enables the simulation of the dynamical model of the exchanges among the points, and that allows also to project the observed points into a weighted graph with direct connections, based on the parameters defined by the distance matrix $\theta_{i,j}$.
i. The Theta prototypical transitions, that represents the main attractors of the interactions (information exchanges) among the observed points (see Tables 5 and 7).

This is the basic toolkit of the topological approach.

## 6. Advanced theory

### 6.1 The TWC and the Meta-distance

Despite their interesting properties analyzed in the previous sections, the quantities generated by means of the basic topological approach also suffer from some limitations. In the first place, their robustness with reference to statistical errors or sampling anomalies in the observed distribution of points has to be evaluated. Moreover, the TWC algorithms present some intrinsic limitations, such as for instance the fact that the TWC points are constrained to lie within the convex hull of the observed points (see Figure 27).

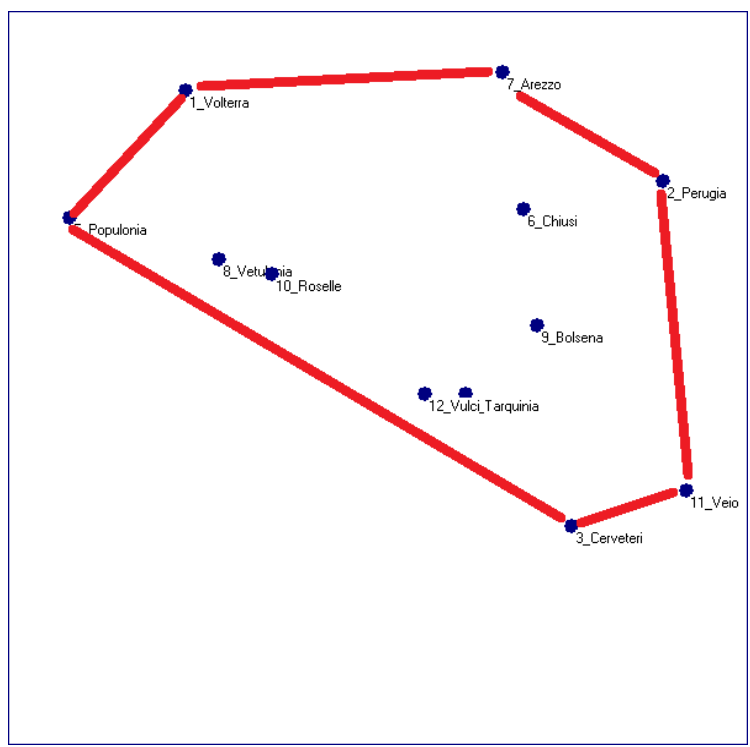

*Figure 27: The convex hull of the twelve Etruscan towns. TWC algorithms (α, β, γ, and θ) are constrained to find points within this perimeter.*

The convex hull limitation is not exclusive of the topological approach but is a common feature of most algorithms designed with the same purposes. For instance, gravitational (or geometric) approaches share the same problem [13, 14]: The estimated outbreak may not fall outside the convex hull for intrinsic mathematical reasons. Likewise for the algebraic approach: The more the estimation of a new point is further away from the convex hull of the assigned set of points, the larger will be the error. However, the topological approach may be further upgraded to overcome this latter limitation, by introducing the concept of *Meta-Distance*.

To illustrate this new concept, let us start again from the familiar notion of Euclidean distance with its well-known characteristic properties:

$$(30) \quad d_{i,i} = 0;$$

$$(31) \quad d_{i,j} = d_{j,i};$$

$$(32) \quad d_{i,j} + d_{i,k} \geq d_{j,k}.$$

The new notion of Meta-distance may be simply put as follows: A Meta-distance among points is the distance matrix of their distances. To fix ideas, consider a two-dimensional space and N points distributed across it. The Euclidean distance between any two of them will be

$$(33) \quad d_{i,j} = \sqrt{(x_i - x_j)^2 + (y_i - y_j)^2};$$

*where*:

$i, j \in N$

*and*

$N$ = Number of Points.

The Meta-distance matrix among the distances of the N points will then be computed as

$$(35) \quad d_{i,j}^{[t+1]} = \sqrt{\sum_{k}^{N} \left( d_{i,k}^{[t]} - d_{j,k}^{[t]} \right)^2};$$

*where*:

$t \in \{1, 2, ..., \infty\}.$

In other words, we consider the familiar Euclidean distance (33) as the first iteration of a recursive process. If we repeatedly apply (34), the sequence will diverge because of the classical Pythagorean inequality. However, if each distance is linearly scaled in the [0,1] lattice, the sequence will now converge.

Specifically, we have

(35) $d_{i,j}^{[t=1]} = \sqrt{(x_i - x_j)^2 + (y_i - y_j)^2}$;

/* first distance matrix between N points in 2 dimensional space*/

(36) $m_{i,j}^{[t=1]} = f(d_{i,j}^{[t=1]})$;

/* Linear scaling of the first distance matrix into the interval [0,1] */

(37) $d_{i,j}^{[t+1]} = \sqrt{\sum_{k}^{N}(m_{i,k}^{[t]} - m_{j,k}^{[t]})^2}$;

/* first Meta Distance matrix of N points in N dimensional space*/

(38) $m_{i,j}^{[t+1]} = f(d_{i,j}^{[t+1]})$;

/* Linear scaling of the first Meta Distance matrix into the interval [0,1] */

*where*:

(39) $1 \leq t \leq \infty$;

(40) $m_{i,j}^{[t]} = f(d_{i,j}^{[t]}) = scale \cdot d_{i,j}^{[t]} + offset$;

(41) $scale = \dfrac{1}{Max(d_{i,j}^{[t]}) - Min(d_{i,j}^{[t]})}$;

(42) $offset = -\dfrac{Min(d_{i,j}^{[t]})}{Max(d_{i,j}^{[t]}) - Min(d_{i,j}^{[t]})}$.

A simple computational procedure to generate the entire set of Meta-distances of N points in a two-dimensional space may be fleshed out as follows:

1. Apply equation (35) to the N original points, to generate the basic distance matrix;
2. Apply equation (36) to the basic distance matrix, to scale the matrix linearly in the interval [0,1];
3. Apply equation (37) to generate the Meta-distance matrix;
4. Apply equation (38), to scale the Meta-distance matrix linearly in the interval [0,1];
5. **If** $\left|m_{i,j}^{[t+1]} - m_{i,j}^{[t]}\right| < \varepsilon, where\ \varepsilon = 0.00000001$, **then** terminate, **else** *t=t+1* and go back to step 3.

The cost function of the Meta-distance algorithm is defined as:

(43)
$$\delta m_{i,j}^{[t+1,t]} = \sum_{i=1}^{N-1}\sum_{j=i+1}^{N}\left|m_{i,j}^{[t]} - m_{i,j}^{[t+1]}\right|;$$
$$\lim_{t \to \infty} \delta m_{i,j}^{[t+1,t]} = 0.$$

The fact that this sequence converges in the unit lattice, and the fact that the matrix generated through this iterative process is not empty, are both interesting in themselves. In fact:

$$(44) \quad \sum_{i=1}^{N-1}\sum_{j=i+1}^{N} m_{i,j}^{[t+1]} > 0 \quad if \left( \sum_{i=1}^{N-1}\sum_{j=i+1}^{N} \left| m_{i,j}^{[t]} - m_{i,j}^{[t+1]} \right| = 0. \right)$$

At this point, since $\forall m(t) \leq 1$, on the basis of fixed-point theorem we can write:

$$(45) \quad \sum_{t=1}^{\infty}\sum_{i=1}^{N-1}\sum_{j=i+1}^{N} \left( m_{i,j}^{[t+1]} - m_{i,j}^{[t]} \right) = \Lambda;$$

$$where: \Lambda \in \square.$$

Therefore, for N given points in a two-dimensional space, the above described iterative process generates a closed set of Meta-distances, due to the convergence of the cost function. Each of the generated Meta-distances clearly has the same cardinality $\|N\|$, where N is also the cardinality of the given set of points in two-dimensional space. Let's denote by $m_{i,j}^{[T]}$ the matrix that is obtained at the last stage of the iterative process (which, as pointed out above, converges). But $m_{i,j}^{[T]}$ is also the locus of the fixed-point, $m^*$, such that $f(m) = m$; in fact:

*stop condition of iteration series :*

$$(46) \quad d_{i,j}^{[t]} = \sqrt[2]{\sum_{k}^{N} \left( m_{i,k}^{[t-1]} - m_{j,k}^{[t-1]} \right)^2};$$

$(47) \quad F(X^{[t]}) = l\left(d_{i,j}^{[t]}\right);$ *linear scaling for* $\forall i, j \in N;$

$(48) \quad G(X^{[t+1]}) = F(X^{[t]});$

$(49) \quad F(X^{[t]}) = G(X^{[t]}).$

The properties of the Meta-distance matrix can be characterized as follows.

### 6.2 Properties of the Meta-distance matrix

We have denoted by $m_{i,j}^{[T]}$ the matrix obtained at the *T*-th step of the Meta-distance process. Its most notable properties are:

a. It is a distance matrix, and all its entries are nonzero;
b. It is a distance matrix that clusters all the points into two sets (see Table 8 below);
c. It is possible to project all the original points in a two-dimensional map, with an optimal approximation, using a Multi-Dimensional-Scaling (MDS) technique.

The two-dimensional map generated by the $m_{i,j}^{[T]}$ matrix re-groups all the original points into two clusters of new mapped points. We name the centroids of these two clusters 'vanishing points', in that they seem to function as the vanishing points of the linear perspective 2D representation of 3D scenes. In other words, they function as vantage points from which the overall spatial organization of the projected points is maximally organized. Each one of these new points may lie outside the convex hull defined by the distribution of the original data points (see Table 8 and Figure 25 below). Once the Meta-distance process has unfolded, we are therefore left with 2xN points: N original points, and N projected points. We can now define to define the strength of the membership that associates each original point to the centroid of its cluster by means the following equation (see Table 9 for an application to the Etruscan towns data set):

$$(50) \quad \mu_{c_k}(p_i) = 1 - \frac{d_{i,c_k}}{\sum_{j}^{G} d_{i,c_j}}.$$

*where*:

$\mu_{c_k}(p_i) = $ membership of the i-th point with respect to the k-th cluster;

$d_{i,c_k} = $ distance bewteen the i-th point and the centroid of the k-cluster;

$d_{i,c_j} = $ distance bewteen the i-th point and the centroid of the j-cluster;

$G = $ Number of Clusters.

In the Etruscan towns case, the two vanishing points that we obtain cluster the twelve towns into the Western, Tuscan group and the Eastern group, that bundles together the Valdichiana and Latium towns. Interestingly, the vanishing point of the Tuscan group lies outside the convex hull, and in particular it sits in the Tyrrhenian sea. This is pretty coherent with the fact that the main source of wealth and power of the Tuscan group of cities, due to their marginal geographical positioning and the physical barrier of the Maremma swamp, lied in sea commerce and exchanges. Significantly, the town in this cluster with the highest membership is by far Populonia, that is, the port city of the cluster which mediated maritime trade to and from the other cities. More generally, strength of membership within the cluster decreases with the distance from the coast. The vanishing point of the second cluster lies close to Bolsena, that represents the natural distribution point between the remote Valdichiana cluster of small towns and the major Latium hub, at similar distances between the Vulci/Tarquinia and Veio/Cerveteri urban dyads. The vanishing point of this cluster is internal to the convex hull, as after the decadence of Veio and the ascent of Rome the hub of economic and social exchange was mostly determined by internal interactions and not, as for the other, by external ones. Here, the highest strength of membership is associated with closeness to the North-South axis passing through the vanishing point – a fact that seems to suggest that whereas the rationale of the first cluster is related to maritime trade, the rationale of the second is in the connection between the Valdichiana and Latium groups of towns, and not in the connection between the latter groups and the Tuscan block. As a consequence, even Tarquinia and Vulci, where as we know the strongest activation momentum is concentrated, present a relatively weak membership in the group in that their geographical location would be more relevant in terms of the mediation between the Eastern and Western clusters, whereas their location is marginal for the North-South connection. This seems to suggest that the momentum of Vulci-Tarquinia is not mainly generated by their function as a spatial hub for the Dodecapolis, but rather by their intrinsic drive that profits from the decadence of the larger, traditionally stronger nearby towns now under threat.

The application of the Meta-distance algorithm causes a significant repopulation of the original set of points. By means of the MDS, to each point of the original set corresponds a projected point through the Meta-distance end matrix $m_{i,j}^{[T]}$. Moreover, the more accurate the projection, the more representative is the position of the new points on the map in terms of the information embodied in the original data set of points, and the projected points belong to the same mathematical function to which the original dataset points belong (see equation 49 above). We also conjecture that the convergence matrix $m_{i,j}^{[T]}$ attains the maximum variance of its components, although a rigorous proof is not currently available:

(50) $\sigma\left(m_{i,j}^{[T]}\right) = \max\left(\sigma\left(m_{i,j}^{[t]}\right)\right); \quad t \in \{0,1,2,...,T\}.$

Figure 29 below synthesizes the steps to generate the convergence matrix $m_{i,j}^{[T]}$.

| Original Distance Matrix | 1_Volterra | 2_Perugia | 3_Cerveteri | 4_Tarquinia | 5_Populonia | 6_Chiusi | 7_Arezzo | 8_Vetulonia | 9_Bolsena | 10_Roselle | 11_Veio | 12_Vulci |
|---|---|---|---|---|---|---|---|---|---|---|---|---|
| 1_Volterra | 0 | 393.636292 | 473.195923 | 334.966339 | 140.575958 | 290.279663 | 257.412384 | 139.772125 | 342.457123 | 164.831848 | 520.021179 | 314.223022 |
| 2_Perugia | 393.636292 | 0 | 291.513641 | 236.440094 | 482.593292 | 115.514542 | 157.645874 | 365.482086 | 156.118576 | 326.580139 | 253.151123 | 259.287933 |
| 3_Cerveteri | 473.195923 | 291.513641 | 0 | 138.46077 | 478.858856 | 261.327911 | 375.486938 | 359.430542 | 166.631607 | 318.63678 | 97.290421 | 160.394806 |
| 4_Tarquinia | 334.966339 | 236.440094 | 138.46077 | 0 | 351.673706 | 157.434586 | 264.841675 | 227.640182 | 80.519928 | 184.817154 | 195.612442 | 32.058411 |
| 5_Populonia | 140.575958 | 482.593292 | 478.858856 | 351.673706 | 0 | 368.603302 | 371.63855 | 126.280243 | 389.029114 | 170.382172 | **547.284668** | 322.673248 |
| 6_Chiusi | 290.279663 | 115.514542 | 261.327911 | 157.434586 | 368.603302 | 0 | 114.159248 | 250.040619 | 94.803307 | 211.150116 | 264.114197 | 169.880035 |
| 7_Arezzo | 257.412384 | 157.645874 | 375.486938 | 264.841675 | 371.63855 | 114.159248 | 0 | 276.400848 | 208.929825 | 250.375656 | 372.702057 | 270.45871 |
| 8_Vetulonia | 139.772125 | 365.482086 | 359.430542 | 227.640182 | 126.280243 | 250.040619 | 276.400848 | 0 | 262.910675 | 44.107174 | 422.838806 | 200.173203 |
| 9_Bolsena | 342.457123 | 156.118576 | 166.631607 | 80.519928 | 389.029114 | 94.803307 | 208.929825 | 262.910675 | 0 | 218.970779 | 181.167389 | 105.965294 |
| 10_Roselle | 164.831848 | 326.580139 | 318.63678 | 184.817154 | 170.382172 | 211.150116 | 250.375656 | 44.107174 | 218.970779 | 0 | 379.495789 | 158.490082 |
| 11_Veio | 520.021179 | 253.151123 | 97.290421 | 195.612442 | **547.284668** | 264.114197 | 372.702057 | 422.838806 | 181.167389 | 379.495789 | 0 | 225.313705 |
| 12_Vulci | 314.223022 | 259.287933 | 160.394806 | 32.058411 | 322.673248 | 169.880035 | 270.45871 | 200.173203 | 105.965294 | 158.490082 | 225.313705 | 0 |

| Final Distance Matrix | 1_Volterra | 2_Perugia | 3_Cerveteri | 4_Tarquinia | 5_Populonia | 6_Chiusi | 7_Arezzo | 8_Vetulonia | 9_Bolsena | 10_Roselle | 11_Veio | 12_Vulci |
|---|---|---|---|---|---|---|---|---|---|---|---|---|
| 1_Volterra | 0 | 547.284607 | 547.284302 | 547.284607 | 0 | 547.284668 | 547.284302 | 0.000007 | 547.284607 | 0.064063 | 547.284302 | 547.284607 |
| 2_Perugia | 547.284607 | 0 | 0.001076 | 0.000005 | 547.284607 | 0.00001 | 0.001079 | 547.284607 | 0.000002 | 547.233154 | 0.001078 | 0.000005 |
| 3_Cerveteri | 547.284302 | 0.001076 | 0 | 0.001078 | 547.284302 | 0.001078 | 0.000006 | 547.284302 | 0.001075 | 547.232788 | 0.00001 | 0.001078 |
| 4_Tarquinia | 547.284607 | 0.000005 | 0.001078 | 0 | 547.284607 | 0.00001 | 0.00108 | 547.284607 | 0.000005 | 547.233154 | 0.001079 | 0 |
| 5_Populonia | 0 | 547.284607 | 547.284302 | 547.284607 | 0 | 547.284668 | 547.284302 | 0.000007 | 547.284607 | 0.064063 | 547.284302 | 547.284607 |
| 6_Chiusi | 547.284668 | 0.00001 | 0.001078 | 0.00001 | 547.284668 | 0 | 0.001081 | 547.284607 | 0.00001 | 547.233154 | 0.001079 | 0.00001 |
| 7_Arezzo | 547.284302 | 0.001079 | 0.000006 | 0.00108 | 547.284302 | 0.001081 | 0 | 547.284302 | 0.001078 | 547.232788 | 0.00001 | 0.00108 |
| 8_Vetulonia | 0.000007 | 547.284607 | 547.284302 | 547.284607 | 0.000007 | 547.284607 | 547.284302 | 0 | 547.284607 | 0.064063 | 547.284302 | 547.284607 |
| 9_Bolsena | 547.284607 | 0.000002 | 0.001075 | 0.000005 | 547.284607 | 0.00001 | 0.001078 | 547.284607 | 0 | 547.233154 | 0.001077 | 0.000005 |
| 10_Roselle | 0.064063 | 547.233154 | 547.232788 | 547.233154 | 0.064063 | 547.233154 | 547.232788 | 0.064061 | 547.233154 | 0 | 547.232788 | 547.233154 |
| 11_Veio | 547.284302 | 0.001078 | 0.00001 | 0.001079 | 547.284302 | 0.001079 | 0.00001 | 547.284302 | 0.001077 | 547.232788 | 0 | 0.001079 |
| 12_Vulci | 547.284607 | 0.000005 | 0.001078 | 0 | 547.284607 | 0.00001 | 0.00108 | 547.284607 | 0.000005 | 547.233154 | 0.001079 | 0 |

*Table 8: The Etruscan Dodecapolis. Above: The distance matrix of the original points (in red, the maximum value). Below: the final distance matrix, after 70 iterations of the Meta-distance process (in red and green, the two sets of distances).*

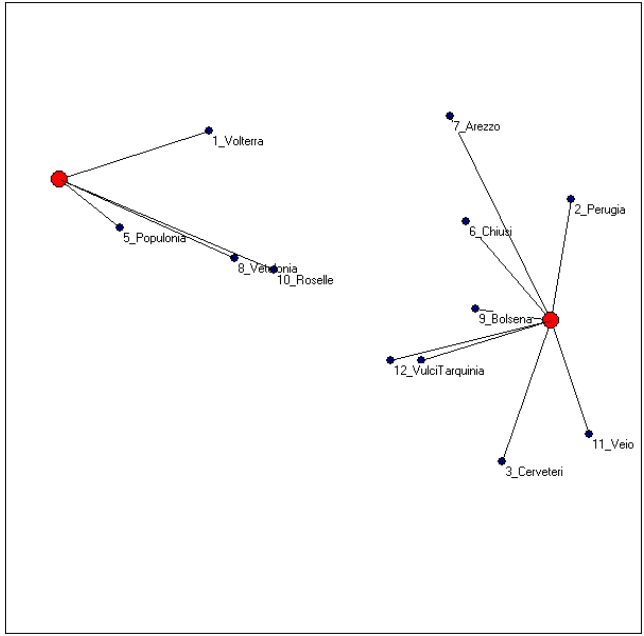

*Figure 28: The twelve Etruscan cities. The two vanishing points (red) projected from the final distance matrix after 70 iterations by means of MDS. The Meta-distance process clusters the data points (blue) into two sets, one of which lies outside the convex hull.*

| Cluster_#1 | | Cluster_#2 | |
| --- | --- | --- | --- |
| Name | Membership | Name | Membership |
| 1_Volterra | 0.7203 | 2_Perugia | 0.7938 |
| 5_Populonia | 0.8586 | 3_Cerveter | 0.7849 |
| 8_Vetulonia | 0.6353 | 4_Tarquini | 0.7461 |
| 10_Roselle | 0.5556 | 6_Chiusi | 0.7425 |
| | | 7_Arezzo | 0.6218 |
| | | 9_Bolsena | 0.8418 |
| | | 11_Veio | 0.84 |
| | | 12_Vulci | 0.693 |

*Table 9: Fuzzy clustering (see equation 50) applied to the Etruscan towns dataset according to the end-matrix of the Meta-distance process.*

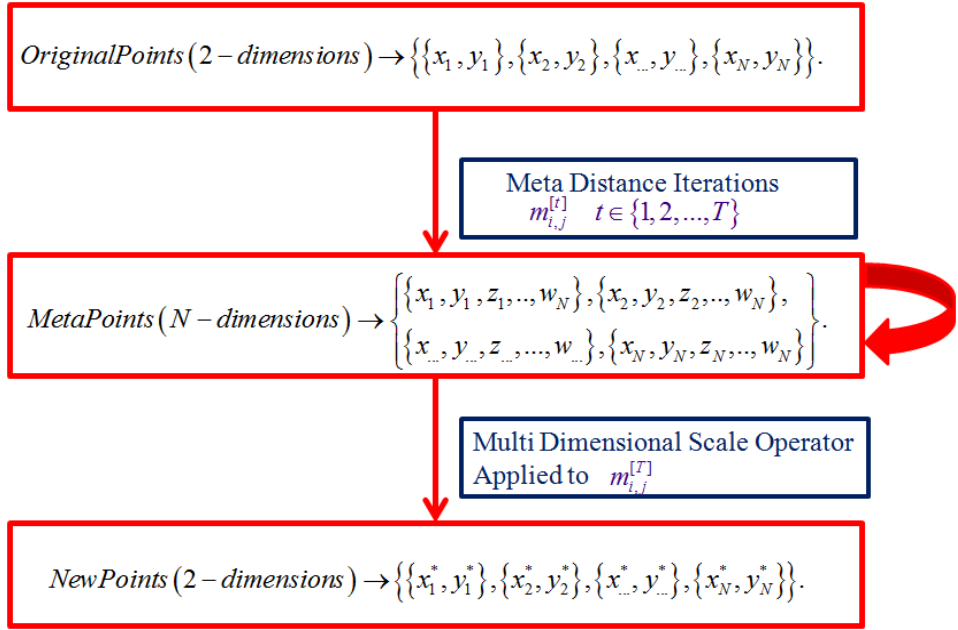

*Figure 29: Steps of the Meta-distance computation process.*

### 6.3 The optimal matrix of the Meta-distance algorithm

As the Meta-distance algorithm generates a distance matrix, $m_{i,j}^{[t]}$, at each iteration, MDS yields a different projection (image) of each point of the data set at each of them. If we aim at expanding the number of points of the original data set by means of a suitable projection and of its re-mapping, we need a criterion to choose the 'best' matrix in this sequence. In fact, the end matrix of the Meta-distance process need not be the one that best preserves and highlights the information content of the original data set. Let us denote by $m_{i,j}^*$ the optimal matrix from the Meta-distance process, determined with respect to the suitable cost function. Once the optimal $m_{i,j}^*$ matrix has been selected, we apply to it the MDS operator, to produce a new set of points in addition to the original one. The cost function $J$ we choose to determine the optimal matrix is given by equations (52-57) below. As we are working on a Meta-distance, the cost criterion that is chosen is the sum of the matrix entries:

$$(52) \quad S^{[t]} = \sum_{i=1}^{N-1} \sum_{j=i+1}^{N} m_{i,j}^{[t]}; \quad t > 1;$$

$$(53) \quad \Delta S^{[t]} = \sum_{i=1}^{N-1} \sum_{j=i+1}^{N} \left( m_{i,j}^{[t]} - m_{i,j}^{[t-1]} \right);$$

$$(54) \quad \begin{cases} \delta = 1 & \Delta S^{[t]} > 0; \\ \delta = -1 & \Delta S^{[t]} < 0; \end{cases}$$

$$(55) \quad J^{[t]} = \log \left( \frac{S^{[t]}}{\left| \Delta S^{[t]} \right|} \right) \cdot \delta;$$

$$(56) \quad J^{[t*]} = \underset{t}{Max} \left\{ J^{[t]} \right\};$$

$$(57) \quad m_{i,j}^{[t*]} = m_{i,j}^{[t]} \ .$$

Simulation analysis reveals that the optimal distance matrix $m_{i,j}^{[t*]}$ often coincides with the end matrix $m_{i,j}^{[T]}$, but not necessarily. When the matrix $m_{i,j}^{[t]*}$ differs from $m_{i,j}^{[T]}$, the projected vanishing points may be more than two, thus defining a more complex spatial clustering pattern. Once the optimal matrix is selected, the MDS will project $m_{i,j}^{[t*]}$ into a two-dimensional space, where each original point will be coupled with its projection. Therefore, at the end of this process, the N points of the original data set will be augmented with N new points.

*6.4 Meta-distance quantities: Summing up*

To sum up, the Meta-distance process defines the following quantities, which prove to be particularly useful for an advanced application of the topological approach, as we will see in the following section:

a. A set of new points, generated by the optimization of the J function (equations 52-57), that doubles the number of points of the original data set. The new data set will be composed of 2 x N points (see Figure 29). It is of course an open point at this stage whether this augmented data set contains more information than the original one. We call this new data set *Meta-distance data set*.
b. A set of centroids of the clusters of the new projected points (see equation 50 and Figure 25). In this case, we need to assess whether such new meta-points carry important information about the original data points distribution. We call such centroids *Meta-cluster Points* or *Vanishing Points*.
c. A *scalar field* generated by the distance of each generic point of the two-dimensional space from each projected point of the Meta-distance algorithm, weighted by the distance that each projected point has from its original point. The farther away a projected point from its original point, the stronger the weight:

*Legend* :

$N =$ Number of Source and Projected Points;

$k \in [1, 2, ..., N]$;

$k =$ index for Source and Projected Points;

$Ps_{k_{(x,y)}} =$ Longitude and Latitude of the $k$-th SourcePoint;

$MaxD =$ Maximum of Euclidean distance among the Source Points;

$Pp_{k_{(x,y)}} =$ Longitude and Latitude of the $k$-th Projected Point;

$Md_k =$ Distance between the k-th source point and the k-th its projection;

$M =$ Number of the points of the entire plane;

$i =$ index for each point of the plane;

$i \in [1, 2, ..., M]$;

$Pg_{i_{(x,y)}} =$ Longitude and Latitude of the i-th point of the plane;

$Dw_i =$ Weighted Summation of the Distances of the i-th point of the plane from each of the projected points;

$Act_i =$ Activation of each point of the entire plane;

$$(58) \quad Md_k = \sqrt[2]{\left(Ps_{k_{(x)}} - Pp_{k_{(x)}}\right)^2 + \left(Ps_{k_{(y)}} - Pp_{k_{(y)}}\right)^2};$$

$$(59) \quad Dw_i = \sum_{k=1}^{N}\left(1 - \frac{Md_k}{\sum_{j=1}^{N} Md_j}\right) \cdot \frac{\sqrt[2]{\left(Pp_{k_{(x)}} - Pg_{i_{(x)}}\right)^2 + \left(Pp_{k_{(y)}} - Pg_{i_{(y)}}\right)^2}}{MaxD};$$

$$(60) \quad A_i = e^{-Dw_i}.$$

We call this new object *Meta-cluster Field*.

## 7. Applying an advanced topological approach: Meta-distance data sets

Is the additional set of data generated by the Meta-distance process an useful addition to the original data set in terms of information content? To provide a first assessment of this hypothesis, we have tested both the original data sets and the Meta-distance data sets generated for different cases of known epidemics for which the benchmarks can be clearly defined and performance unambiguously evaluated.

In our opinion Meta Distance data set should be richer than the original data set. In order to assess this hypothesis we have tested the original data set and the Meta Distance data set using different examples of known epidemics [21].

*7.1. Chikungunya Fever, Italy 2007.*

Chikungunya fever is a toga viral illness spread by the bite of infected Aedes Aegypti mosquitoes. The disease is debilitating but rarely life-threatening. The virus was first isolated between 1952-1953 from both man and mosquitoes during an epidemic of fever that was considered clinically indistinguishable from dengue, in Tanzania. Periodic outbreaks have taken place in Asia, Africa and the latest in Indonesia in 1999. Up to 2007, no autochthonous cases took place outside these areas, but between July and August 2007, 205 cases of Chikungunya fever were recorded in the environs of the two villages of Castiglione di Cervia and Castiglione di Ravenna in Northern Italy. The spatial topography of the epidemic had a decreasing concentric gradient, with fewer cases taking place moving away from the epicenter. The probable index case was identified as a traveler from an area of India where an epidemic of Chikungunya fever was underway, and the highest concentration of cases was identified as the village of Castiglione di Cervia [32].

Figure 30 reports the 14 observations of the epidemics we have analyzed. Figure 31 shows the estimation of the outbreak using TWC($\alpha^*$) by means two different data sets: The original data set of points mapping the epidemics in its late phase of development, and the new augmented data set generated by the application of the Meta-distance algorithm via the optimization of the *J* cost function. TWC ($\alpha^*$) estimates the same correct outbreak in both cases, but the Alpha Vector is quite different.

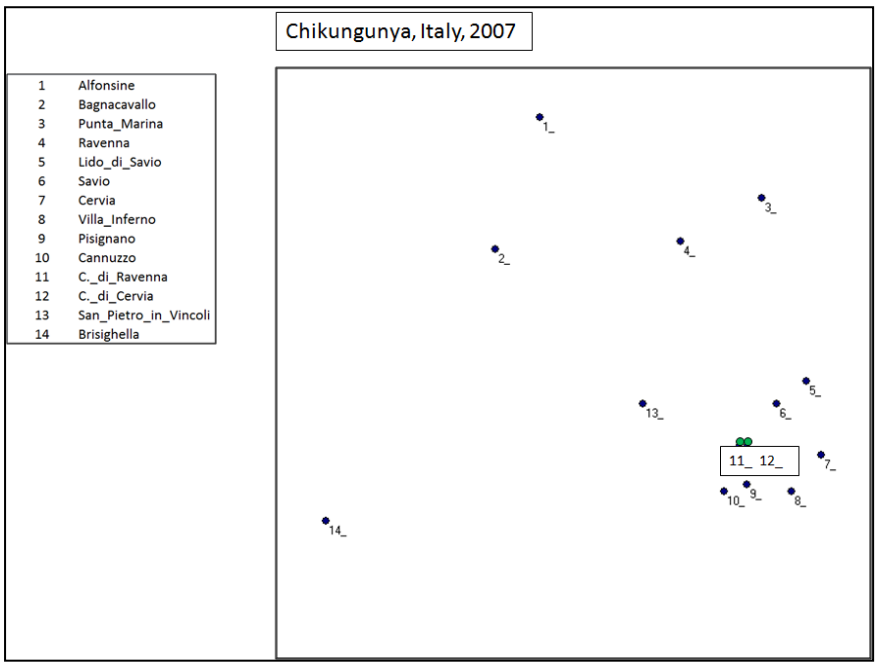

*Figure 30: The 14 points of the Chikungunya data set (in green color, the known outbreak).*

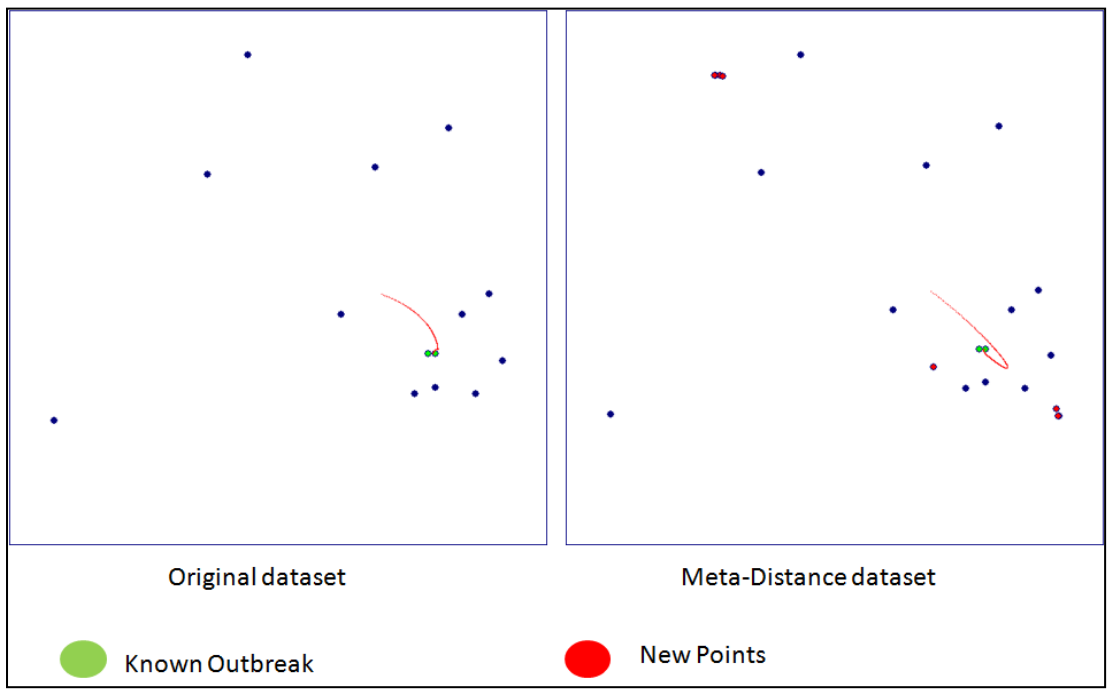

*Figure 31: Chikungunya fever: The two hypotheses for TWC(α$^*$) with the original data set (left) and with the augmented data set (right).*

### 7.2 Foot and Mouth Disease, UK 1968-69

The animal epidemic caused by the Foot and Mouth disease (FMD) of 1968-69 was responsible for the death of 2,000 animals and the compulsory slaughter of half a million more. It wrought economic chaos to the UK cattle and meat producing industries. The origin of the epidemic was traced to infected pig swill at Bryn Farm, near Oswestry in the English county of Shropshire [33]. FMD was not a new pathology in the UK in 1967–68. Throughout the 1950s and 1960s, it was not unusual to witness outbreaks of FMD, sometimes even 2-3 in a same year. From 1954 to 1967, apart from the 1967–68 outbreak, there were 1002 outbreaks, with an average of 75 cases every year, and there were only two years – 1963 and 1964 – when no cases were observed, the longest period without the disease since 1908. In the 1967-68 outbreak, when the first case was diagnosed at Bryn Farm on October 25, a Wednesday, the normal Oswestry market was taking place and two cows from that farm had gone to the market that morning. After the disease was recognized in the infected cattle, the State Veterinary Service immediately decided the closing down of the market. The two cows were examined on the next day, October 26 and found free of the disease, but despite this they were included in the slaughter. Some animals traded at the market on that day had left, and gone as far as Banffshire, Scotland, and Devon in the South-West before precautionary measures were enacted. All such animals were traced and found healthy, so it was decided not to prescribe compulsory slaughter of all the animals in the market. On Monday, October 30, the situation changed dramatically, as 9 new cases were confirmed, six of which close to the original outbreak location. The other 3 were at considerably distant locations (12, 35 and 100 miles, the furthest in Lancashire). From October 30 onwards, there was a dramatic escalation of contagion, which in the end involved about 2,000 animals.

For our test, we have used the coordinates of the points corresponding to the 20 farms involved in the first week of the epidemics. Figure 32 shows how the TWC($\alpha^*$) finds the same correct outbreak in both cases, with slightly different Alpha Vectors.

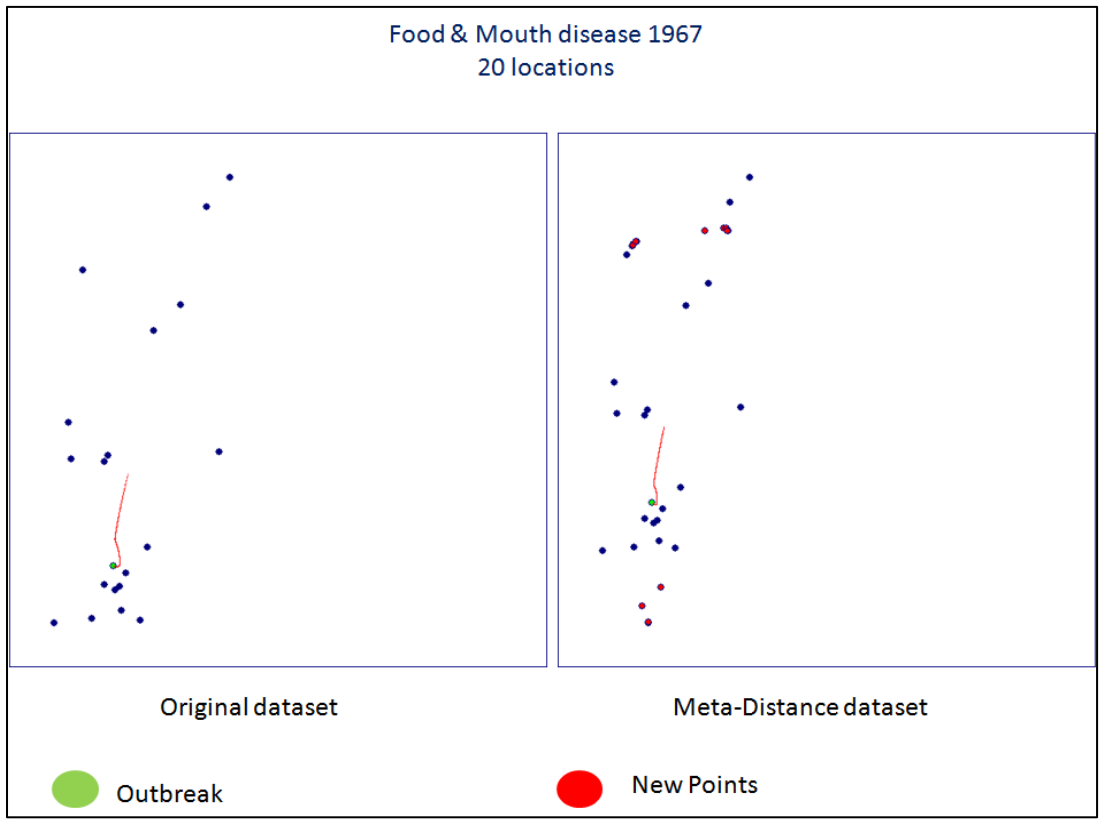

*Figure 32: Foot and Mouth disease: The two hypotheses for TWC($\alpha^*$) with the original data set (left) and with the augmented data set (right).*

### 7.3 Dengue fever, Brazil 2001

The Dengue fever epidemics in the State of São Paulo, Brazil, in 2001 was monitored for seven weeks as reported in [34-35]. The map of the epidemics diffusion shows the actual pattern of diffusion of contagion over 7 months in the 54 counties where at least one case was reported. At the end of the monitoring time, a map of fever propagation was built. From these data, a digital map has been elaborated, which reports the main town center of all the 54 counties that have witnessed at least one case of the epidemics. This map is our original database, which we use to analyze the actual propagation dynamics, according to the TWC($\alpha^*$) for the two data sets: The original one with 54 points, and the augmented one with 108 points. This time, unlike the previous cases, as shown in Figure 33, the TWC($\alpha^*$) using the original data fails to find the origin of the outbreak, whereas the TWC($\alpha^*$) using the augmented data set finds the correct outbreak (Rio Prieto). Also the corresponding Alpha Vectors are consequently very different.

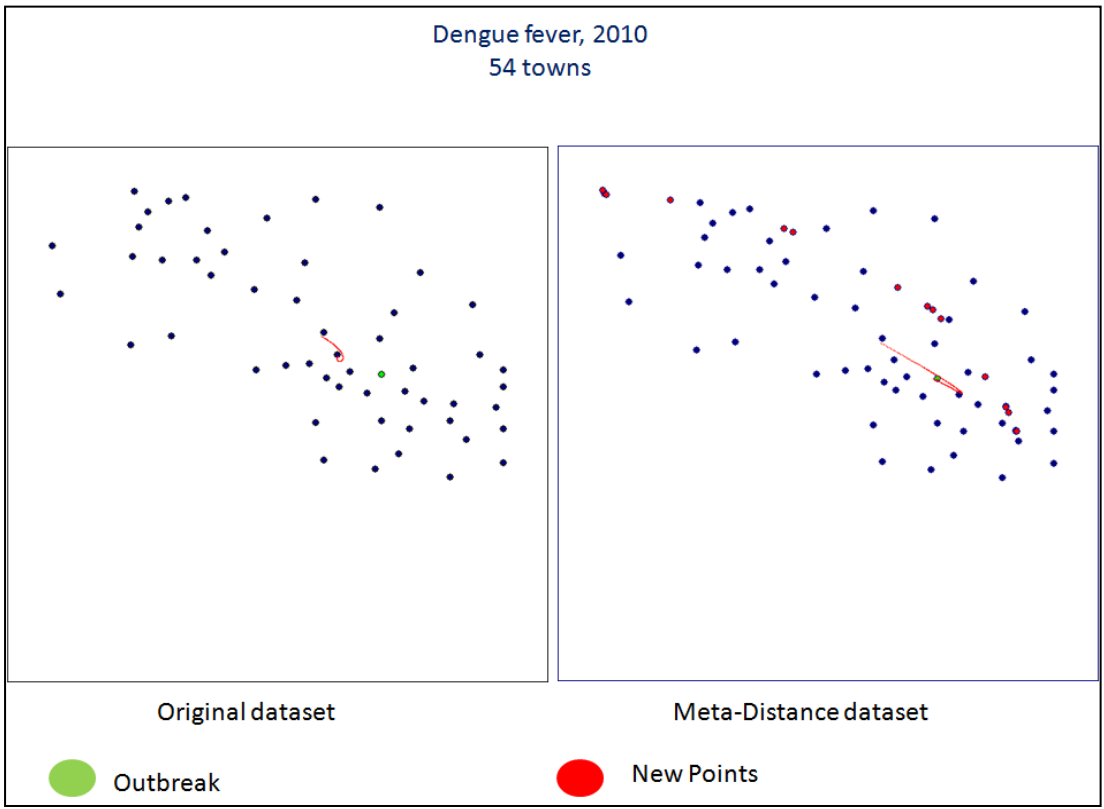

*Figure 33: Dengue Fever: The two hypotheses for TWC($\alpha^*$) with the original data set (left) and with the augmented data set (right).*

### 7.4 Food epidemic, Oahu 2010

Data from a food epidemic happened in Oahu, Hawaii (1,245 cases) were collected systematically for 12 months (for more details see [21]). We have received this dataset from Al Bronstein, director of the Rocky Mountain Poison Center at the time. Figure 34 shows the geographical distribution of cases of the epidemic at the end of 2010, and Table 10 shows the distribution of the new cases for each month of 2010.

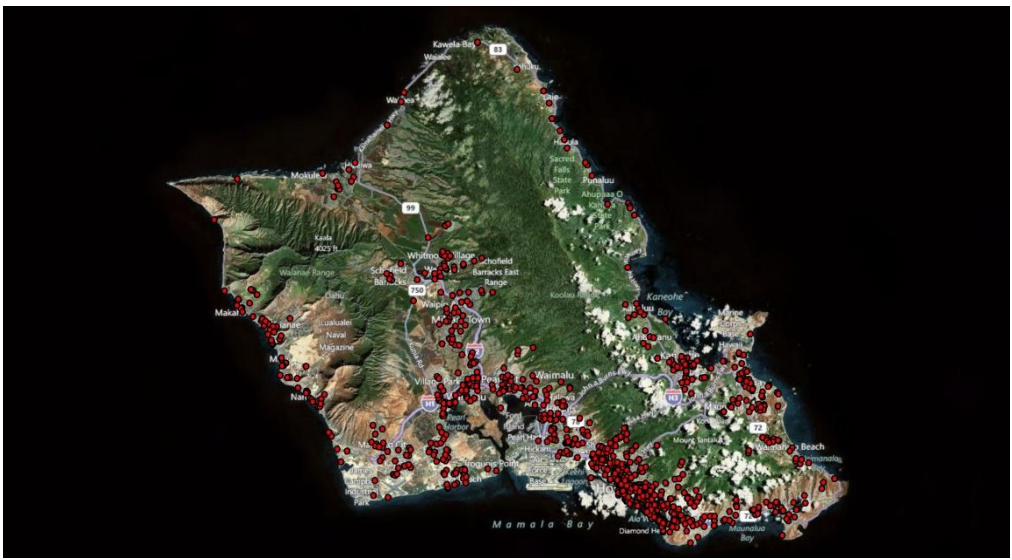

*Figure 34: Global distribution of the Oahu, Hawaii food epidemic in 2010 (Red circles).*

| Oahu 2010: Number of Cases Each Month | |
|---|---|
| Jan | 108 |
| Feb | 109 |
| March | 114 |
| April | 98 |
| May | 79 |
| June | 79 |
| July | 93 |
| August | 109 |
| Sept | 92 |
| Oct | 134 |
| Nov | 114 |
| Dec | 108 |

*Table 10: Dynamics of the Oahu food epidemic by month, 2010.*

In Figure 35, we apply the Meta-distance approach to the cases of January 2010, and find the two vanishing points and the global Meta-cluster field for the January data only.

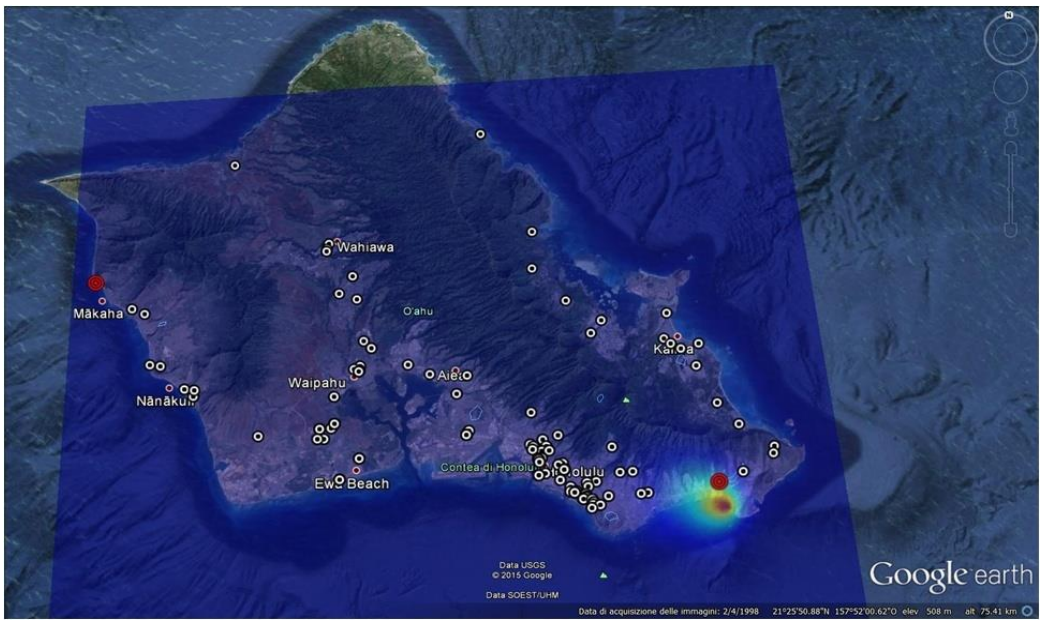

*Figure 35: The 108 cases during January 2010 (white circles), and the two Meta-clusters found by the Meta-distance algorithm (Red Spirals), and the Meta-cluster field (Color Map Overlays).*

From Figure 35, we see that the two Vanishing Points are located in areas with no or few cases. If we compare this map with the map reporting all the cases occurred up to the end of 2010, we easily check how the two Meta-clusters mark two areas characterized by the emergence of a high concentration of new cases after January 2010 (see Figures 36a-b).

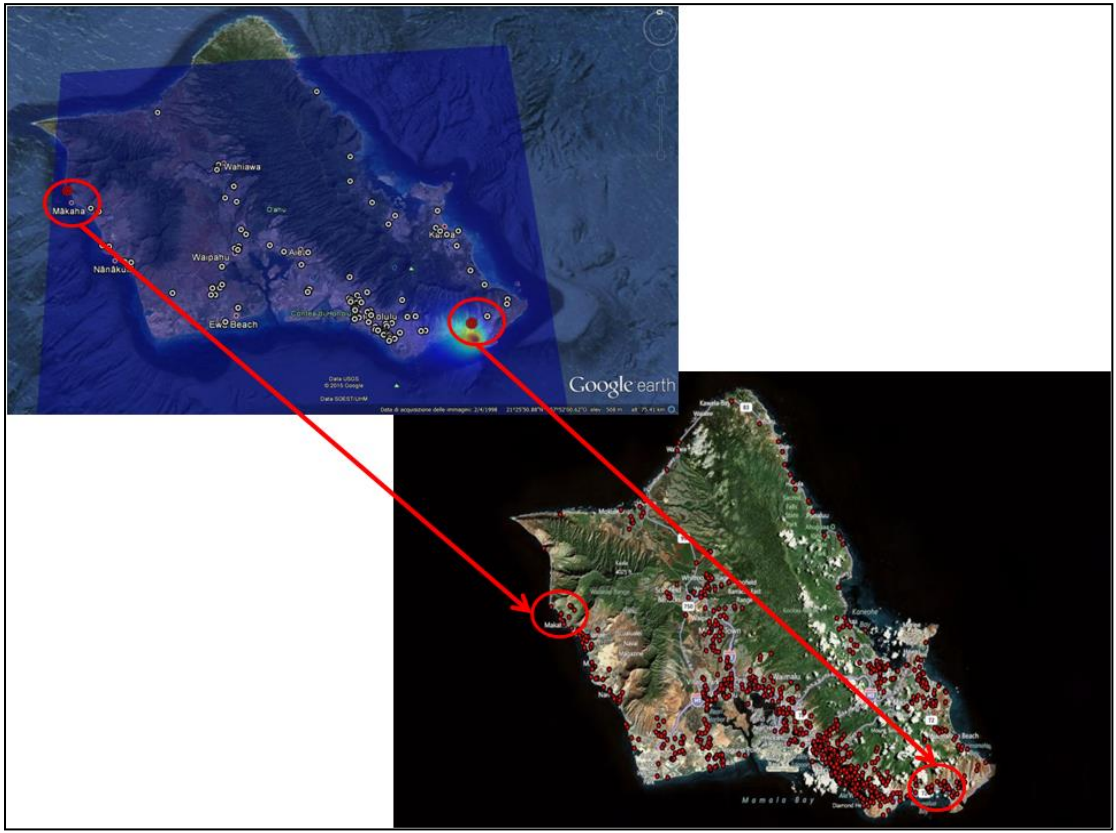

*Figure 36a: Up: New cases in January 2010 (white points) with the two Meta-clusters (red points). Down: the total cases at the end of 2010.*

This application clearly shows the potential of Vanishing Points as predictive tools for an ongoing epidemic phenomenon. Together with the examples in the previous subsections, we have therefore a first corroboration that the meta-distance approach can find more accurately than the basic topological approach the outbreak point of an epidemics in particularly difficult cases, while at the same time providing analogous results to the basic approach in simpler cases (although with different reconstructions of the underlying dynamics). Moreover, the Meta-distance approach may also function as a predictive tool for the future evolution of the epidemics also in regions outside the convex hull of observed cases, thus extending its applicability beyond the feasibility range of the basic topological approach.

These preliminary tests only provide an initial corroboration, and much more research and testing is needed to establish whether the Meta-distance approach may be considered as a robust, reliable tool for the analysis of a vast variety of epidemic and pseudo-epidemic processes. However, these first results are encouraging enough to warrant further investigation.

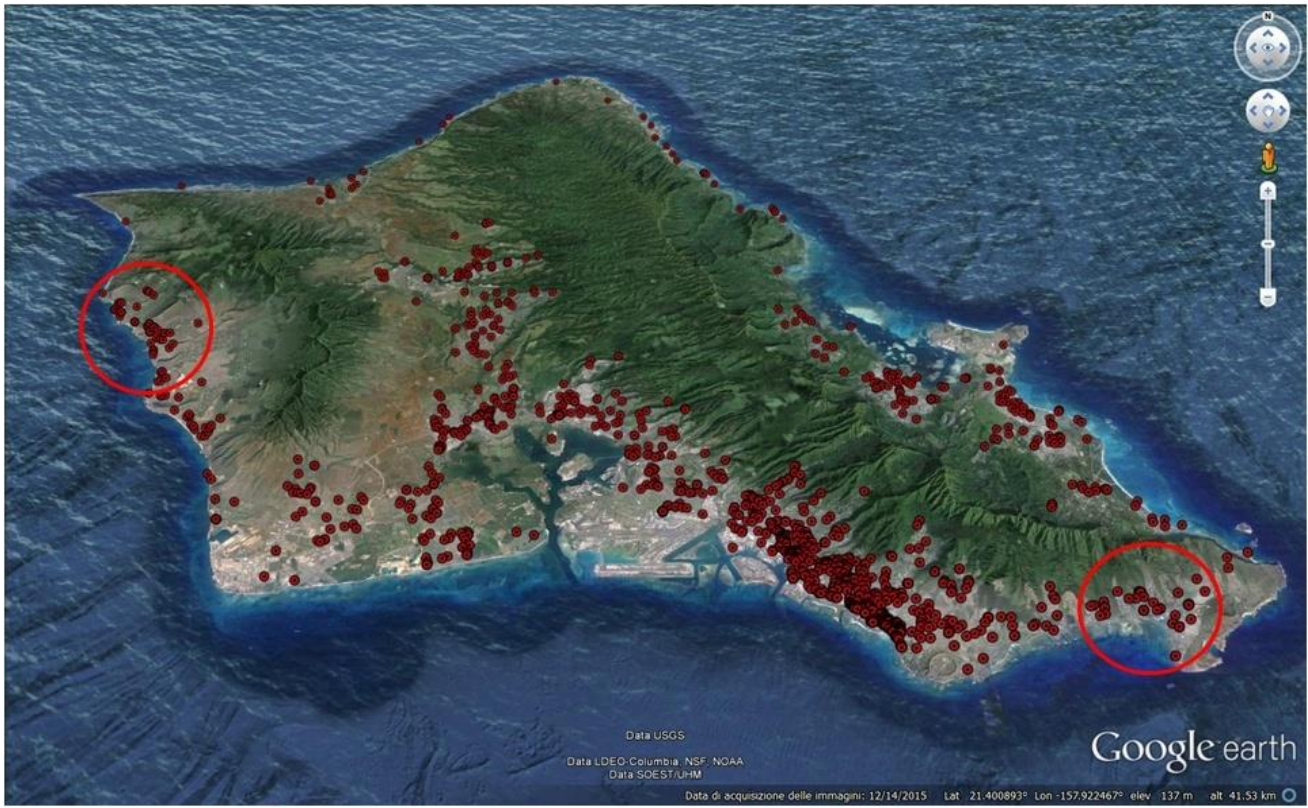

*Figure 36b: The full distribution of cases at the end of 2010, with the new contagion areas predicted in January by the location of the two Meta-clusters.*

## 8. Conclusions

In this paper, we have presented for the first time in a systematic way the topological approach to the spatial analysis of epidemic and pseudo-epidemic processes. Our approach presents a significant innovation with respect to existing ones, but due to this it has to be considered still a work in progress to a large extent. Here, we have presented a still fragmentary ensemble of formal definitions and demonstrations, conjectures and experimental tests. However, we feel that this first step has enabled us to make some interesting progress with respect to the literature standards, so as to encourage further work in this vein in the future. Moreover, the range of applicability of these techniques seems to extend beyond epidemic processes to a diverse range of pseudo-epidemic ones, where it allows to develop fresh insights with respect to the standard methods [36-38].

On the other hand, the topological approach presented here establishes some new tenets for the spatial analysis of epidemic and pseudo-epidemic processes which can direct future research toward new, promising lines:

a. Discrete manifestations of an ongoing epidemic (pseudo-epidemic) process contain spatially encoded, key information about its past and future unfolding. In other words, in the GIS positions of events from a same process, time seems to be squashed into space configurations, that can therefore be read in terms of an intrinsic space grammar that abstract from the actual sequence and frequency of observations.

b. Such information may be retrieved through a bottom-up analysis of the spatial configurations of the observed points, as available in a spatial snapshot. The retrieving algorithms are rigorously data driven.
c. The concept of pseudo-distance is a new concept with interesting features. Pseudo-distance is a quasi-metric (due to symmetry violation) and it is also a semi metric (due to triangular inequality violation) distance. Pseudo-distance is a completely relational distance: That is, the distance of any point from another depends on its distance from all the other points of the spatial distribution. Under the action of the pseudo-distance, space is conveniently warped and the classic Euclidean distances among points may be expanded or contracted according to the context. In general, groups of neighboring points move apart, while relatively distant points move closer, bringing together those points which have structural features in common in the context of the spatial distribution of points and distancing those which are diverse among each other.
d. From our analysis, a new concept of 'information mass of a point' also emerges. In any distribution of spatially referenced points, each point exerts a force of attraction on others, that is proportional to its cumulative similarities with the other points of the distribution. The idea that an abstract point may have an information mass may be important also in the theory of physical networks, to determine how each node of an assigned graph warps the space around it. In other words, the concept of information mass transforms the flat plane of the Euclidean distances into a space with local warping.
e. The Meta-distance is also a new concept that emerges from this paper. We have shown how the new points generated from a recursive application of the Meta-distance algorithm seem to provide a meaningful expansion of the original data set of points. Moreover, the Meta-distance algorithm clusters the dataset in a fuzzy way that enlightens subtle structural properties which seem to have a predictive value on the future spatial evolution of the process itself.

The topological approach gives therefore a new, more profound meaning to the geo-localization of data and opens new ways of analysis and interpretation of geo-localized data. Such analysis, as shown in our reference case study of the spatial distribution of ancient Etruscan towns, may also allow to re-interpret already existing knowledge or to integrate incomplete or fragmentary knowledge of a spatial (pseudo-epidemic) process in useful ways. The domain of applicability of the topological approach therefore potentially extends to a vast range of disciplines in the natural and in the social sciences. Testing it in different contexts and conditions represents therefore a very stimulating challenge, that we look forward to pursuing in future research.